\documentclass[11pt]{article}
\usepackage{amsmath}
\usepackage{amssymb}
\usepackage{epsfig,graphics,color,graphicx,graphpap}

\usepackage{theorem}
\numberwithin{equation}{section}
\numberwithin{equation}{subsection}

\newcommand{\RR} {\mathbb R}
\newcommand{\CC} {\mathbb C}
\newcommand{\ZZ} {\mathbb Z}

\newcommand{\pa} {\partial}
\newcommand{\Cal} {\mathcal}

\newcommand\dv {\operatorname{div}}
\newcommand\ep {\varepsilon}
\newcommand\rot {\operatorname{rot}}
\newcommand\Ker {\operatorname{Ker}}
\newcommand\bmo {\operatorname{bmo}}
\newcommand\fh {\mathfrak{h}}
\newcommand\wtp {\widetilde{\Pi}}
\newcommand\Spec {\operatorname{Spec}}
\newcommand\Span {\operatorname{Span}}

\renewcommand\Im {\operatorname{Im}}
\renewcommand\Re {\operatorname{Re}}

\newcommand\tu {\tilde{u}}
\newcommand\tw {\tilde{w}}
\newcommand\sk {$\text{}$ \newline}
\newcommand\nsk {\noindent}

\newcommand{\beq} {\begin{equation}}
\newcommand{\eeq} {\end{equation}}
\newcommand{\demo} {\noindent {\it Proof. }}
\newcommand{\qed} {\hfill$\Box$ {\newline $\text{}$}}

\theoremstyle{plain}
\newtheorem{theorem}{Theorem}[subsection]
\newtheorem{proposition}[theorem]{Proposition}
\newtheorem{corollary}[theorem]{Corollary}
\newtheorem{lemma}[theorem]{Lemma}

\title{Euler Equation on a Rotating Surface}

\author{Michael Taylor 
\footnote{2010 Math Subject Classification: 35Q31, 35R09
{\it Key words.} Euler equation, Coriolis force, vorticity, stability.
Work supported by NSF grants DMS-1500817, DMS-1312874}}

\date{}

\begin{document}

\maketitle

\centerline{\bf With an Appendix by Jeremy Marzuola and Michael Taylor}

\begin{abstract}
We study 2D Euler equations on a rotating surface, subject to the
effect of the Coriolis force, with an emphasis on surfaces of revolution.
We bring in conservation laws that yield long time estimates on solutions
to the Euler equation, and examine ways in which the solutions behave
like zonal fields, building on work of B.~Cheng and A.~Mahalov, examining how
such 2D Euler equations can account for the observed band structure of rapidly
rotating planets.  Specific results include both an analysis of time averages
of solutions and a study of stability of stationary zonal fields.  The
latter study includes both analytical and numerical work.
\end{abstract}

$$\text{}$$
{\bf Contents}

\nsk
1. Introduction

1.1. Further properties of the operator $B$

\nsk
2. Basic existence results

2.1. Short time existence

2.2. Vorticity equation

2.3. Long time existence

\nsk
3. Bodies with rotational symmetry

3.1. Stationary solutions

3.2. Time averages of solutions

3.3. Another conservation law

3.4. Computation of $\chi$ and $\xi$

3.5. Smoothness issues

\nsk
4. Stability of stationary solutions

4.1. Arnold-type stability results

4.2. Linearization about a stationary solution

4.3. Further results on linearization

\nsk
5. Appendix, by J.~Marzuola and M.~Taylor: Matrix approach

and numerical study of linear instability

5.1. Matrix analysis for $f(x)=cP_2(x_3)$

5.2. Matrix analysis for $f(x)=cP_3(x_3)$

5.3. Numerical study of truncated matrices

\section{Introduction}\label{c1}

Let $M=\pa\Cal{O}$ be a surface in $\RR^3$, rotating about the $x_3$-axis at
constant angular velocity $\omega=-\Omega/2$.  A natural class of such bodies
would be those that are rotationally symmetric about the $x_3$-axis, and we
will eventually settle into the study of that class, but initially we will not
make that assumption.  We will assume $M$ is diffeomorphic to the standard
unit sphere $S^2$.
We aim to study 2D incompressible Euler flows on $M$.

The approach of Rossby to the effect of the Coriolis force on flows on $M$,
described on p.~21 of \cite{Ped}, yields the Euler equation
\beq
\frac{\pa u}{\pa t}+\nabla_u u=\Omega \chi(x)Ju-\nabla p,\quad \dv u=0,
\label{1.1}
\eeq
where
\beq
\chi(x)=e_3\cdot\nu(x),
\label{1.2}
\eeq
$\nu(x)$ being the unit outward pointing normal to $M$ at $x$.  Here $u$ is the
flow velocity, a tangent vector field to $M$, and $J:T_xM\rightarrow T_xM$
is counterclockwise rotation by $90^\circ$.  In case $M=S^2$, we have $\chi(x)=
x_3$.  For more general $M$ that are rotationally symmetric about the $x_3$-axis
and that have positive Gauss curvature, we have $\chi(x)=\chi(x_3)$.

Bringing in the 1-form $\tu$, arising from $u$ via the isomorphism $T_xM
\approx T^*_xM$ determined by the metric tensor on $M$, we can rewrite
\eqref{1.1} as
\beq
\frac{\pa\tu}{\pa t}+\nabla_u\tu=\Omega\chi*\tu-dp,\quad \delta\tu=0.
\label{1.3}
\eeq
We can eliminate $p$ from \eqref{1.3} via the Leray projection $P$, the
orthogonal projection of $L^2(M,\Lambda^1)$ onto the subspace where $\delta
\tu=0$.  We get
\beq
\frac{\pa\tu}{\pa t}+P\nabla_u\tu=\Omega B\tu,\quad \tu=P\tu,
\label{1.4}
\eeq
where
\beq
B\tu=P(\chi*P\tu).
\label{1.5}
\eeq

We mention a few essential properties of $B$, which will facilitate the
analysis of \eqref{1.4}.  First,
one can deduce from the Hodge decomposition that, on 1-forms,
\beq
P=-\delta \Delta_2^{-1}d,
\label{1.6}
\eeq
where $\Delta_2^{-1}$ denotes the inverse of the Hodge Laplacian on 2-forms,
defined to annihilate the area form.  Thus
\beq
\aligned
B\tu&=-\delta\Delta_2^{-1}d(\chi*P\tu) \\
&=-\delta\Delta_2^{-1}(d\chi\wedge*P\tu),
\endaligned
\label{1.7}
\eeq
since $d*P\tu=0$.  We deduce that
\beq
\text{$B$ is a compact, skew-adjoint operator on } L^2(M,\Lambda^1).
\label{1.8}
\eeq
In fact, $B\in OPS^{-1}(M)$, i.e., $B$ is a pseudodifferential operator of
order $-1$.  The skew adjointness is a direct consequence of the formula
\eqref{1.5}, the skew adjointness of the Hodge star operator $*$, and the
commutativity of $*$ and multiplication by $\chi$.
Further results on $B$ can be found in \S{\ref{c1s1}}.

Our interest in the equation \eqref{1.1}
was stimulated by the recent paper \cite{CM}
of B.~Cheng and A.~Mahalov, investigating the case where $M$ is the standard
sphere $S^2$.  That paper took \eqref{1.1} as a model of the behavior of the
atmosphere of a rotating planet, and investigated how it might account for
observed band structure, particularly on rapidly rotating planets, such as
Jupiter.  This involved a study of zonal flows, i.e., velocity fields of the
form $J\nabla f$, where $f=f(x_3)$ is a zonal function.  The paper looks at
time averages,
\beq
\Cal{A}_{S,T}u=\frac{1}{T}\int_S^{S+T} u(t)\, dt,
\label{1.9}
\eeq
for a solution $u$ to \eqref{1.1}.  The main result (Theorem 1.1 of \cite{CM})
is that, if $u_0\in H^k(S^2)$,
$k\ge 3,\ \dv u_0=0$, there exists $T_0>0$, independent of $\Omega$, such that
\eqref{1.1} has a unique solution for $t\in [0,T_0/\|u_0\|_{H^k}]$, satisfying
$u(0)=u_0$, and, for $0\le S< S+T\le T_0/\|u_0\|_{H^k}$, one has
\beq
\|(I-\Pi)\Cal{A}_{S,T}u\|_{H^{k-3}}=O(|\Omega|^{-1}),
\label{1.10}
\eeq
where $\Pi$ is a projection of the space of divergence-free velocity fields
on $S^2$ onto the space of zonal fields.

$\text{}$

In the present paper, we push the study of
\eqref{1.1} in the following directions.

\sk
(A) Investigate a larger class of rotating bodies, beyond the class of
rotating spheres.

\sk
(B) Establish estimates on $\Cal{A}_{S,T}u$ that are uniform in $S,T\in
(0,\infty)$, without restrictions on their size.

\sk
(C) Investigate another way that large $|\Omega|$ enhances band formation,
namely by enhancing the stability of zonal fields as stationary solutions to
\eqref{1.1}.

$\text{}$

These are natural directions to pursue.  Rapidly rotating planets are flattened
at the poles and bulging at the equator.  Also, a planet like Jupiter has
been rotating for a very long time.  Of course, of major interest to us is
the set of interesting new mathematical challenges that arise in addressing
these issues.

We proceed as follows.  In \S{\ref{c2}} we produce basic existence results,
starting with short time existence in \S{\ref{c2s1}}.
Results of \S{\ref{c2}} apply to any surface
$M$ diffeomorphic to $S^2$, with no symmetry hypothesis on the geometry.
To go from short time to long time existence, we derive in \S{\ref{c2s2}} an
equation for the vorticity $w=\rot u$, namely
\beq
\frac{\pa}{\pa t}(w-\Omega\chi)+\nabla_u(w-\Omega\chi)=0.
\label{1.11}
\eeq
This is a conservation law, yielding a uniform bound on $\|w(t)\|_{L^\infty}$
on any time interval on which \eqref{1.1} has a sufficiently smooth solution.
Using this, we adapt the classical Beale-Kato-Majda argument \cite{BKM} to
establish existence for all $t$ of a solution to \eqref{1.1}, provided $u(0)=u_0$
is divergence-free and belongs to $H^s(M)$ for some $s>2$.  This is accompanied
by the estimate
\beq
\|u(t)\|^2_{H^s}\le C\|u(0)\|^2_{H^s}\,\exp \, \exp
\Bigl(C_s(\|w(0)\|_{L^\infty}+C|\Omega|)t\Bigr).
\label{1.12}
\eeq

In \S{\ref{c3}} we specialize to the class of smooth compact surfaces $M\subset\RR^3$
that are invariant under the group of rotations about the $x_3$-axis, and that
in addition have positive Gauss curvature everywhere.  This hypothesis will
be in effect for the rest of the paper.
As already mentioned, this symmetry hypothesis implies $\chi=\chi(x_3)$ in
\eqref{1.1}.  In \S{\ref{c3s1}} we show that when $f$ is a zonal function, the associated
zonal vector field $u=J\nabla f$ is a stationary solution to \eqref{1.1}.  We
also give examples of stationary solutions that are not zonal fields.
In \S{\ref{c3s2}}, we study time averages of the form \eqref{1.9} and establish
estimates of the form
\beq
\aligned
&\|(I-\Pi)\Cal{A}_{S,T}u\|_{H^{-3,q(\theta)}} \\
&\le \frac{C_\theta}{|\Omega|}\Bigl\{T^{-1}\|u(S+T)-u(S)\|_{L^2}
+C\|u(0)\|^{2-\theta}\bigl(\|w(0)\|_{L^\infty}+2|\Omega|\bigr)^\theta\Bigr\},
\endaligned
\label{1.13}
\eeq
for $0<1<\theta$, with $q(\theta)=1/(1-\theta)$.  This is valid for all
$S,T\in (0,\infty)$.  It should be expected that the norm on the left side of
\eqref{1.13} is weaker than that in \eqref{1.10}.  In any case, having a weak norm seems
consistent with the appearance of complicated eddies within the bands of a
planet like Jupiter.

We proceed in \S{\ref{c3s3}} to derive an additional conservation
law, of the form
\beq
\frac{\pa}{\pa t} \int\limits_M \xi(x)w(t,x)\, dS(x)=0,
\label{1.14}
\eeq
for solutions to \eqref{1.1} on our radially symmetric domain.  We discuss
computations of $\chi$ and $\xi$ in \S{\ref{c3s4}} and technical smoothness
results in \S{\ref{c3s5}}.

In \S{\ref{c4}} we take up the study of stability of stationary zonal solutions
to \eqref{1.1}, assuming $M$ is radially symmetric and has positive Gauss
curvature.  In \S{\ref{c4s1}}, we apply an Arnold-type approach, and deduce that
a sufficient condition for stability in $H^1(M)$ is that
\beq
w(\xi)-\Omega\chi(\xi)\ \text{ is strictly monotone in }\ \xi,
\label{1.15}
\eeq
where $w=\Delta f$ is the vorticity.  In \S{\ref{c4s2}} we study
the linearization at a steady zonal solution of \eqref{1.1},
or more precisely of the vorticity equation \eqref{1.11}, obtaining
a linear equation of the form $\pa\zeta/\pa t=\Gamma\zeta$.  We establish a
version of the Rayleigh criterion, namely, if the spectrum of $\Gamma$ is not
contained in the imaginary axis, then
\beq
w'(\xi)-\Omega\chi'(\xi)\ \text{ must change sign.}
\label{1.16}
\eeq
Note that \eqref{1.15} and \eqref{1.16}
are almost perfectly complementary.  Nevertheless,
the criterion \eqref{1.16}, while necessary for failure of $\Spec \Gamma\subset
i\RR$, is not sufficient.  This matter is discussed in \S{\ref{c5}}.

In \S{\ref{c5}} we specialize to $M=S^2$ and perform some specific computations,
taking $f(x)=cP_\nu(x_3),\ \nu=2,3,4$.  We make use of classical results on
spherical harmonics to produce infinite matrix representations of the linear
operator $\Gamma$.  We present some analytical results for $\nu=2$ and
some numerical results for $\nu=3$ and $4$, indicating how the Rayleigh-type
criterion \eqref{1.16} is not definitive as a criterion for linear instability.
We also discuss the extent to which stability seems to depend monotonically
on $\Omega$ (or, sometimes, not).

\subsection{Further properties of the operator $B$}\label{c1s1}

The operator $B$, defined in \eqref{1.5},
arose in the form \eqref{1.4} of the Euler equation.  By \eqref{1.7}, we have
\beq
B\in OPS^{-1}(M),
\label{1.1.1}
\eeq
the class of pseudodifferential operators on order $-1$ on $M$.  We record some
other properties of $B$, which will be useful later on.

Since $\delta\tu=0$ on $M\Rightarrow \tu=*df$ for a scalar function $f$
(known as the stream function), uniquely determined up to an additive
constant, it is useful to compute
\beq
B(*df)=\delta\Delta_2^{-1}(d\chi\wedge df).
\label{1.1.2}
\eeq
We have (with $\alpha$ denoting the area form on $M$)
\beq
\aligned
d\chi\wedge df&=-**d\chi\wedge df \\ &=df\wedge*(*d\chi) \\
&=\langle df,*d\chi\rangle\alpha \\ &=\langle df,J\nabla\chi\rangle\alpha \\
&=*Zf,
\endaligned
\label{1.1.3}
\eeq
with the vector field $Z$ given by
\beq
Z=J\nabla\chi.
\label{1.1.4}
\eeq
Note that
\beq
\dv Z=0.
\label{1.1.5}
\eeq
The formula \eqref{1.1.2} yields
\beq
\aligned
B(*df)&=\delta\Delta_2^{-1}*Zf \\
&=*d\Delta_0^{-1}Zf.
\endaligned
\label{1.1.6}
\eeq
Note that \eqref{1.1.5} implies that $Z$ is skew-adjoint and that
$\int_M Zf\, dS=0$.  We see from \eqref{1.1.6} that
\beq
V\cap\Ker B=\{*df:f\in H^1(M),\, Zf=0\},
\label{1.1.7}
\eeq
where
\beq
V=\{\tu\in L^2(M,\Lambda^1):\delta\tu=0\}.
\label{1.1.8}
\eeq

When $M=S^2$, we have the following result, observed in \cite{CM}.

\begin{proposition} \label{p1.1.1}
If $M=S^2$, then $B$ commutes with $\Delta_1$,
the Hodge Laplacian on 1-forms.
\end{proposition}

\demo
In such a case, we have (1.1.4) with $\chi(x)=x_3$, hence $Z=X_3$,
the vector field generating the $2\pi$-periodic rotation about the $x_3$-axis.
Since the flow generated by $X_3$ consists of isometries on $S^2$, $X_3$ and
$\Delta_0$ commute.  Then (by \eqref{1.1.6})
\beq
\aligned
\Delta_1B(*df)&=*d\Delta_0\Delta_0^{-1}X_3f \\
&=*d\Delta_0^{-1}X_3\Delta_0f \\
&=B(*d\Delta_0f) \\&=B\Delta_1*df.
\endaligned
\label{1.1.9}
\eeq
\qed

\section{Basic existence results}\label{c2}

Here we establish existence of solutions $\tu(0)$ to \eqref{1.3}, given
$\tu_0\in H^s(M)$, $s>2$, and $\delta\tu_0=0$,
and we produce estimates on such solutions.
We begin in \S{\ref{c2s1}} with short time existence results.  In
preparation for long time existence results, we derive a vorticity equation
in \S{\ref{c2s2}}.  We show that if $u$ solves \eqref{1.1} and $w(t)=\rot u(t)$, then
\beq
\frac{\pa}{\pa t}(w-\Omega\chi)+\nabla_u(w-\Omega\chi)=0.
\label{2.1}
\eeq
This is a conservation law.  We use it, together with a method pioneered by
\cite{BKM}, in \S{\ref{c2s3}} to establish long time existence of solutions
to \eqref{1.3}.  We show these solutions satisfy the estimate
\beq
\|\tu(t)\|^2_{H^s}\le C\|\tu(0)\|^2_{H^s} \exp\, \exp\,
\Bigl(C_s(\|w(0)\|_{L^\infty}+C|\Omega|)t\Bigr).
\label{2.2}
\eeq
The path from \eqref{2.1} to \eqref{2.2} passes through the estimate
\beq
\|u(t)\|_{\fh^{1,\infty}}\le C\|w(t)\|_{L^\infty}\le C\bigl(\|w(0)\|_{L^\infty}
+2|\Omega|\bigr),
\label{2.3}
\eeq
which will see futher use in \S{\ref{c3}}.  Here,
\beq
\fh^{1,\infty}(M)=\{\tu\in C(M):\nabla \tu\in\bmo(M)\}.
\label{2.4}
\eeq

\subsection{Short time existence}\label{c2s1}

Our approach to the short time solvability of \eqref{1.1},
or equivalently \eqref{1.4}, i.e.,
\beq
\frac{\pa\tu}{\pa t}+P\nabla_u\tu=\Omega B\tu,\quad \tu=P\tu,
\label{2.1.1}
\eeq
with initial data
\beq
\tu(0)=\tu_0\in H^s(M),\quad \delta\tu_0=0,
\label{2.1.2}
\eeq
is to take a mollifier $J_\ep=\varphi(\ep\Delta_1),\ \varphi$ real valued and
in $C_0^\infty(\RR)$, with $\varphi(0)=1$ ($\Delta_1$ the Hodge Laplacian on
1-forms), and solve
\beq
\gathered
\frac{\pa\tu_\ep}{\pa t}+PJ_\ep\nabla_{u_\ep}J_\ep\tu_\ep
=\Omega J_\ep BJ_\ep\tu_\ep, \\
P\tu_\ep=\tu_\ep,\quad \tu_\ep(0)=J_\ep\tu_0.
\endgathered
\label{2.1.3}
\eeq
Compare the treatment in \S{2}, Chapter 17, of \cite{T} (for $\Omega=0$).
Given $\ep>0$, the short-time solvability of \eqref{2.1.3}
is elementary, since \eqref{2.1.3}
is essentially a finite system of ODEs.  The goal is to obtain estimates of
$\tu_\ep(t)$ in $H^s(M)$, for $t$ in some interval, independent of $\ep$, and
pass to the limit.

To start, we have
\beq
\frac{1}{2}\, \frac{d}{dt}\|\tu_\ep(t)\|^2_{L^2}=-(PJ_\ep\nabla_{u_\ep}J_\ep
\tu_\ep,\tu_\ep)+\Omega(J_\ep BJ_\ep\tu_\ep,\tu_\ep).
\label{2.1.4}
\eeq
As in \cite{T}, the first term on the right is $0$ (cf.~(2.3)--(2.5) in
Chapter 17 of \cite{T}).  By \eqref{1.8},
so is the second term on the right side of \eqref{2.1.4}.  Hence
\beq
\|\tu_\ep(t)\|_{L^2}\equiv \|J_\ep \tu_0\|_{L^2}.
\label{2.1.5}
\eeq
This is enough to guarantee global existence of solutions to \eqref{2.1.3},
for each $\ep>0$.

To estimate higher-order Sobolev norms, we bring in
\beq
A=(-\Delta_1)^{1/2},\quad \|\tu\|_{H^s}=\|A^s\tu\|_{L^2}.
\label{2.1.6}
\eeq
Then
\beq
\aligned
\frac{1}{2}\, \frac{d}{dt}\|\tu_\ep(t)\|^2_{H^s}=\
&-(A^sPJ_\ep\nabla_{u_\ep}J_\ep\tu_\ep,A^s\tu_\ep) \\
&+\Omega(A^sJ_\ep BJ_\ep\tu_\ep,A^s\tu_\ep).
\endaligned
\label{2.1.8}
\eeq
Now, by \eqref{1.8},
\beq
\aligned
(A^sJ_\ep BJ_\ep\tu_\ep,A^s\tu_\ep)=\
&(BJ_\ep A^s\tu_\ep,J_\ep A^s\tu_\ep) \\
&+([A^s,B]J_\ep\tu_\ep,J_\ep A^s\tu_\ep) \\
=\ &([A^s,B]J_\ep\tu_\ep,J_\ep A^s\tu_\ep).
\endaligned
\label{2.1.9}
\eeq
Furthermore, since $A^s$ has scalar principal symbol,
\beq
[A^s,B]\in OPS^{s-2}(M).
\label{2.1.10}
\eeq
It follows that the second term on the right side of \eqref{2.1.8} is
\beq
\le C|\Omega|\cdot\|\tu_\ep\|^2_{H^{s-1}}.
\label{2.1.11}
\eeq

We can write the first term on the right side of \eqref{2.1.8} as
\beq
-(A^sPJ_\ep\nabla_{u_\ep}J_\ep\tu_\ep,A^s\tu_\ep)
=-(A^s\nabla_{u_\ep}J_\ep\tu_\ep,A^sJ_\ep\tu_\ep).
\label{2.1.12}
\eeq
In order to make use of the identity $(\nabla_{u_\ep}v,v)=0$, we need to
analyze the commutator $[A^s,\nabla_{u_\ep}]$.  We claim that
\beq
\|[A^s,\nabla_{u_\ep}]J_\ep\tu_\ep\|_{L^2}\le C\|\tu_\ep(t)\|_{C^1}
\|\tu_\ep(t)\|_{H^s},
\label{2.1.13}
\eeq
with $C$ independent of $\ep$.
If $s=2k$ is a positive even integer, this is a
Moser estimate, as in (2.11)--(2.13) of \cite{T}, Chapter 17.  For general real
$s>0$, this is a Kato-Ponce estimate, established in \cite{KP} in the Euclidean
space setting, and in greater generality (directly applicable to the setting
here) in \S{3.6} of \cite{T2}.

In more detail, the KP-estimate gives, for $s>0$,
\beq
\|A^s(fv)-fA^sv\|_{L^2}\le C\|f\|_{C^1}\|v\|_{H^{s-1}}
+C\|f\|_{H^s}\|v\|_{L^\infty}.
\label{2.1.14}
\eeq
We take $v=Xu$, where $X$ is a first-order differential operator, and write
\beq
A^s(fXu)-fX(A^su)=A^s(fXu)-fA^s(Xu)+f[A^s,X]u.
\label{2.1.15}
\eeq
Then \eqref{2.1.14} applies to estimate the first two terms on the right side
of \eqref{2.1.15}, and the $L^2$-norm of the last term is bounded by
$C\|f\|_{L^\infty}\|u\|_{H^s}$.  Thus
\beq
\|[A^s,fX]u\|_{L^2}\le C\|f\|_{C^1}\|u\|_{H^s}+C\|f\|_{H^s}\|u\|_{C^1},
\label{2.1.16}
\eeq
which in turn yields \eqref{2.1.13}.

From \eqref{2.1.13}, we bound the absolute value of \eqref{2.1.12} by
$C\|\tu_\ep(t)\|_{C^1}\|\tu_\ep(t)\|^2_{H^s}$.  Together with \eqref{2.1.11},
this gives
\beq
\frac{d}{dt}\|\tu_\ep(t)\|^2_{H^s}\le C\|\tu_\ep(t)\|_{C^1}\|\tu_\ep(t)\|^2_{H^s}
+C|\Omega|\cdot\|\tu_\ep(t)\|^2_{H^{s-1}}.
\label{2.1.17}
\eeq

On the 2D manifold $M$, $\|\tu\|_{C^1}\le C_s\|\tu\|_{H^s}$, as long as
$s>2$, so \eqref{2.1.17} implies
\beq
\frac{d}{dt}\|\tu_\ep(t)\|^2_{H^3}\le C\|\tu_\ep(t)\|^3_{H^3}
+C|\Omega|\cdot\|\tu_\ep(t)\|^2_{H^3}.
\label{2.1.18}
\eeq
By Gronwall's inequality, we have, for $t\ge 0$,
\beq
\|\tu_\ep(t)\|^2_{H^3}\le y(t),
\label{2.1.19}
\eeq
where $y(t)$ solves
\beq
\frac{dy}{dt}=C(y^{3/2}+|\Omega|y),\quad y(0)=\|\tu_\ep(0)\|^2_{H^3}.
\label{2.1.20}
\eeq
In particular, $\{\tu_\ep(t):0\le t<T_m\}$ is uniformly bounded, in $H^3(M)$,
independent of $\ep\in (0,1]$, as long as
\beq
T_m<\frac{1}{C} \int_{y(0)}^\infty \frac{dy}{y^{3/2}+|\Omega|y}.
\label{2.1.21}
\eeq
For a more explicit (though cruder) upper bound, we can say that
\beq
\|\tu_\ep(t)\|^2_{H^3}+1\le z(t),
\label{2.1.22}
\eeq
where $z(t)$ solves
\beq
\frac{dz}{dt}=C(1+|\Omega|)z^{3/2},\quad z(0)=\|\tu_\ep(0)\|^2_{H^3}+1.
\label{2.1.23}
\eeq
Explicit integration gives
\beq
z(t)=\Bigl(z(0)^{-1/2}-C_1(1+|\Omega|)t\Bigr)^{-2},\quad \text{for }\
0\le t<C_1^{-1}z(0)^{-1/2}(1+|\Omega|)^{-1}.
\label{2.1.24}
\eeq
Consequently,
\beq
\|\tu_\ep(t)\|_{C^1}\le N_\Omega(t)=\frac{C_2}{A-C_1(1+|\Omega|)t},\quad
\text{for }\ 0\le t<T_m,
\label{2.1.25}
\eeq
with
\beq
T_m=\frac{A}{C_1(1+|\Omega|)},\quad A=(\|\tu_0\|^2_{H^3}+1)^{-1/2}.
\label{2.1.26}
\eeq
This plugs into \eqref{2.1.17} to yield
\beq
\aligned
\frac{d}{dt}\|\tu_\ep(t)\|^2_{H^s}&\le CN_\Omega(t)\|\tu_\ep(t)\|^2_{H^s}
+C|\Omega|\cdot \|\tu_\ep(t)\|^2_{H^{s-1}} \\
&\le C(N_\Omega(t)+|\Omega|)\|\tu_\ep(t)\|^2_{H^s},
\endaligned
\label{2.1.27}
\eeq
for $t\in [0,T_m)$, which in turn yields
\beq
\|\tu_\ep(t)\|^2_{H^s}\le \|\tu_\ep(0)\|^2_{H^s}\, \exp\,
C\int_0^t(N_\Omega(s)+|\Omega|)\, ds,\quad \text{for }\ t\in [0,T_m),
\label{2.1.28}
\eeq
an estimate that is uniform in $\ep\in (0,1]$.

With these uniform estimates in hand, one can apply standard techniques,
discussed in Chapter 17 of \cite{T},
to obtain a solution $\tu(t)$ to \eqref{1.3}
in $C([0,T_m),H^s(M))$, given initial data $\tu_0\in H^s(M)$ (satisfying
$\delta\tu_0=0$) as long as $s\ge 3$.  Here $T_m$ is as in \eqref{2.1.26}, and
the solution $\tu(t)$ satisfies an estimate parallel to \eqref{2.1.28}.  Also,
estimates parallel to those produced above establish uniqueness of the
solution $\tu(t)$ and continuous dependence on the initial data $\tu_0$.

\sk
{\sc Remark.} One could replace $H^3$ in \eqref{2.1.18} by $H^{s_0}$ for any
$s_0>2$, and have a local existence result for initial data $\tu_0\in
H^s(M)$ for any $s\ge s_0$.

\sk
{\bf Improved estimates for $M=S^2$}
\newline {}

The estimates \eqref{2.1.26} and \eqref{2.1.28} for the existence time and size of
solutions to \eqref{1.3} exhibit a dependence on $|\Omega|$.  It was observed in
\cite{CM} that one has estimates independent of $\Omega$ when $M=S^2$ is the
standard sphere.  We note the changes in the arguments above that yield this.

The key modification arises in the estimate \eqref{2.1.11} on the second term
on the right side of \eqref{2.1.8}.  If $M=S^2$, then $B$ commutes with $\Delta_1$
(cf.~Proposition \ref{p1.1.1}), hence with $A^s$ and $J_\ep$, so
\beq
(A^sJ_\ep BJ_\ep\tu_\ep,A^s\tu_\ep)=(BA^sJ_\ep\tu_\ep,A^sJ_\ep\tu_\ep)=0,
\label{2.1.29}
\eeq
the latter identity by the skew adjointness of $B$.
Hence \eqref{2.1.11} is replaced
by $0$, and \eqref{2.1.17} is improved to
\beq
\frac{d}{dt}\|\tu_\ep(t)\|^2_{H^s}\le C\|\tu_\ep(t)\|_{C^1}
\|\tu_\ep(t)\|^2_{H^s},
\label{2.1.30}
\eeq
provided $M=S^2$ and $s\ge 3$.  In this case, Gronwall's inequality yields
\eqref{2.1.19} where $y(t)$ solves
\beq
\frac{dy}{dt}=Cy^{3/2},\quad y(0)=\|\tu_\ep(0)\|^2_{H^3}.
\label{2.1.31}
\eeq
This has the explicit solution
\beq
y(t)=(y(0)^{-1/2}-C_1t)^{-2},\quad \text{for }\ 0\le t<C_1^{-1}y(0)^{-1/2}.
\label{2.1.32}
\eeq
Consequently, \eqref{2.1.25} is improved 
\beq
\|\tu_\ep(t)\|_{C^1}\le N(t)=\frac{C_2}{A-C_1t},\quad \text{for }\ 0\le t<T_m,
\label{2.1.33}
\eeq
with
\beq
T_m=C_1^{-1}A,\quad A=\|\tu_0\|^{-1}_{H^3},
\label{2.1.34}
\eeq
and \eqref{2.1.28} is improved to
\beq
\|\tu_\ep(t)\|^2_{H^s}\le \|\tu_\ep(0)\|^2_{H^s}\,
\exp\, C\int_0^t N(s)\, ds,\quad \text{for }\ t\in [0,T_m),
\label{2.1.35}
\eeq
given $s\ge 3$, an estimate that is uniform in both $\ep\in (0,1]$ and $\Omega
\in\RR$.  From here, arguments parallel to those indicated above give a
unique solution to \eqref{1.3}, with initial data $\tu_0\in H^s(S^2)$, for $s\ge 3$,
on $t\in [0,T_m)$, satisfying an estimate parallel to \eqref{2.1.35}.  This result
is similar to Theorem 5.1 in \cite{CM}, except that here (thanks to Moser-type
estimates) the $t$ interval is independent of $s\ge 3$ (and $s$ is not required
to be an integer, and also we can actually fix $s_0>2$ and take $s\ge s_0$).

\subsection{Vorticity equation}\label{c2s2}

The Euler equation \eqref{1.3} can be rewritten as
\beq
\frac{\pa\tu}{\pa t}+\Cal{L}_u\tu=\Omega\chi*\tu+d\Bigl(\frac{1}{2}|u|^2
-p\Bigr),\quad \delta\tu=0,
\label{2.2.1}
\eeq
where $\Cal{L}$ is the Lie derivative.  This follows from the general identity
\beq
\nabla_u\tu=\Cal{L}_u\tu-\frac{1}{2}d|u|^2.
\label{2.2.2}
\eeq
Compare \cite{T}, Chapter 17, \S{1}.  We obtain an equation for the vorticity $w$,
given by
\beq
d\tu=\tw=w\alpha,\quad w=\rot u,
\label{2.2.3}
\eeq
where $\alpha$ is the area form on $M$, by applying the exterior derivative
to \eqref{2.2.1}:
\beq
\aligned
\frac{\pa\tw}{\pa t}+\Cal{L}_u\tw&=\Omega d(\chi*\tu) \\
&=\Omega(d\chi\wedge*\tu) \\&=\Omega \langle d\chi,u\rangle\alpha,
\endaligned
\label{2.2.4}
\eeq
hence (since $\Cal{L}_u\alpha=0$), we have the vorticity equation
\beq
\aligned
\frac{\pa w}{\pa t}+\nabla_u w&=\Omega\langle d\chi,u\rangle \\
&=\Omega\nabla_u\chi.
\endaligned
\label{2.2.5}
\eeq
We can rewrite \eqref{2.2.5} as
\beq
\frac{\pa}{\pa t}(w-\Omega\chi)+\nabla_u(w-\Omega\chi)=0,
\label{2.2.6}
\eeq
which is a conservation law.

It is useful to know that we can reverse the path from \eqref{2.2.1} to
\eqref{2.2.4}.

\begin{proposition} \label{p2.2.1}
Assume $\delta\tu=0$ and set $\tw=d\tu$.
If $\tw$ satisfies \eqref{2.2.4}, then $\tu$ satisfies \eqref{2.2.1}.
\end{proposition}

\demo
For such $\tu$, the Hodge decomposition on $M$ allows us to write
\beq
\frac{\pa\tu}{\pa t}+\Cal{L}_u\tu-\Omega\chi*\tu=dF+\widetilde{G},
\label{2.2.7}
\eeq
where $\widetilde{G}$ is a 1-form on $M$ (for each $t$) satisfying
\beq
\delta\widetilde{G}=0.
\label{2.2.8}
\eeq
Applying the exterior derivative to \eqref{2.2.7} yields
\beq
\frac{\pa\tw}{\pa t}+\Cal{L}_u\tw=\Omega(\nabla_u\chi)\alpha+d\widetilde{G}.
\label{2.2.9}
\eeq
If \eqref{2.2.4} holds, we deduce that
\beq
d\widetilde{G}=0,
\label{2.2.10}
\eeq
which, in concert with \eqref{2.2.8}, implies $\widetilde{G}=0$, since the hypothesis
that $M$ is diffeomorphic to $S^2$ implies $H^1(M,\RR)=0$.
\qed

Also the identity $H^1(M,\RR)=0$ allows us to write
\beq
\delta\tu=0\Longrightarrow \tu=*df,\ \text{ hence }\ u=J\nabla f,
\label{2.2.11}
\eeq
with scalar $f$ (the stream function) determined uniquely up to an additive
constant, which we can specify uniquely by requiring
\beq
\int\limits_M f\, dS=0.
\label{2.2.12}
\eeq
Note that
\beq
w=\Delta f,
\label{2.2.13}
\eeq
and we have
\beq
\tu=*d\Delta^{-1}w,
\label{2.2.14}
\eeq
with $\Delta^{-1}$ defined on scalar functions to annihilate constants and have
range satisfying \eqref{2.2.12}.
We can rewrite the vorticity equation \eqref{2.2.5} as
\beq
\frac{\pa w}{\pa t}+\langle J\nabla f,\nabla(w-\Omega\chi)\rangle=0.
\label{2.2.15}
\eeq

\subsection{Long time existence}\label{c2s3}

As seen in \S{\ref{c2s1}}, if $\tu_0\in H^s(M),\ s\ge 3$ (or even $s>2$) and
$\delta \tu_0=0$, then \eqref{1.3} has a solution $\tu\in C([0,T_0),H^s(M))$,
satisfying
\beq
\aligned
\frac{d}{dt}\|\tu(t)\|^2_{H^s}&\le C\|\tu(t)\|_{C^1}\|\tu(t)\|^2_{H^s}
+C|\Omega|\cdot \|\tu(t)\|^2_{H^{s-1}} \\
&\le C\|\tu(t)\|_{C^1}\|\tu(t)\|^2_{H^s}+C|\Omega|\cdot \|\tu(t)\|^2_{H^{s}},
\endaligned
\label{2.3.1}
\eeq
for some $T_0>0$.  The analysis behind short time existence shows that if
$[0,T_0)\subset\RR^+$ is the maximal interval of existence for $\tu$, with
such regularity, and $T_0<\infty$, then $\|\tu(t)\|_{H^s}$ cannot remain
bounded as $t\nearrow T_0$.

Our goal here is to demonstrate global existence of such a solution.  We use
the method of \cite{BKM} to obtain such long time existence.
(An alternative approach could proceed as in the analysis in \cite{Yud}.)
A key ingredient is
the vorticity equation \eqref{2.2.6}, which is a conservation law.  It implies that,
for all $t\in [0,T_0)$,
\beq
\|w(t)-\Omega\chi\|_{L^\infty}=\|w(0)-\Omega\chi\|_{L^\infty},
\label{2.3.2}
\eeq
where $w(t)=\rot u(t)$.  It follows that
\beq
\|w(t)\|_{L^\infty}\le \|w(0)\|_{L^\infty}+2|\Omega|,
\label{2.3.3}
\eeq
since, by \eqref{1.2}, $\|\chi\|_{L^\infty}=1$.  Now $\|w(t)\|_{L^\infty}$ does
not bound $\|\tu(t)\|_{C^1}$, but, since
\beq
\tu(t)=-\Delta_1^{-1}\delta *w(t),\quad \Delta_1^{-1}\delta\in OPS^{-1}(M),
\label{2.3.4}
\eeq
we have
\beq
\|\tu(t)\|_{C^1_*}\le \|\tu(t)\|_{\fh^{1,\infty}}\le C\|w(t)\|_{L^\infty},
\label{2.3.5}
\eeq
where
\beq
\fh^{1,\infty}(M)=\{\tu\in C(M):\nabla \tu\in \bmo(M)\},
\label{2.3.6}
\eeq
and $C^1_*(M)$ is a Zygmund space.  A variant of the analysis of \cite{BKM}
(See \cite{T2}, Appendix B) gives
\beq
\|\tu\|_{C^1}\le C\|\tu\|_{C^1_*}
\Bigl(1+\log \frac{\|\tu\|_{C^{1+r}}}{\|\tu\|_{C^1_*}}\Bigr),\quad r>0.
\label{2.3.7}
\eeq
Hence
\beq
\|\tu\|_{C^1}\le C\|\tu\|_{C^1_*}(1+\log^+\|u\|_{H^s}),
\label{2.3.8}
\eeq
provided $s>2$.  Plugging into \eqref{2.3.1}, we get, for $s>2$,
\beq
\frac{d}{dt}\|\tu(t)\|^2_{H^s}\le C_s(\|w(0)\|_{L^\infty}+C|\Omega|)
(1+\log^+\|\tu(t)\|^2_{H^s})\|\tu(t)\|^2_{H^s}.
\label{2.3.9}
\eeq
Now, with
\beq
y(t)=\|\tu(t)\|^2_{H^s},\quad A=\|w(0)\|_{L^\infty}+C|\Omega|,
\label{2.3.10}
\eeq
\eqref{2.3.9} says
\beq
\frac{dy}{dt}\le C_sA(1+\log^+y(t))y(t),
\label{2.3.11}
\eeq
so
\beq
\int_{y(0)}^{y(t)} \frac{d\eta}{\eta(1+\log^+\eta)}\le C_sAt.
\label{2.3.12}
\eeq
Now, for $y>e$,
\beq
\int_e^y \frac{d\eta}{\eta\log\eta}=\log\log\eta,
\label{2.3.13}
\eeq
so
\beq
y(t)\le\exp \Bigl(e^{C_sAt}\Bigr),\quad \text{if }\ y(0)=e.
\label{2.3.14}
\eeq
From this we can deduce that
\beq
\|\tu(t)\|^2_{H^s}\le C\|\tu(0)\|^2_{H^s}\, \exp\,\exp
\Bigl(C_s(\|w(0)\|_{L^\infty}+C|\Omega|)t\Bigr).
\label{2.3.15}
\eeq
This estimate implies that $\|\tu(t)\|_{H^s}$ is bounded on $[0,T_0)$ for
all $T_0<\infty$, so we have global existence, with the global estimate
\eqref{2.3.15}.

\sk
{\sc Remark.} As seen in \S{\ref{c2s1}}, when $M=S^2$ one has an improvement on
\eqref{2.3.1}, namely, the term on the right side containing $|\Omega|$ can
be dropped.  However, the term containing $|\Omega|$ in \eqref{2.3.3} persists,
so one does not get a significant improvement on \eqref{2.3.9},
or on \eqref{2.3.15}, in this case.

\section{Bodies with rotational symmetry}\label{c3}

Here we assume $M\subset\RR^3$ is invariant under the group
\beq
R_s=\begin{pmatrix} \cos s & -\sin s & {}\\ \sin s & \cos s & {}\\
{} & {} & 1\end{pmatrix}
\label{3.1}
\eeq
of rotations about the $x_3$-axis.  We also assume that $M$ is diffeomorphic
to $S^2$ and has positive Gauss curvature everywhere.  We consider special
properties of the Euler equation \eqref{1.1} under this hypothesis.  Note that,
if $\chi$ is given by \eqref{1.2}, then
\beq
X_3\chi=0,
\label{3.2}
\eeq
where $X_3$ denotes the vector field on $M$ generating the flow \eqref{3.1}.
In fact, we can write
\beq
\chi(x)=\chi(x_3),
\label{3.3}
\eeq
where on the right side $\chi\in C^\infty([-a,a])$, assuming
\beq
x_3:M\rightarrow [-a,a],\quad a=\max\limits_M\, x_3,\quad
-a=\min\limits_M\, x_3,
\label{3.4}
\eeq
which can be arranged by a translation.  In such a case, the vector field
\beq
Z=J\nabla\chi,
\label{3.5}
\eeq
arising in \eqref{1.1.4}, is parallel to $X_3$.  In fact,
\beq
Z=\Phi X_3,\quad \Phi(x)=-\chi'(x_3).
\label{3.6}
\eeq
If $M=S^2$, then $\chi(x)=x_3$, and $Z=-X_3$.

In \S{\ref{c3s1}} we study stationary solutions to \eqref{1.1},
i.e., solutions that are independent of $t$.
We show that if $f\in C^\infty(M)$ is a zonal function,
i.e., $X_3f=0$, then the associated divergence-free vector field $u=J\nabla f$
(which we call a zonal field) is a stationary solution to \eqref{1.1}, for all
$\Omega$.  We also give examples of stationary solutions that are not zonal
fields.

In \S{\ref{c3s2}} we return to time-dependent solutions and study time averages
\beq
\Cal{A}_{S,T}\tu=\frac{1}{T}\int_S^{S+T}\tu(t)\, dt,
\label{3.7}
\eeq
where $\tu$, solving \eqref{1.3},
is the 1-form counterpart to the vector field $u$, solving \eqref{1.1}.
We construct a projection $\Pi$ from forms satisfying $\delta
\tu=0$ onto the subspace of zonal forms, and produce estimates on
\beq
\|(I-\Pi)\Cal{A}_{S,T}\tu\|_{H^{-3,q}},
\label{3.8}
\eeq
for $1<q<\infty$, involving negative powers of $|\Omega|$ (see \eqref{3.2.36}).
Our interest in such estimates was stimulated by the paper \cite{CM}, which
produced estimates on $\|(I-\Pi)\Cal{A}_{S,T}\tu\|_{H^{k,2}}$ for positive $k$,
valid however over a limited range of $S,T$.  Our goal was to produce
uniform estimates, valid for all time.  One key to this was to use the
long-time existence and estimates, from \S{\ref{c2s3}}, in place of short-time
existence and estimates, from \S{\ref{c2s1}}.
Also, the results of \cite{CM} were
derived for a rotating sphere.  Since fast rotating planets have noticeable
bulges at the equator, we were motivated to treat more general rotationally
symmetric cases $M$.

In \S{\ref{c3s3}}, we produce another conservation law.  Namely, with $\xi\in
C^\infty(M)$ satisfying
\beq
X_3=-J\nabla\xi,
\label{3.9}
\eeq
if $u$ satisfies \eqref{1.1} and $w=\rot u$ is the associated vorticity, then
\beq
\int\limits_M \xi(x)w(t,x)\, dS(x)
\label{3.10}
\eeq
is independent of $t$.  Note that $M=S^2\Rightarrow \xi(x)=\chi(x)=x_3$.
Such a conservation law appears in [CaM] in the special case $M=S^2$.
We derive it here for a similar reason as [CaM], as a tool to use in an
Arnold-type analysis of stability of stationary solutions to \eqref{1.1};
see \S{\ref{c4}}.

In \S{\ref{c3s4}} we discuss computations of $\chi$ and $\xi$, first for a
general surface of revolution
\beq
x_1^2+x_2^2=r(x_3)^2,
\label{3.11}
\eeq
and then, more explicitly, for ellipsoids of revolution
\beq
x_1^2+x_2^2+\Bigl(\frac{x_3}{a}\Bigr)^2=1.
\label{3.12}
\eeq
Section \ref{c3s5} establishes smoothness of various functions, such as
$\chi(x_3)$ and $\xi(x_3)$, making essential use of the positive Gauss
curvature assumption.

\subsection{Stationary solutions}\label{c3s1}

A stationary solution to \eqref{1.1} is one for which $\pa u/\pa t=0$.
In such a case, $w=\rot u$ satisfies \eqref{2.2.5}
with $\pa w/\pa t=0$.  Hence, by \eqref{2.2.15},
\beq
\langle J\nabla f,\nabla(w-\Omega\chi)\rangle=0,
\label{3.1.1}
\eeq
where
\beq
w=\Delta f,\quad \tu=*df,
\label{3.1.2}
\eeq
which determines $f$ uniquely, up to an additive constant.  The equation
\eqref{3.1.1} is equivalent to
\beq
\nabla(\Delta f-\Omega\chi)\, \| \, \nabla f.
\label{3.1.3}
\eeq
By Proposition \ref{p2.2.1}, whenever $f$ satisfies \eqref{3.1.3},
which implies \eqref{3.1.1}, then $\tu$, defined by \eqref{3.1.2},
is a stationary solution to \eqref{1.3}.  This gives
the following class of stationary solutions.  We say $f\in C^\infty(M)$ is a
zonal function if $X_3=0$, where the vector field $X_3$ generates $2\pi$-periodic
rotation about the $x_3$-axis.  Then we say $u=J\nabla f$ is a zonal velocity
field.

\begin{proposition} \label{p3.1.1}
Assume $M\subset\RR^3$ is a smooth, compact surface,
with positive Gauss curvature, and radially symmeric about the $x_3$-axis.
If $f\in C^\infty(M)$ is a zonal function, then $u=J\nabla f$ is a stationary
solution to \eqref{1.1}, for all $\Omega$.
\end{proposition}

\demo
Under our hypothesis, we have $\chi=\chi(x_3),\ f=f(x_3)$, and
$w=w(x_3)$, so \eqref{3.1.3} holds.
\qed

\begin{corollary} \label{c3.1.2}
With $M$ as in Proposition \ref{p3.1.1},
\beq
\text{each $\tu\in V\cap\Ker B$ is a stationary solution to \eqref{1.3}.}
\label{3.1.4}
\eeq
\end{corollary}

\demo
Recall that $V\cap \Ker B$ is given by \eqref{1.1.7}, with $Z=J\nabla
\chi$, as in \eqref{1.1.4}.  The geometrical hypothesis on $M$
made above implies
\beq
Z=\Phi X_3,
\label{3.1.5}
\eeq
for some nowhere vanishing $\Phi\in C^\infty(M)$, which yields \eqref{3.1.4}.
\qed

While our study of stationary solutions to the Euler equation will focus on
zonal functions, we mention that there are stationary solutions to \eqref{1.3}
that are not of the form \eqref{3.1.4}.
We give examples when $M=S^2$, the standard sphere.
To get started, note that \eqref{3.1.3} holds whenever there is a smooth
$\psi:\RR\rightarrow\RR$ such that
\beq
\Delta f=\psi(f)+\Omega x_3\quad \text{(given $M=S^2$)}.
\label{3.1.6}
\eeq
We will apply this with $\psi(f)=-\lambda_kf$, where $\lambda_k$ is chosen from
\beq
\Spec (-\Delta)=\{\lambda_k=k^2+k:k=0,1,2,3,\dots\}.
\label{3.1.7}
\eeq
Note that $x_3$ is an eigenfunction of $-\Delta$ with eigenvalue $\lambda_1=2$.
Thus we assume $k\ge 2$.  Then \eqref{3.1.6} becomes
\beq
(\Delta+\lambda_k)f=\Omega x_3.
\label{3.1.8}
\eeq
As long as $\lambda_k\neq 2$, \eqref{3.1.8} has solutions, and the general
solution is of the form
\beq
f=\frac{\Omega}{\lambda_k-2}x_3+g_k,\quad g_k\in\Ker (\Delta+\lambda_k).
\label{3.1.9}
\eeq
Thus $g$ is the restriction to $S^2$ of a harmonic polynomial, homogeneous of
degree $k$.  For example, we can take
\beq
\aligned
k&=2,\quad \lambda_k=6,\quad g_k(x)=x_1^2-x_2^2, \\
k&=3,\quad \lambda_k=12,\quad g_k(x)=\Re(x_1+ix_2)^3,
\endaligned
\label{3.1.10}
\eeq
etc.  For such $f$ as in \eqref{3.1.9}, we have
\beq
u=-\frac{\Omega}{\lambda_k-2}X_3+J\nabla g_k
\label{3.1.11}
\eeq
as a stationary solution to \eqref{1.1}.

\subsection{Time averages of solutions}\label{c3s2}

As before, $M\subset\RR^3$ is a smooth compact surface of positive Gauss
curvature that is radially symmetric about the $x_3$-axis.  We take $u$ to
be a smooth solution to the Euler equation \eqref{1.1}, so $\tu$ solves
\eqref{1.3}, or equivalently \eqref{1.4}, i.e.,
\beq
\frac{\pa\tu}{\pa t}=\Omega B\tu-P\nabla_u\tu,\quad \delta\tu=0.
\label{3.2.1}
\eeq
Given $S,T\in (0,\infty)$, we want to investigate the time-averaged
field
\beq
\Cal{A}_{S,T}\tu=\frac{1}{T} \int_S^{S+T} \tu(t)\, dt.
\label{3.2.2}
\eeq
In particular, we investigate the extent to which it can be shown that
$\Cal{A}_{S,T}\tu$ is close to a zonal field, particularly for
large $\Omega$.

We start by integrating \eqref{3.2.1} over $t\in [S,S+T]$, obtaining
\beq
\frac{\tu(S+T)-\tu(S)}{T}=\Omega B\Cal{A}_{S,T}\tu-\Cal{A}_{S,T}P\nabla_u\tu,
\label{3.2.3}
\eeq
or
\beq
B\Cal{A}_{S,T}\tu=\frac{1}{\Omega}\Bigl\{\frac{\tu(S+T)-\tu(S)}{T}
+\Cal{A}_{S,T}P\nabla_u\tu\Bigr\}.
\label{3.2.4}
\eeq
We want to show that if the left side of \eqref{3.2.4} is small (in some
sense) then $\Cal{A}_{S,T}\tu$ is close to being a zonal field.

To do this, we will produce an operator $\Pi$ with the property that,
for each $s\in\RR,\ p\in (1,\infty)$, $\Pi$ is a projection of
\beq
V^{s,p}=\{\tu\in H^{s,p}(M,\Lambda^1):\delta\tu=0\}\ \text{ onto }\
V^{s,p}\cap\Ker B.
\label{3.2.5}
\eeq
As in \eqref{1.1.7},
\beq
V^{s,p}\cap \Ker B=\{*df:f\in H^{s+1,p}(M),\ Zf=0\},
\label{3.2.6}
\eeq
where, as in \eqref{1.1.4},
\beq
Z=J\nabla\chi.
\label{3.2.7}
\eeq
As noted in \eqref{3.1.5}, the geometrical hypothesis on $M$ implies
$Z=\Phi X_3$ for some nowhere vanishing $\Phi\in C^\infty(M)$.  Consequently,
\beq
V^{s,p}\cap\Ker B=\{*df:f\in H^{s+1,p}(M),\ X_3f=0\}.
\label{3.2.8}
\eeq
We will define $\Pi$ by
\beq
\Pi(*df)=*d\wtp f,
\label{3.2.9}
\eeq
where
\beq
\wtp:H^{s+1,p}(M)\longrightarrow H^{s+1,p}(M)
\label{3.2.10}
\eeq
is the projection onto $\Ker X_3$ given by
\beq
\wtp f(x)=\frac{1}{2\pi}\int_0^{2\pi} f(R_sx)\, ds,
\label{3.2.11}
\eeq
with
\beq
R_s=\begin{pmatrix} \cos s & -\sin s & {}\\ \sin s & \cos s & {}\\
{} & {} & 1\end{pmatrix}.
\label{3.2.12}
\eeq
The following result is the key to exploiting \eqref{3.2.4}.

\begin{proposition} \label{p3.2.1}
If $M$ is a body of rotation about the $x_3$-axis,
with positive Gauss curvature, then, for $q\in (1,\infty),\ s\in\RR$,
\beq
\|(I-\Pi)u\|_{H^{s-2,q}}\le C_{q,s}\|Bu\|_{H^{s,q}}.
\label{3.2.13}
\eeq
\end{proposition}

Given $B(*df)=*d\Delta_0^{-1}Zf$, from \eqref{1.1.6}, and given
\eqref{3.2.9}, it suffices to prove the following.

\begin{proposition} \label{p3.2.2}
In the setting of Proposition 3.2.1, for
$p\in (1,\infty),\ \sigma\in\RR$,
\beq
\|(I-\wtp)f\|_{H^{\sigma,p}}\le C\|Zf\|_{H^{\sigma,p}}.
\label{3.2.14}
\eeq
\end{proposition}
\demo
Given \eqref{3.2.7}, it suffices to show that
\beq
\|(I-\wtp)f\|_{H^{\sigma,p}}\le C\|X_3f\|_{H^{\sigma,p}}.
\label{3.2.15}
\eeq
The formula \eqref{3.2.11} implies that $\wtp$ is bounded on $L^p(M)$ for all
$p\in (1,\infty)$ and commutes with $\Delta_0$.  Also $X_3$ commutes with
$\Delta_0$, so it suffices to establish
\beq
\|(I-\wtp)f\|_{L^p}\le C\|X_3f\|_{L^p}.
\label{3.2.16}
\eeq
To get this, it suffices to construct a bounded map $T$ on $L^p$ such that
\beq
TX_3=I-\wtp.
\label{3.2.17}
\eeq
We construct $T$ in the form
\beq
Tg=\frac{1}{2\pi} \int_{-\pi}^\pi \psi(s)e^{sX_3}g\, ds,
\label{3.2.18}
\eeq
where $e^{sX_3}g(x)=g(R_s x)$.  We have
\beq
\aligned
TX_3f&=\frac{1}{2\pi}\int_{-\pi}^\pi \psi(s)e^{sX_3}X_3f\, ds \\
&=\frac{1}{2\pi}\int_{-\pi}^\pi \psi(s)\,\frac{d}{ds}e^{sX_3}f\, ds \\
&=-\frac{1}{2\pi}\int_{-\pi}^\pi \psi'(s)e^{sX_3}f\, ds,
\endaligned
\label{3.2.19}
\eeq
provided $\psi(-\pi)=\psi(\pi)$.  To get \eqref{3.2.17}, we want
\beq
-\psi'(s)=2\pi\delta(s)-1.
\label{3.2.20}
\eeq
This is achieved by
\beq
\aligned
\psi(s)=\ &s+\pi,\quad -\pi<s<0, \\
&s-\pi,\quad 0<s<\pi.
\endaligned
\label{3.2.21}
\eeq
Since $\psi\in L^1(-\pi,\pi)$, $T$ is bounded on each $L^p(M)$.  This establishes
\eqref{3.2.17}, hence \eqref{3.2.16},
hence \eqref{3.2.14}, hence \eqref{3.2.13}.  The proof of
Proposition \ref{p3.2.1} is complete.
\qed

Applying Proposition \ref{p3.2.1} to \eqref{3.2.4}, we have
\beq
\aligned
&\|(I-\Pi)\Cal{A}_{S,T}\tu\|_{H^{s-2,q}} \\
&\le \frac{C}{|\Omega|}\Bigl\{T^{-1}\|\tu(S+T)-\tu(S)\|_{H^{s,q}}
+\Cal{A}_{S,T}\|\nabla_u\tu\|_{H^{s,q}}\Bigr\},
\endaligned
\label{3.2.22}
\eeq
when $\tu$ solves \eqref{3.2.1}.  Now (cf.~\cite{T}, Chapter 17, (2.23)),
\beq
\dv u=0\Longrightarrow \nabla_u u=\dv(u\otimes u),
\label{3.2.23}
\eeq
so
\beq
\|\nabla_{u(t)}\tu(t)\|_{H^{s,q}}\le C\|\tu(t)\otimes\tu(t)\|_{H^{s+1,q}}.
\label{3.2.24}
\eeq
Meanwhile, as seen in \S{\ref{c2s3}}, with $w(t)=\rot u(t)$,
\beq
\aligned
\|\tu(t)\|_{\fh^{1,\infty}}&\le C\|w(t)\|_{L^\infty}
\le C(\|w(0)\|_{L^\infty}+2|\Omega|), \\
\|u(t)\|_{L^2}&=\|u(0)\|_{L^2}.
\endaligned
\label{3.2.25}
\eeq
We produce further estimates by interpolation.  For starters,
\beq
\|u(t)\|_{H^{1/2,4}}\le C\|u(0)\|_{L^2}^{1/2}(\|w(0)\|_{L^\infty}
+2|\Omega|)^{1/2}.
\label{3.2.26}
\eeq
In formulas below, in order to simplify the notation, we set
\beq
\Omega_\theta=\|u(0)\|_{L^2}^{1-\theta}(\|w(0)\|_{L^\infty}+2|\Omega|)^\theta.
\label{3.2.27}
\eeq
Then \eqref{3.2.26} becomes
\beq
\|u(t)\|_{H^{1/2,4}}\le C\Omega_{1/2}.
\label{3.2.28}
\eeq
Since \eqref{3.2.24} is quadratic in $u(t)$,
we want to interpolate \eqref{3.2.28} with
the $L^2$ estimate in \eqref{3.2.25}.  We get, for $0\le\theta\le 1$,
\beq
\aligned
\|u(t)\|_{H^{\theta/2,p(\theta)}}&\le C\|u(t)\|_{L^2}^{1-\theta}
\|u(t)\|_{H^{1/2,4}}^\theta \\
&\le C\Omega_{\theta/2},
\endaligned
\label{3.2.29}
\eeq
with
\beq
\frac{1}{p(\theta)}=\frac{1-\theta}{2}+\frac{\theta}{4}=\frac{2-\theta}{4}.
\label{3.2.30}
\eeq
Now, for $0\le\theta<1$,
\beq
H^{\theta/2,p(\theta)}(M)\subset L^{r(\theta)}(M),
\label{3.2.31}
\eeq
with
\beq
r(\theta)=\frac{2p(\theta)}{2-\theta p(\theta)/2}
=\frac{4}{\frac{4}{p(\theta)}-\theta}=\frac{4}{2-2\theta}=\frac{2}{1-\theta},
\label{3.2.32}
\eeq
so
\beq
\|u(t)\|_{L^{2/(1-\theta)}}\le C_\theta \Omega_{\theta/2},\quad 0\le\theta<1.
\label{3.2.33}
\eeq
Hence
\beq
\|u(t)\otimes u(t)\|_{L^{1/(1-\theta)}}\le C_\theta \Omega^2_{\theta/2},
\label{3.2.34}
\eeq
so, by (3.2.24),
\beq
\|\nabla_{u(t)}\tu(t)\|_{H^{-1,q(\theta)}}\le C_\theta \Omega^2_{\theta/2},
\quad q(\theta)=\frac{1}{1-\theta}.
\label{3.2.35}
\eeq
This indicates taking $s=-1$ in \eqref{3.2.22}, and leads to the following result.

\begin{proposition} \label{p3.2.3}
Let $\tu$ be a global solution to \eqref{1.3}, with
initial data in $H^s(M)$, $s>2$.  Then there exists $C_\theta<\infty$,
independent of $S,T,\Omega\in (0,\infty)$, such that
\beq
\aligned
\|(I-\Pi)\Cal{A}_{S,T}&\tu\|_{H^{-3,q(\theta)}} \\
\le \frac{C_\theta}{|\Omega|}&\Bigl\{
T^{-1}\|\tu(S+T)-\tu(S)\|_{H^{-1,q(\theta)}} \\
&+C\|u(0)\|_{L^2}^{2-\theta}(\|w(0)\|_{L^\infty}+2|\Omega|)^\theta\Bigr\},
\endaligned
\label{3.2.36}
\eeq
for $0<\theta<1$, with $q(\theta)$ as in \eqref{3.2.35}.
\end{proposition}

\noindent
{\sc Remark.}  We have
$$
\aligned
\|\tu(S+T)-\tu(S)\|_{H^{-1,q(\theta)}}
&\le C_\theta \|\tu(S+T)-\tu(S)\|_{L^2} \\
&\le 2C_\theta \|u(0)\|_{L^2}.
\endaligned
$$

Let us look at some special cases to which Proposition \ref{p3.2.3} applies.
First, if $\tu$ is a zonal field, then, as seen in \S{\ref{c3s1}},
$\tu$ is a stationary solution to \eqref{1.3},
and consequently the left side of \eqref{3.2.36} vanishes.
By contrast, recall the non-zonal stationary solutions on $S^2$ given by
\eqref{3.1.11}, i.e.,
\beq
u=-\frac{\Omega}{\lambda_k-2}X_3+J\nabla g_k,
\label{3.2.37}
\eeq
with $\lambda_k=k^2+k>2$ an eigenvalue of $-\Delta$ and $g_k$ a non-zonal
$\lambda_k$-eigenfunction on $S^2$, as in \eqref{3.1.10}.  In such a case,
\beq
\Cal{A}_{S,T}\tu\equiv \tu,\quad \text{so }\ (I-\Pi)\Cal{A}_{S,T}\tu=*dg_k.
\label{3.2.38}
\eeq
To make contact with the estimate \eqref{3.2.36},
let us suppose that $\|\tu\|_{L^2}
=\|\tu(0)\|_{L^2}\approx 1$.  Then
\beq
\frac{\Omega}{\lambda_2-2}\le C\quad \text{and }\ \|\lambda_k^{1/2}g_k\|_{L^2}
\le C.
\label{3.2.39}
\eeq
It follows that
\beq
\aligned
\|(I-\Pi)\Cal{A}_{S,T}\tu\|_{H^{-s}}&\le C\|dg_k\|_{H^{-s}} \\
&\le C\lambda_k^{(1-s)/2} \|g_k\|_{L^2} \\
&\le C\lambda_k^{-s/2} \\ &\le C\Omega^{-s/2}.
\endaligned
\label{3.2.40}
\eeq
For $s\in (0,1)$, this estimate is stronger than \eqref{3.2.36},
but of a similar flavor.  Of course, since \eqref{3.2.37}
covers only a special class of stationary
solutions to \eqref{1.1}, it is not surprising that estimates here are better
than the general estimates guaranteed by \eqref{3.2.36}.

\subsection{Another conservation law}\label{c3s3}

As usual, $M\subset\RR^3$ is a surface, diffeomorphic to $S^2$, with positive
Gauss curvature, and invariant under
the group of rotations about the $x_3$-axis generated by $X_3$.  
As a consequence,
\beq
\text{$\chi$ is a smooth function of $x_3$ and }\ \frac{d\chi}{dx_3}\ge\alpha
>0\ \text{for }\ x\in [-1,1].
\label{3.3.1}
\eeq
Next, since $X_3$ generates a flow by isometries on $M$, we have $\dv X_3=0$ on
$M$, so there exists $\xi\in C^\infty(M)$ such that
\beq
J\nabla\xi=-X_3.
\label{3.3.2}
\eeq
Clearly $X_3\xi=0$.  As a further consequence of our geometric hypothesis,
\beq
\text{$\xi$ is a smooth function of $x_3$ and }\ \frac{d\xi}{dx_3}\ge\alpha
>0\ \text{ for }\ x_3\in [-1,1].
\label{3.3.3}
\eeq
(If $M=S^2$, then $\chi=\xi=x_3$.)

We now aim to establish the following conservation law.

\begin{proposition} \label{p3.3.1}
Under the hypotheses on $M$ made above, if $u(t)$
solves \eqref{1.1} and $\rot u=w$, then
\beq
\int\limits_M \xi(x)w(t,x)\, dS(x)\ \text{ is independent of $t$.}
\label{3.3.4}
\eeq
\end{proposition}

\demo
From the vorticity equation \eqref{2.2.5}, we have
\beq
\aligned
\frac{d}{dt} \int\limits_M \xi w(t)\, dS
&=\int\limits_M \xi\frac{\pa w}{\pa t}\, dS \\
&=\int\limits_M \xi\nabla_u(\Omega\chi-w)\, dS \\
&=-\int\limits_M \xi(\nabla_uw)\, dS+\Omega\int\limits_M\xi\nabla_u\chi\, dS.
\endaligned
\label{3.3.5}
\eeq
Note that \eqref{3.3.1}--\eqref{3.3.3}
imply $\xi$ is a smooth function of $\chi$; write
$\xi=\xi(\chi)$.  Then $\xi\nabla_u\chi=\nabla_uG(\chi)$ where $G'(\chi)=
\xi(\chi)$.  Hence
\beq
\int\limits_M \xi\nabla_u\chi\, dS=\int\limits_M\nabla_uG(\chi)\, dS=0,
\label{3.3.6}
\eeq
since $\dv u=0$ implies $\nabla_u$ is skew adjoint, and $\nabla_u1=0$.
Next,
\beq
\aligned
-\int\limits_M \xi(\nabla_uw)\, dS
&=\int\limits_M (\nabla_u\xi)w\, dS \\
&=\int\limits_M \langle J\nabla f,\nabla\xi\rangle(\Delta f)\, dS \\
&=\int\limits_M (X_3f)(\Delta f)\, dS \\
&=(X_3f,\Delta f).
\endaligned
\label{3.3.7}
\eeq
Now, since $X_3$ commutes with $\Delta$ and is skew-adjoint,
\beq
(X_3f,\Delta f)=-(X_3(-\Delta)^{1/2}f,(-\Delta)^{1/2}f)=0.
\label{3.3.8}
\eeq
It follows that
\beq
\frac{d}{dt}\int\limits_M \xi w(t)\, dS=0,
\label{3.3.9}
\eeq
proving Proposition \ref{p3.3.1}.
\qed

\subsection{Computation of $\chi$ and $\xi$}\label{c3s4}

Let the surface of revolution $M\subset\RR^3$ be given by
\beq
x_1^2+x_2^2=r(x_3)^2,
\label{3.4.1}
\eeq
i.e., $u(x)=0$ with $u(x)=x_1^2+x_2^2-r(x_3)^2$.  We have
\beq
\nabla u(x)=2(x_1,x_2,-r(x_3)r'(x_3)),
\label{3.4.2}
\eeq
so the unit outward normal to $M$ is
\beq
N(x)=\frac{\nabla u(x)}{|\nabla u(x)|}
=\frac{1}{r(x_3)\sqrt{1+r'(x_3)^2}}(x_1,x_2,-r(x_3)r'(x_3)).
\label{3.4.3}
\eeq
Hence
\beq
\chi(x)=N(x)\cdot e_3=-\frac{r'(x_3)}{\sqrt{1+r'(x_3)^2}}.
\label{3.4.4}
\eeq

We next look for $\xi\in C^\infty(M)$, satisfying $J\nabla\xi=-X_3$.
Clearly $\xi$ is to be a function of $x_3$, and then the desired condition is
$|\nabla\xi|=|X_3|=r(x_3)$.  Now
\beq
\nabla\xi=\xi'(x_3)\nabla x_3,
\label{3.4.5}
\eeq
and $\nabla x_3$ is the orthogonal projection onto $T_x M$ of $e_3$, so
\beq
|\nabla x_3|^2=1-(e_3\cdot N(x))^2=1-\chi(x_3)^2.
\label{3.4.6}
\eeq
Hence $\xi$ is defined by the condition
\beq
\xi'(x_3)=\frac{r(x_3)}{\sqrt{1-\chi(x_3)^2}}.
\label{3.4.7}
\eeq
Bringing in \eqref{3.4.4}, we obtain
\beq
\chi(x_3)^2=\frac{r'(x_3)^2}{1+r'(x_3)^2},
\label{3.4.8}
\eeq
hence
\beq
1-\chi(x_3)^2=\frac{1}{1+r'(x_3)^2},
\label{3.4.9}
\eeq
so
\beq
\xi'(x_3)=r(x_3)\sqrt{1+r'(x_3)^2}.
\label{3.4.10}
\eeq
Note that this yields an interesting geometrical interpretation of $\xi$.
Namely, up to an additive constant, $2\pi \xi(x_3)$ is the area of
$$
\{(x,y,z)\in M:z\le x_3\}.
$$

\sk
{\bf Special case: ellipsoids of revolution}
\newline {}

We specialize our calculations to the case where $M$ is given by
\beq
x_1^2+x_2^2+\Bigl(\frac{x_3}{a}\Bigr)^2=1,\quad a>0.
\label{3.4.11}
\eeq
Thus $r(x_3)=\sqrt{1-x_3^2/a^2}$ in \eqref{3.4.1}.  It follows that
\beq
r'(x_3)=-\frac{x_3}{a^2}\Bigl(1-\frac{x_3^2}{a^2}\Bigr)^{-1/2},
\label{3.4.12}
\eeq
hence
\beq
1+r'(x_3)^2=\frac{1-\beta x_3^2}{1-x_3^2/a^2},
\label{3.4.13}
\eeq
with
\beq
\beta=\frac{1}{a^2}-\frac{1}{a^4}.
\label{3.4.14}
\eeq
Thus, by \eqref{3.4.4},
\beq
\chi(x_3)=\frac{x_3}{a^2}\, \frac{1}{\sqrt{1-\beta x_3^2}},
\label{3.4.15}
\eeq
and, by \eqref{3.4.10},
\beq
\xi'(x_3)=\sqrt{1-\beta x_3^2}.
\label{3.4.16}
\eeq
For these ellipsoids, $x_3\in [-a,a]$, and we have $\beta x_3^2<1$, so the
formulas \eqref{3.4.15}--\eqref{3.4.16} clearly exhibit
$\chi$ and $\xi$ as elements of
$C^\infty([-a,a])$.  Note that
\beq
0<a<1\Rightarrow \beta<0,\quad a=1\Rightarrow \beta=0,\quad a>1\Rightarrow
\beta\in (0,1).
\label{3.4.17}
\eeq
In case $M$ is the unit sphere, so $a=1$, we get
\beq
\chi(x_3)=\xi(x_3)=x_3,
\label{3.4.18}
\eeq
as expected.  Ellipsoidal planets that bulge at the equator have $a<1$.

\subsection{Smoothness issues}\label{c3s5}

As we have seen, \eqref{3.4.15}--\eqref{3.4.16} exhibit $\chi$ and $\xi$ as
elements of $C^\infty([-a,a])$ when $M$ is an ellipsoid of the form
\eqref{3.4.11}.  To extend this, assume
\beq
x_3:M\longrightarrow[-a,a],\quad a=\max\limits_M\, x_3,\quad
-a=\min\limits_M\, x_3.
\label{3.5.1}
\eeq
We claim that if $M\subset\RR^3$ is a surface of revolution about the
$x_3$-axis, diffeomorphic to $S^2$, and with positive Gauss curvature,
the following holds:
\beq
\aligned
&\text{If $\Phi\in C^\infty(M)$ is invariant under rotation about the
$x_3$-axis, then} \\
&\Phi(x)=\varphi(x_3)\ \text{ with }\ \varphi\in C^\infty([-a,a]).
\endaligned
\label{3.5.2}
\eeq
Recall that $\chi,\xi\in C^\infty(M)$.  We first show how, under the
condition \eqref{3.5.2}, the conclusion
$\chi,\xi\in C^\infty([-a,a])$ is manifested
in the formulas \eqref{3.4.4} and \eqref{3.4.10}.
First note that, by \eqref{3.4.1}, $r(x_3)^2$
belongs to $C^\infty(M)$, so, under the condition \eqref{3.5.2},
\beq
r(x_3)^2=\rho(x_3),\quad \rho\in C^\infty([-a,a]).
\label{3.5.3}
\eeq
Let us bring in the hypothesis that the curvature of $M$ is nonzero at the
poles $(0,0,\pm a)$, so
\beq
\rho'(\pm a)\neq 0.
\label{3.5.4}
\eeq
Note how these results can be directly verified in case \eqref{3.4.11}.
Generally, we have
\beq
r'(x_3)=\frac{1}{2}\rho'(x_3)\rho(x_3)^{-1/2}.
\label{3.5.5}
\eeq
Hence
\beq
\sqrt{1+r'(x_3)^2}=\frac{1}{2} \sqrt{4+\frac{\rho'(x_3)^2}{\rho(x_3)}}.
\label{3.5.6}
\eeq
Thus, by \eqref{3.4.4},
\beq
\chi(x_3)=-\frac{\rho'(x_3)}{\sqrt{4\rho(x_3)+\rho'(x_3)^2}},
\label{3.5.7}
\eeq
and, by \eqref{3.4.10},
\beq
\xi'(x_3)=\frac{1}{2}\sqrt{4\rho(x_3)+\rho'(x_3)^2}.
\label{3.5.8}
\eeq
By virtue of \eqref{3.5.4}, these formulas clearly give
$\chi,\xi'\in C^\infty([-a,a])$.

Geometric hypotheses guaranteeing that the condition \eqref{3.5.2} holds
are pretty straightforward away
from the extreme values $x_3=\pm a$.
Let us verify \eqref{3.5.3} under the following
explicit hypothesis on $M$ near the poles $(0,0,\pm a)$.  Namely, we asume
$M$ is given near the poles as                             
\beq
x_3=\pm a+\varphi_{\pm}(x_1^2+x_2^2),
\label{3.5.9}
\eeq
for $x_1^2+x_2^2<\delta$, with
\beq
\varphi_{\pm}\in C^\infty([0,\delta)),\quad \varphi_{\pm}(0)=0,\quad
\mp \varphi'_{\pm}(0)>0.
\label{3.5.10}
\eeq
Then, by \eqref{3.4.1},
\beq
x_3=\pm a+\varphi_{\pm}(\rho(x_3)),
\label{3.5.11}
\eeq
so
\beq
1=\varphi'_{\pm}(\rho(x_3))\rho'(x_3).
\label{3.5.12}
\eeq
Hence
\beq
\frac{d\rho}{\varphi'_{\pm}(\rho)}=dx_3,
\label{3.5.13}
\eeq
which yields $\rho$ satisfying \eqref{3.5.3}--\eqref{3.5.4} near $x_3=\pm a$.

After these observations, we are now ready to prove a clean smoothness
result.

\begin{proposition} \label{p3.5.1}
Let $M\subset\RR^3$ be a smooth, compact surface,
invariant under rotation about the $x_3$-axis.
Assume \eqref{3.5.1} holds.  Assume that the Gauss curvature
$K(x)>0$ for all $x\in M$.  Then \eqref{3.5.2}
holds.  That is, if $\Phi\in C^\infty(M)$ is invariant under rotation about
the $x_3$-axis, then $\Phi(x)=\varphi(x_3)$ with $\varphi\in C^\infty([-a,a])$.
\end{proposition}

\demo
The conclusion about $\Phi(x)=\varphi(x_3)$ is straightforward except
for smoothness at $x_3=\pm a$, so we concentrate on that.  Near the poles
$(0,0,\pm a)$, $(x_1,x_2)$ serves as a smooth coordinate system on $M$, so
\beq
\Phi(x)=\psi_{\pm}(x_1,x_2),
\label{3.5.14}
\eeq
with $\psi_{\pm}$ smooth on a disk $D_\delta(0)\subset\RR^2$ and invariant under
rotations.  It is a very special case of results of [Ma] that $\psi_{\pm}$
are smooth in $x_1^2+x_2^2$, so
\beq
\Phi(x)=\gamma_{\pm}(x_1^2+x_2^2),\quad \gamma_{\pm}\in C^\infty([0,\delta)).
\label{3.5.15}
\eeq
This observation applies in particular to $x_3\in C^\infty(M)$, so we have
\eqref{3.5.9}, with $\varphi_{\pm}$ as in \eqref{3.5.10}.
The last item of \eqref{3.5.10},
$\mp \varphi'_{\pm}(0)>0$, follows from $K>0$ at the poles of $M$.  Thus the
analysis \eqref{3.5.11}--\eqref{3.5.13} applies, and we get
\beq
\rho(x_3)=x_1^2+x_2^2\bigr|_M\Rightarrow \rho\in C^\infty([-a,a]),\
\rho'(\pm a)\neq 0,
\label{3.5.16}
\eeq
so, near the poles $(0,0,\pm a)$,
\beq
\Phi(x)=\gamma_{\pm}(\rho(x_3)),
\label{3.5.17}
\eeq
smooth in $x_3\in [-a,a]$.
\qed

We take a further look at an ingredient
in the proof of Proposition \ref{p3.5.1},
namely the following.  Let $D_\delta(0)=\{(x_1,x_2)\in\RR^2:x_1^2+x_2^2<
\delta^2\}$.

\begin{lemma} \label{l3.5.2}
If $\psi\in C^\infty(D_\delta(0))$ is invariant under
rotations, then there exists $\gamma\in C^\infty([0,\delta^2))$ such that
\beq
\psi(x_1,x_2)=\gamma(x_1^2+x_2^2).
\label{3.5.18}
\eeq
\end{lemma}

We present a direct proof of this, not appealing to the general (and rather
deep) work of [Ma].  It is clear that, if $\psi$ is rotationally invariant,
then \eqref{3.5.18} holds with
\beq
\gamma(s)=\psi(s^{1/2},0),\quad s\in [0,\delta^2).
\label{3.5.19}
\eeq
The crux of the matter is to show that such $\gamma$ is $C^\infty$ on $[0,
\delta^2)$, and of course such smoothness is clear except at $s=0$.  To
restate \eqref{3.5.19}, we have
\beq
\gamma(s)=\tilde{\psi}(s^{1/2}),\quad \text{with }\ \tilde{\psi}(t)=
\psi(t,0).
\label{3.5.20}
\eeq
We have
\beq
\tilde{\psi}\in C^\infty((-\delta,\delta)),\quad
\tilde{\psi}(-t)=\tilde{\psi}(t),
\label{3.5.21}
\eeq
and we want to deduce from this that $\gamma$ is smooth at $s=0$.
  
Now \eqref{3.5.21} implies that the formal power series of $\tilde{\psi}$
has the form
\beq
\sum\limits_{k=0}^\infty a_k t^{2k},
\label{3.5.22}
\eeq
with only even powers of $t$ appearing.  Consider the formal power series
\beq
\sum\limits_{k=0}^\infty a_k s^k.
\label{3.5.23}
\eeq
A theorem of Borel guarantees that there exists $\tilde{\gamma}\in
C^\infty((-\delta^2,\delta^2))$ whose formal power series is given by
\eqref{3.5.23}.  Thus $\gamma(t^2)=\psi(t,0)$ and $\tilde{\gamma}(t^2)$ both
have the same formal power series, namely \eqref{3.5.22}.  Thus
\beq
\gamma(t^2)-\tilde{\gamma}(t^2)=u(t),
\label{3.5.24}
\eeq
with
\beq
u\in C^\infty((-\delta,\delta)),\quad u^{(j)}(0)=0,\ \forall\, j.
\label{3.5.25}
\eeq
It then follows from the chain rule that
\beq
v(s)=u(s^{1/2})\Longrightarrow v\in C^\infty([0,\delta^2)).
\label{3.5.26}
\eeq
Since
\beq
\gamma(s)=\tilde{\gamma}(s)+v(s),
\label{3.5.27}
\eeq
this proves the desired smoothness of $\gamma$ at $s=0$.

\section{Stability of stationary solutions}\label{c4}

In this section we examine stability of
stationary zonal solutions of \eqref{1.1},
again assuming $M$ is radially symmetric and has positive Gauss curvature.
First, in \S{\ref{c4s1}}, we look at an Arnold-type approach to stability,
bringing in functionals
\beq
\Cal{H}(u)=\int\limits_M \Bigl\{\frac{1}{2}|u|^2+\varphi(w-\Omega\chi)
+\gamma\xi w\Bigr\}\, dS,
\label{4.1}
\eeq
for various functions $\varphi$ and real constants $\gamma$.  Given a
stationary solution $J\nabla f$ and $w=\Delta f$, we see that if
\beq
\text{$w(\xi)-\Omega\chi(\xi)$ is strictly monotone in $\xi$,}
\label{4.2}
\eeq
then one can find $\varphi$ and $\gamma$ such that $u=J\nabla f$ is a
critical point of \eqref{4.1}, with positive definite second derivative.
Stability in $H^1(M)$ is a consequence.
Note that, for fixed $f$ (hence fixed $w$),
\eqref{4.2} holds for all sufficiently large $\Omega$.

In \S{\ref{c4s2}}, we linearize \eqref{1.1}
about a stationary zonal solution $J\nabla f$.
More precisely, we linearize the associated vorticity equation, obtaining
a linear equation of the form
\beq
\frac{\pa\zeta}{\pa t}=\Gamma\zeta,\quad \Gamma\zeta=-\nabla_{J\nabla f}\zeta
+\nabla_{J\nabla(w-\Omega\chi)}\Delta^{-1}\zeta.
\label{4.3}
\eeq
The symmetry hypothesis on $M$ allows us to write
\beq
\Gamma=\bigoplus\limits_k \Gamma_k,\quad \Gamma_k:V_k\rightarrow V_k,\quad
V_k=\{\zeta\in L^2(M):X_3\zeta=ik\zeta\},
\label{4.4}
\eeq
and deduce that $\Gamma$ has spectrum off the imaginary axis if and only
if some $\Gamma_k\ (k\neq 0)$ has an eigenvalue off the imaginary axis.
In this setting, we derive a version of the Rayleigh criterion, namely,
if $\Gamma$ has an eigenvalue with nonzero real part, then
\beq
w'(\xi)-\Omega\chi'(\xi)\ \text{ must change sign.}
\label{4.5}
\eeq
Note how this interfaces with the criterion \eqref{4.2} for Arnold-type
stability.  We see that the Arnold-type criterion for proving stability and the
Rayleigh-type condition for the lack of proof of linear instability are
almost equivalent.

This is not at all to say that the criterion \eqref{4.2} nails stability.
Just when stability holds and when it fails remains a subtle
question.  The rest of this paper is aimed at formulating some attacks on
this queston.  In \S{\ref{c4s3}} we set things up for some specific calculations,
which will be continued in \S{\ref{c5}}.  At this point, we will want to make use of
classical results on spherical harmonics, so in \S{\ref{c4s3}} and
\S{\ref{c5}} we will specialize to the case $M=S^2$.

In \S{\ref{c4s3}}, we look at \eqref{4.4} with
\beq
\Gamma_k=ikM_k,\quad M_k=M\bigr|_{V_k},\quad
M\zeta=A(x_3)\zeta+B(x_3)\Delta^{-1}\zeta.
\label{4.6}
\eeq
In the setting of \S{\ref{c4s2}},
$A(x_3)=f'(x_3)$ and $B(x_3)=\Omega-w'(x_3)$.
We present some results on Spec $M_k$, particularly when
\beq
A(x_3)=\alpha f'(x_3),\quad B(x_3)=\Omega+\lambda_\nu \alpha f'(x_3).
\label{4.7}
\eeq
These results will have further use in \S{\ref{c5}}.

\subsection{Arnold-type stability results}\label{c4s1}

We use the following variant of the Arnold stability method (cf.~\cite{AK},
pp.~89--94, \cite{MP}, pp.~106--111) for producing stable, stationary solutions to
the 2D Euler equations, in case $M$ is rotationally symmetric, and has positive
Gauss curvature.  Namely, we
look for stable critical points of a functional
\beq
\Cal{H}(u)=\int\limits_M \Bigl\{\frac{1}{2}|u|^2+\varphi(w-\Omega\chi)
+\gamma\xi w\Bigr\}\, dS,
\label{4.1.1}
\eeq
with $w=\rot u$ and $\varphi$ and $\gamma$ tuned to the specific steady
solution $u$.  The functions $\chi$ and $\xi$ are as in \eqref{1.2} and
\eqref{3.3.2}.  See also \eqref{3.4.4} and \eqref{3.4.10}.
Such a functional is independent of $t$ when
applied to a solution $u(t)$ to \eqref{1.1}.  Taking
\beq
u=J\nabla f,\quad \text{so }\ w=\Delta f,
\label{4.1.2}
\eeq
we rewrite \eqref{4.1.1} as
\beq
H(f)=\int\limits_M \Bigl\{\frac{1}{2}|\nabla f|^2+\varphi(\Delta f-\Omega\chi)
+\gamma\xi\Delta f\Bigr\}\, dS.
\label{4.1.3}
\eeq
Then
\beq
\aligned
\pa_sH(f+sg)=\int\limits_M \Bigl\{\langle \nabla f,\nabla g\rangle
+s|\nabla g|^2+\varphi'(\Delta f+s\Delta g-\Omega\chi)\Delta g& \\
+\gamma\xi\Delta g\Bigr\}&\, dS,
\endaligned
\label{4.1.4}
\eeq
so
\beq
\aligned
&\pa_sH(f+sg)\bigr|_{s=0} \\
&=\int\limits_M \Bigl\{\langle\nabla f,\nabla g\rangle
+\varphi'(\Delta f-\Omega\chi)\Delta g+\gamma\xi\Delta g\Bigr\}\, dS \\
&=\int\limits_M \Bigl\{-\Delta f+\Delta \varphi'(\Delta f-\Omega\chi)
+\gamma\Delta\xi\Bigr\}g\, dS.
\endaligned
\label{4.1.5}
\eeq
This is 0 for all $g$ if and only if $f-\varphi'(\Delta f-\Omega\chi)
-\gamma\xi$ is constant, and since the stream function $f$ is determined only
up to an additive constant, we can write
\beq
f=\varphi'(\Delta f-\Omega\chi)+\gamma\xi,
\label{4.1.6}
\eeq
as the condition for $f$ to be a critical point of $H$ in \eqref{4.1.3}.
Note that \eqref{4.1.6} implies that, if
\beq
\nabla f\, \| \, \nabla\xi,
\label{4.1.7}
\eeq
then
\beq
\nabla(\Delta f-\Omega\chi)\, \| \, \nabla f,
\label{4.1.8}
\eeq
hence
\beq
\langle J\nabla f,\nabla(w-\Omega\chi)\rangle=0,
\label{4.1.9}
\eeq
so by Proposition \ref{p3.1.1}, such $f$ produces a stationary solution to
\eqref{1.1}, provided $f$ is a zonal function.
If $f$ is not a zonal function, one would need to
take $\gamma=0$ in \eqref{4.1.1} in order for \eqref{4.1.9} to hold.
(Consequently, the Arnold method apparently produces much weaker stability
results for non-zonal stationary solutions than for zonal stationary solutions.)

To proceed, we apply $\pa_s$ to \eqref{4.1.4} and evaluate at $s=0$, to get
\beq
\pa_s^2H(f+sg)\bigr|_{s=0}=\int\limits_M \Bigl\{|\nabla g|^2
+\varphi''(\Delta f-\Omega\chi)(\Delta g)^2\Bigr\}\, dS.
\label{4.1.10}
\eeq

Now, if we are given a zonal function $f$, we want to find $\varphi$ such
that \eqref{4.1.6} holds, and then check \eqref{4.1.10} to see if this is a
coercive quadratic form in $g$.  Let us write also
\beq
\Delta f=w(\xi),\quad \chi=\chi(\xi).
\label{4.1.11}
\eeq
Then \eqref{4.1.6} takes the form
\beq
f(\xi)=\varphi'(w(\xi)-\Omega\chi(\xi))+\gamma\xi,
\label{4.1.12}
\eeq
or
\beq
\varphi'(w(\xi)-\Omega\chi(\xi))=f(\xi)-\gamma\xi.
\label{4.1.13}
\eeq
Given arbitrary $\gamma\in\RR$, this identity uniquely specifies $\varphi'$,
provided
\beq
w(\xi)-\Omega\chi(\xi)\ \text{ is strictly monotone in }\ \xi,
\label{4.1.14}
\eeq
that is,
\beq
w'(\xi)-\Omega\chi'(\xi)\ \text{ is bounded away from }\ 0.
\label{4.1.15}
\eeq
With $\varphi'$ determined, in turn $\varphi$ is determined, up to an
additive constant, which would not affect the critical points of
\eqref{4.1.3}.  Then, applying $d/d\xi$ to \eqref{4.1.13} yields
\beq
\varphi''(w-\Omega\chi)=\frac{\gamma-f'(\xi)}{\Omega\chi'(\xi)-w'(\xi)}.
\label{4.1.16}
\eeq
Substitution into \eqref{4.1.10} gives
\beq
\pa_s^2H(f+sg)\bigr|_{s=0}=\int\limits_M \Bigl\{|\nabla g|^2+
\frac{\gamma-f'(\xi)}{\Omega\chi'(\xi)-w'(\xi)}(\Delta g)^2\Bigr\}\, dS.
\label{4.1.17}
\eeq
By calculations of \S\S{\ref{c3s4}--\ref{c3s5}},
as long as the Gauss curvature of $M$ is everywhere
positive, both $\chi$ and $\xi$ are smooth, strictly monotonic functions of
$x_3$, with positive $x_3$-derivatives, so
\beq
\chi'(\xi)\ge a>0\ \text{ on }\ M.
\label{4.1.18}
\eeq
As long as the hypothesis \eqref{4.1.14}--\eqref{4.1.15} holds, then either
\beq
\aligned
\Omega\chi'(\xi)-w'(\xi)&\ge b>0,\ \ \text{or } \\
\Omega\chi'(\xi)-w'(\xi)&\le -b<0,
\endaligned
\label{4.1.19}
\eeq
on $M$.  In the first case, we can make
\beq
K(\xi)=\frac{\gamma-f'(\xi)}{\Omega\chi'(\xi)-w'(\xi)}\ge c>0
\label{4.1.20}
\eeq
on $M$ by taking $\gamma>0$ large enough, and in the second case we can arrange
\eqref{4.1.20} by taking $\gamma$ sufficiently negative.  Both cases yield
\beq
\pa_s^2H(f+sg)\bigr|_{s=0}\ge \|\nabla g\|^2+C\|\Delta g\|^2_{L^2},
\label{4.1.21}
\eeq
with $C>0$, for all $g\in H^2(M)$.  This implies stability of $f$ in $H^2(M)$
as a critical point of \eqref{4.1.3}
(recall that $f$ is defined only up to an additive constant),
hence stability of $u$ in $H^1(M)$ as a critical point of
\eqref{4.1.1}.  We summarize.

\begin{theorem} \label{t4.1.1}
Given a smooth $f(\xi)$, $u=J\nabla f$ is a stable
stationary solution to (1.1), in $H^1(M)$, as long as $\Omega$ is such that
\eqref{4.1.14}--\eqref{4.1.15} hold, where $w=\Delta f$.
\end{theorem}

Note that $w=\Delta f$ implies
\beq
w(\xi)=f'(\xi)\Delta\xi+f''(\xi)|\nabla\xi|^2,
\label{4.1.22}
\eeq
if $f=f(\xi)$.

\subsection{Linearization about a stationary solution}\label{c4s2}

Let $M\subset\RR^3$ be a compact surface, rotationally symmetric about the
$x_3$-axis, with positive Gauss curvature,
and let $u=J\nabla f$ be a stationary
solution to \eqref{1.1}.  We derive an equation for the linearization at $u$.
More precisely, we work with the vorticity equation \eqref{2.2.15}, i.e.,
\beq
\frac{\pa w}{\pa t}+\langle J\nabla f,\nabla(w-\Omega\chi)\rangle=0.
\label{4.2.1}
\eeq
Let us set
\beq
f_\ep(t)=f+\ep\eta(t)+\cdots,\quad w_\ep(t)=w+\ep\zeta(t)+\cdots,\quad
\zeta=\Delta\eta.
\label{4.2.2}
\eeq
Inserting these into the analogue of \eqref{4.2.1},
using \eqref{4.2.1} and discarding
higher powers of $\ep$ produces the linearized equation        
\beq
\pa_t\zeta+\langle J\nabla f,\nabla\zeta\rangle+\langle J\nabla\eta,
\nabla(w-\Omega\chi)\rangle=0.
\label{4.2.3}
\eeq
Now
\beq
\aligned
\langle J\nabla\eta,\nabla(w-\Omega\chi)\rangle
&=-\langle\nabla\eta, J\nabla(w-\Omega\chi)\rangle \\
&=-\nabla_{J\nabla(w-\Omega\chi)}\eta.
\endaligned
\label{4.2.4}
\eeq
Also, since $\zeta$ integrates to $0$ on $M$, we can write
\beq
\eta=\Delta^{-1}\zeta,
\label{4.2.5}
\eeq
where, here and below, we define $\Delta^{-1}$ to annihilate constants and to
have range orthogonal to constants.
Then \eqref{4.2.3} becomes the linear equation
\beq
\frac{\pa\zeta}{\pa t}=\Gamma\zeta,
\label{4.2.6}
\eeq
where
\beq
\Gamma\zeta=-\nabla_{J\nabla f}\zeta+\nabla_{J\nabla(w-\Omega\chi)}
\Delta^{-1}\zeta.
\label{4.2.7}
\eeq
The question of linear stability is the question of whether $\Gamma$ generates
a uniformly bounded group of operators on
\beq
L^2_b(M)=\Bigl\{\zeta\in L^2(M):\int\limits_M \zeta\, dS=0\Bigr\}.
\label{4.2.8}
\eeq

Under our hypotheses, we have $\chi=\chi(\xi)$,
with $\xi$ as in \eqref{3.3.2}, i.e.,
$J\nabla\xi=-X_3$.  Let us also assume $f$ is a zonal function, i.e.,
$X_3f=0$, so $f=f(\xi)$.  This also implies $X_3w=0$, hence $w=w(\xi)$.  Then
\beq
\aligned
J\nabla f&=-f'(\xi)X_3, \\
J\nabla(w-\Omega\chi)&=[\Omega\chi'(\xi)-w'(\xi)]X_3,
\endaligned
\label{4.2.9}
\eeq
and \eqref{4.2.7} becomes
\beq
\Gamma\zeta=f'(\xi)X_3\zeta+(\Omega\chi'(\xi)-w'(\xi))X_3\Delta^{-1}\zeta.
\label{4.2.10}
\eeq
In such a case, $\Gamma$ commutes with $X_3$.  hence we can decompose
\beq
L^2_b(M)=\bigoplus\limits_k V_k,
\label{4.2.11}
\eeq
where, for $k\in\ZZ$,
\beq
V_k=\{\zeta\in L^2_b(M):X_3\zeta=ik\zeta\},             
\label{4.2.12}
\eeq
and we have
\beq
\Gamma=\bigoplus\limits_k \Gamma_k,\quad \Gamma_k:V_k\rightarrow V_k,
\label{4.2.13}
\eeq
where
\beq
\Gamma_k\zeta=ik\bigl[f'(\xi)\zeta+(\Omega\chi'(\xi)-w'(\xi))\Delta^{-1}\zeta].
\label{4.2.14}
\eeq
Note that
\beq
\Delta^{-1}:V_k\longrightarrow V_k\ \text{ is compact,}
\label{4.2.15}
\eeq
for each $k$, so each $\Gamma_k$ is a compact perturbation of a bounded,
skew-adjoint operator on $V_k$.  In light of this, basic analytic Fredholm
theory yields the following.

\begin{proposition} \label{p4.2.1}
For each $k$,
\beq
\Spec \Gamma_k\subset ik\Sigma\cup S_k,
\label{4.2.16}
\eeq
where
\beq
\Sigma=\{f'(\lambda):\alpha_0\le\lambda\le\alpha_1\},\quad
\alpha_0=\min\limits_M\, \xi,\ \alpha_1=\max\limits_M \, \xi,
\label{4.2.16A}
\eeq
and $S_k$ is a countable set of points in $\CC$ whose accumulation points all
must lie in $ik\Sigma$.  Each $\mu\in S_k$ is an eigenvalue of $\Gamma_k$,
and the associated generalized eigenspace is finite dimensional.
\end{proposition}

In fact, for each $\mu\in \CC\setminus ik\Sigma$, $\Gamma_k-\mu I$ is a bounded
operator on $V_k$ that is Fredholm of index $0$, and it is clearly invertible
for $|\mu|>\|\Gamma_k\|$.

\begin{corollary} \label{c4.2.2}
Assume $\Gamma$ has the form \eqref{4.2.10}.  If $\Spec
\Gamma$ is not contained in the imaginary axis, then some $\Gamma_k$ has an
eigenvalue with nonzero real part.
\end{corollary}

Now having $\Spec \Gamma\subset i\RR$ would not guarantee that $\Gamma$
generates a bounded group of operators on $L^2_b(M)$, but not having this
inclusion definitely guarantees that the associated group of operators
is not uniformly bounded.  Thus Corollary \ref{c4.2.2} points to an approach
to finding cases that are linearly unstable.

Actually establishing such cases of linear instability is not so straightforward.
We proceed to derive some necessary conditions for such linear instability
to hold, i.e., for some $\Gamma_k$ to have an eigenvalue with nonzero real
part.

Of course, $\Gamma_0=0$.  Suppose $k\neq 0$ and
$\Gamma_k$ has an eigenvalue $\mu=ik\beta,\ \beta\notin\RR$.
Then there exists a nonzero $\zeta\in V_k$ such that
\beq
(f'(\xi)-\beta)\zeta=-(\Omega\chi'(\xi)-w'(\xi))\Delta^{-1}\zeta,
\label{4.2.17}
\eeq
hence
\beq
\Delta \eta=\frac{w'(\xi)-\Omega\chi'(\xi)}{f'(\xi)-\beta} \eta,
\label{4.2.18}
\eeq
where $\eta=\Delta^{-1}\zeta$.  Note that if $\beta\notin\RR$, the
denominator on the right side of \eqref{4.2.18} is nowhere vanishing.
In \eqref{4.2.17}--\eqref{4.2.18}, $\zeta$ and $\eta$ would not be real
valued.  Taking the inner product of both
sides of \eqref{4.2.18} with $\eta$ yields
\beq
\aligned
(\Delta\eta,\eta)
&=\int\limits_{S^2}\frac{w'(\xi)-\Omega\chi'(\xi)}{f'(\xi)-\beta} |\eta|^2\,
dS \\
&=\int\limits_{S^2}\frac{w'(\xi)-\Omega\chi'(\xi)}{|f'(\xi)-\beta|^2}           
[f'(\xi)-\overline{\beta}]\,|\eta|^2\, dS.
\endaligned
\label{4.2.19}
\eeq
Now $(\Delta\eta,\eta)$ is real and negative, but $\Im \overline{\beta}\neq
0$.  Hence taking the imaginary part of \eqref{4.2.19} yields
\beq
\int\limits_{S^2} \frac{w'(\xi)-\Omega\chi'(\xi)}{|f'(\xi)-\beta|^2}
|\eta|^2\,dS=0.
\label{4.2.20}
\eeq
Using this in \eqref{4.2.19} gives
\beq
(\Delta\eta,\eta)=\int\limits_{S^2}
\frac{w'(\xi)-\Omega\chi'(\xi)}{|f'(\xi)-\beta|^2}
[f'(\xi)-K]\,|\eta|^2\, dS<0,\ \ \forall\, K\in\RR.
\label{4.2.21}
\eeq
We have \eqref{4.2.20} and \eqref{4.2.21}
as necessary conditions for $\Gamma_k$ to have
an eigenvalue with nonzero real part, with associated eigenfunction $\zeta
=\Delta \eta,\ \eta\in V_k$.  These results in turn imply the following.

\begin{proposition} \label{p4.2.3}
If $\Gamma$ has an eigenvalue with nonzero real
part, then
\beq
w'(s)-\Omega\chi'(s)\ \text{ must change sign in }\ s\in (\alpha_0,\alpha_1),
\label{4.2.22}
\eeq                                                                 
with $\alpha_j$ as in \eqref{4.2.16A}, and
\beq
\forall\, K\in\RR,\ \ \exists\, s\in (\alpha_0,\alpha_1)\ \text{ such that }\
(w'(s)-\Omega\chi'(s))(f'(s)-K)<0.
\label{4.2.23}
\eeq
\end{proposition}

In the setting of planar flows (and with $\Omega=0$), \eqref{4.2.22} is
known as the ``Rayleigh criterion'' for linear instability, and
\eqref{4.2.23} is called the ``Fjortoft criterion.''  See \cite{MP},
pp.~122--123.

Proposition \ref{p4.2.3} is close to Theorem \ref{t4.1.1} in the following
sense.  By Theorem \ref{t4.1.1}, if
\beq
w'(s)-\Omega\chi'(s)\neq 0\ \text{ for all }\ s\in [\alpha_0,\alpha_1],
\label{4.2.24}
\eeq
then the associated stationary solution
$u=J\nabla f$ to \eqref{1.1} is stable, in the sense of \S{\ref{c4s1}}.
Condition \eqref{4.2.22} is a little stronger than the assertion that
\eqref{4.2.24} fails.  Thus, in some sense, the first part of Proposition
\ref{p4.2.3} is almost a corollary of Theorem \ref{t4.1.1}.

\subsection{Further results on linearization}\label{c4s3}

Here we produce some results complementary to those of \S{\ref{c4s2}}.  We
consider operators $\Gamma$ of a more general nature than those in
\S{\ref{c4s2}}, as indicated in \eqref{4.3.4} below.  However, we specialize
from more general surfaces of rotation to the standard sphere $S^2$, in order
to make some explicit computations using spherical harmonics.

To proceed, we investigate matters related to
whether the operator $\Gamma$ generates a
uniformly bounded group on $L^2_b(S^2)=\{f\in L^2(S^2):\int_{S^2}f\, dS=0\}$,
when $\Gamma$ has the following structure:
\beq
\Gamma=\bigoplus\limits_k \Gamma_k,\quad \Gamma_k:V_k\rightarrow V_k,\quad
V_k=\{f\in L^2_b(S^2):X_3f=ikf\},
\label{4.3.1}
\eeq
where $X_3$ is the vector field generating $2\pi$-periodic rotation about the
$x_3$-axis.  We assume
\beq
\Gamma_k=ikM_k,\quad M_k=M\bigr|_{V_k},
\label{4.3.2}
\eeq
where
\beq
M\zeta=A(x_3)\zeta+B(x_3)\Delta^{-1}\zeta.
\label{4.3.3}
\eeq
We assume $A$ and $B$ are smooth and real valued.  In studies of linear
stability of stationary, zonal Euler flows on the rotating sphere, such an
operator arises with
\beq
A(x_3)=f'(x_3),\quad B(x_3)=\Omega-w'(x_3),
\label{4.3.4}
\eeq
with $w=\Delta f=\rot u$, $u$ a steady zonal solution to the Euler equation.

The question we examine is whether $\Spec \Gamma$ is contained in the
imaginary axis.  In view of \eqref{4.3.2}, this is equivalent to the question of
whether $\Spec M_k$ is contained in the real axis, for each $k\neq 0$.
Basic Fredholm theory gives the following.  (Compare Proposition \ref{p4.2.1}.)

\begin{proposition} \label{p4.3.1}
For each $k\neq 0$,
$$
\Spec M_k\subset \Sigma\cup S_k,
$$
where $\Sigma=\{A(x_3):-1\le x_3\le 1\}$ and $S_k$ is a countable set of
points in $\CC$ whose accumulation points all must lie in $\Sigma$.  Each
$\lambda\in S_k$ is an eigenvalue of $M_k$, and the associated generalized
eigenspace is finite dimensional.
\end{proposition}

In fact, for each $\lambda\in\CC\setminus\Sigma$, $M_k-\lambda I$ is a bounded
operator on $V_k$ that is Fredholm, of index $0$, and it is clearly invertible
for $|\lambda|>\|M_k\|$.  The next result is a cousin to Corollary
\ref{c4.2.2}.

\begin{corollary} \label{c4.3.2}
Assume $M$ has the form \eqref{4.3.3}.  If $\Spec \Gamma$
is not contained in the imaginary axis, then some $M_k\ (k\neq 0)$ has an
eigenvalue that is not real.
\end{corollary}

\sk
{\sc Remark.} If $\lambda$ is a non-real eigenvalue of $M_k$, then
$\overline{\lambda}$ is an eigenvalue of both $M_k$ and $M_{-k}$.

$\text{}$

Actually establishing such cases of linear instability is not so straightforward.
We proceed to derive some necessary conditions for such linear instability to
hold, i.e., for some $M_k$ to have a non-real eigenvalue.

If $\lambda\notin\RR$ is an eigenvalue of $M_k$, then there is a nonzero
$\zeta\in V_k$ such that
\beq
(A(x_3)-\lambda)\zeta=-B(x_3)\Delta^{-1}\zeta.
\label{4.3.5}
\eeq
We can take $\eta=\Delta^{-1}\zeta\in V_k$ and write this as
\beq
\Delta\eta=\frac{B(x_3)}{\lambda-A(x_3)} \eta.
\label{4.3.6}
\eeq
Note that if $\lambda\notin\RR$, then, since $A$ is real valued, the denominator
on the right side of \eqref{4.3.6} is nowhere vanishing.
(Note also that $\eta$ is not real valued.)
We take the inner product of both sides of \eqref{4.3.6} with $\eta$,
to get
\beq
\aligned
(\Delta\eta,\eta)&=\int\limits_{S^2}\frac{B(x_3)}{\lambda-A(x_3)}
|\eta|^2\, dS \\
&=\int\limits_{S^2} \frac{B(x_3)}{|\lambda-A(x_3)|^2}
(\overline{\lambda}-A(x_3))|\eta|^2\, dS.
\endaligned
\label{4.3.6A}
\eeq
Now $(\Delta\eta,\eta)$ is real and negative.  Hence the imaginary part
of the last integral is zero.  If $\lambda\notin\RR$, this forces
\beq
\int\limits_{S^2} \frac{B(x_3)}{|\lambda-A(x_3)|^2}|\eta|^2\, dS=0.
\label{4.3.7}
\eeq
Given this, we can then deduce from (4.3.6A) that
\beq
(\Delta\eta,\eta)=\int\limits_{S^2}\frac{B(x_3)}{|\lambda-A(x_3)|^2}
(K-A(x_3))|\eta|^2\, dS<0,\quad \forall\, K\in\RR.
\label{4.3.8}
\eeq

We have \eqref{4.3.7} and \eqref{4.3.8} as necessary conditions for $M_k$ to
have an eigenvalue $\lambda\notin\RR$, with associated eigenfunction $\zeta=
\Delta\eta,\ \eta\in V_k$.  These results imply the following.
(Compare Proposition \ref{p4.2.3}.)

\begin{proposition} \label{p4.3.3}
If $\Gamma$ has a non-imaginary eigenvalue, then
\beq
\text{$B(s)$ must change sign in } s\in (-1,1),
\label{4.3.9}
\eeq
and
\beq
\forall\, K\in\RR,\ \exists\, s\in(-1,1)\ \text{such that}\ B(s)(K-A(s))<0.
\label{4.3.10}
\eeq
\end{proposition}

Condition \eqref{4.3.9} is a version of the ``Rayleigh criterion'' and
\eqref{4.3.10} a version of the ``Fjortoft criterion'' for linear instability.
Compare the remarks after Proposition \ref{p4.2.3}.

Regarding the relation between \eqref{4.3.9} and \eqref{4.3.10},
we mention that there is at least one situation
where the Rayleigh criterion \eqref{4.3.9} holds but the Fjortoft
condition \eqref{4.3.10} fails, namely when $A(x_3)=A$ is constant.
Then \eqref{4.3.10} fails for $K=A$,
but \eqref{4.3.9} holds for many choices of $B(x_3)$.  This result
is equivalent to the statement that
\beq
B(x_3)\Delta^{-1}\ \text{ has real spectrum on each }\ V_k,
\label{4.3.11}
\eeq
for $k\neq 0$.  This fact might seem nontrivial, since $B(x_3)\Delta^{-1}$
is not self adjoint (if $B(x_3)$ is not constant), but this operator acts
on Sobolev scales, and in this framework the operator is similar to, and has
the same spectrum as
\beq
-(-\Delta)^{-1/2}B(x_3)(-\Delta)^{-1/2},                              
\label{4.3.12}
\eeq
which is self adjoint.

Regarding the reverse implication, we have:

\begin{proposition} \label{p4.3.4}
Assume $A$ and $B$ are continuous and real valued on
$[-1,1]$.  Then \eqref{4.3.10} $\Rightarrow$ \eqref{4.3.9}.
\end{proposition}

\demo
There exist $K_1$ and $K_2$ such that $K_1-A(s)>0$ for all
$s\in [-1,1]$ and $K_2-A(s)<0$ for all $s\in [-1,1]$.
Applying \eqref{4.3.10} to $K=K_1$ yields $s_1\in (-1,1)$
such that $B(s_1)<0$ and applying \eqref{4.3.10} to
$K=K_2$ yields $s_2\in (-1,1)$ such that $B(s_2)>0$.
\qed

Here is another case where \eqref{4.3.9} does not imply \eqref{4.3.10}.
Namely, $A(s)=-B(s)$, where \eqref{4.3.10} fails for $K=0$.

It seems not so easy to give examples where
\eqref{4.3.9} holds but \eqref{4.3.10} fails when
$A(s)$ and $B(s)$ have the form \eqref{4.3.4}, with $f$ and $w$ zonal functions
related by $w=\Delta f$.  Suppose, for example, that $f$ is a zonal
eigenfunction of $\Delta$,
\beq
\Delta f=-\lambda_\nu f,\quad \text{so }\ w=-\lambda_\nu f\quad (\lambda_\nu>0).
\label{4.3.13}
\eeq
Then
\beq
\aligned
B(s)(K-A(s))&=(\Omega+\lambda_\nu f'(s))(K-f'(s)) \\
&=-\lambda_\nu \Bigl(-\frac{\Omega}{\lambda_\nu}-f'(s)\Bigr)(K-f'(s)).
\endaligned
\label{4.3.14}
\eeq
To pick $K\in\RR$ violating \eqref{4.3.10}, we need the two factors above to
have the same zeros, to avoid the product changing sign.  This tends to
force $K=-\Omega/\lambda_\nu$.  But then $B(s)(K-A(s))=-\lambda_\nu(K-f'(s))^2$,
which is $<0$ on most of $(-1,1)$.  So this approach fails to produce an
example where \eqref{4.3.9} holds but \eqref{4.3.10} fails.

Having the Fjortoft condition hold along with the Rayleigh condition is certainly
an acceptable state of affairs, and it is of interest to pursue the use of
such $f$ as in \eqref{4.3.13}.  We will find it useful to generalize a little, and
consider the situation
\beq
\Delta f=-\lambda_\nu f,\quad w=-\mu f\quad (\mu>0).
\label{4.3.15}
\eeq
Note that
\beq
\aligned
\Spec (-\Delta)&=\{\lambda_j=j(j+1):j=0,1,2,3,\dots\}, \\
\Spec (-\Delta)\bigr|_{V_k}&=\{\lambda_j:j\ge |k|\}.
\endaligned
\label{4.3.16}
\eeq

Let us take
\beq
f(x_3)=\tilde{\alpha}P_\nu(x_3),\quad \tilde{\alpha}>0,
\label{4.3.17}
\eeq
where $P_\nu$ are Legendre polynomials, given by
$$
P_k(s)=\frac{1}{2^k k!}\Bigl(\frac{d}{ds}\Bigr)^k(s^2-1)^k,
$$
for example,
\beq
\gathered
P_0(s)=1,\quad P_1(s)=s,\quad P_2(s)=\frac{1}{2}(3s^2-1), \\
P_3(s)=\frac{1}{2}(5s^3-3s),\quad P_4(s)=\frac{1}{8}(35s^4-30s^2+3).
\endgathered
\label{4.3.18}
\eeq
Taking $f=\tilde{\alpha}P_0$ produces a trivial flow.  Taking $f=\tilde{\alpha}
P_1$ gives $f'(s)=\tilde{\alpha}$, hence $w'(s)=-6\mu\tilde{\alpha}$, constant,
so $B(s)=\Omega-w'(s)$ does not satisfy \eqref{4.3.9}.
The first choice that might lead to linear instability is
\beq
f(x_3)=\tilde{\alpha}P_2(x_3),
\label{4.3.19}
\eeq
giving
\beq
f'(x_3)=\alpha x_3,\quad w'(x_3)=-\mu\alpha x_3
\label{4.3.20}                                            
\eeq
(with $\alpha=3\tilde{\alpha}$), hence
\beq
A(x_3)=\alpha x_3,\quad B(x_3)=\Omega+\mu\alpha x_3.
\label{4.3.21}
\eeq
Then, given $\Omega\ge 0$, the Rayleigh condition \eqref{4.3.9} holds if and
only if
\beq
0\le\Omega<\mu\alpha.
\label{4.3.22}
\eeq
For $\Omega>\mu\alpha$, the Arnold stability criterion applies.  Hence, by
a limiting argument, for $\Omega=\mu\alpha$, $\Gamma$ will have no
non-imaginary eigenvalues.

When \eqref{4.3.22} holds, we might find that $\Gamma$ does have some non-imaginary
eigenvalues.  That is, $M_k$ might have some non-real eigenvalues, for some
$k\neq 0$.  Let us take a closer look at this issue when $\Omega=0$.  In
such a case,
\beq
\aligned
M_k\zeta&=\alpha x_3+\mu\alpha x_3\Delta^{-1}\zeta \\
&=\alpha x_3(I+\mu\Delta^{-1})\zeta,
\endaligned
\label{4.3.23}
\eeq
for $\zeta\in V_k$.  Recall we are assuming $\mu>0$.  Now, by \eqref{4.3.16},
\beq
\aligned
0<\mu<\lambda_k
&\Longrightarrow I+\mu\Delta^{-1}\ \text{ is positive definite on $V_k$} \\
&\Longrightarrow \Spec M_k=\alpha \Spec (I+\mu\Delta^{-1})^{1/2} x_3
(I+\mu\Delta^{-1})^{1/2}\bigr|_{V_k} \\
&\Longrightarrow \Spec M_k\subset\RR.
\endaligned
\label{4.3.24}
\eeq
A limiting argument gives the last conclusion for $\mu=\lambda_k$.
We record the conclusion.

\begin{proposition} \label{p4.3.5}
In case $A$ and $B$ are given by \eqref{4.3.21} and
$\Omega=0$,
\beq
0<\mu\le\lambda_k\Longrightarrow \Spec M_k\subset\RR.
\label{4.3.25}
\eeq
\end{proposition}

Now let us specialize to the case relevant for Euler flow.  That is to say,
we take $\mu=\lambda_2=6$ in \eqref{4.3.20}--\eqref{4.3.21}:
\beq
A(x_3)=\alpha x_3,\quad B(x_3)=\Omega+\lambda_2\alpha x_3.
\label{4.3.26}
\eeq

\begin{corollary} \label{c4.3.6}
In case $A$ and $B$ are given by \eqref{4.3.26} and
$\Omega=0$,
\beq
k\ge 2\Longrightarrow \Spec M_k\subset\RR,
\label{4.3.27}
\eeq
and ditto for $M_{-k}$.
\end{corollary}

Thus, if linear instability arises in this situation, the only possibility
is that
\beq
\Spec M_1\ \text{ is not contained in }\ \RR,
\label{4.3.28}
\eeq
and ditto for $M_{-1}$.
It is therefore of great interest to investigate whether \eqref{4.3.28} holds.

Let us extend our considerations to nonzero $\Omega$, in the setting of
\eqref{4.3.26}.  Then we have
\beq
\aligned
(M_k+\lambda_2^{-1}\Omega I)\zeta
&=(\alpha x_3+\lambda_2^{-1}\Omega)\zeta+(\Omega+\lambda_2\alpha x_3)
\Delta^{-1}\zeta \\
&=(\alpha x_3+\lambda_2^{-1}\Omega)[I+\lambda_2\Delta^{-1}]\zeta,
\endaligned
\label{4.3.29}
\eeq
for $\zeta\in V_k$.  Again, we see that $M_k$ has real spectrum if $k\ge 2$,
so again our search for non-real eigenvalues of $M_k$ is reduced to
investigating whether \eqref{4.3.28} holds.

Keep in  mind that Arnold stability holds for $\Omega>\lambda_2\alpha$, in
this situation.  Thus we are looking at when \eqref{4.3.28} holds, given
\beq
0\le\Omega<\lambda_2\alpha.
\label{4.3.30}
\eeq

\sk
{\sc Note.} For the purpose of this analysis, there is no loss of generality
in taking $\alpha=1$.

$\text{}$

We next generalize the setting \eqref{4.3.26}, along the lines of
\eqref{4.3.13}.  Thus, in place of \eqref{4.3.26}, we have
\beq
A(x_3)=\alpha f'(x_3),\quad B(x_3)=\Omega+\lambda_\nu \alpha f'(x_3).
\label{4.3.31}
\eeq
Now, in place of \eqref{4.3.29}, we have
\beq
\aligned
(M_k+\lambda_\nu^{-1}\Omega I)\zeta
&=(\alpha f'(x_3)+\lambda_\nu^{-1}\Omega)\zeta
+(\Omega+\lambda_\nu\alpha f'(x_3))\Delta^{-1}\zeta \\
&=(\alpha f'(x_3)+\lambda_\nu^{-1}\Omega)[I+\lambda_\nu \Delta^{-1}]\zeta,
\endaligned
\label{4.3.32}
\eeq
for $\zeta\in V_k$.  This is a composition
\beq
\lambda_\nu^{-1}B(x_3)(I+\lambda_\nu\Delta^{-1}),
\label{4.3.33}
\eeq
and this operator has real spectrum as long as either factor, $B(x_3)$ or
$I+\lambda_\nu\Delta^{-1}$ is positive, as an operator on $V_k$.
If $\Omega$ is such that $B(x_k)$ changes sign, the operator still has real
spectrum as long as $I+\lambda_\nu\Delta^{-1}$ is positive on $V_k$, i.e.,
as long as $\lambda_\nu\le\lambda_{|k|}$.  This produces the following variant
of Proposition \ref{p4.3.5}.

\begin{proposition} \label{p4.3.7}
In case $A$ and $B$ are given by \eqref{4.3.31}, then
\beq
\Spec M_k\subset\RR
\label{4.3.34}
\eeq
provided that either $\Omega+\lambda_\nu\alpha f'(x_3)$ does not change sign
or $\lambda_\nu\le \lambda_{|k|}$.
\end{proposition}

\section{Appendix by Jeremy Marzuola and Michael Taylor:
Matrix approach and numerical study of linear instability}\label{c5}

In \S{\ref{c4}} we saw that the sufficient condition \eqref{4.2} for stability in
$H^1(M)$ of a stationary zonal solution $u=J\nabla f$ and the necessary
condition \eqref{4.5} for the existence of non-imaginary spectrum of the
linearized operator $\Gamma$ in \eqref{4.3} are almost perfectly complementary.
Nevertheless, as we will see here, the spectrum of $\Gamma$ might be confined
to the imaginary axis even when \eqref{4.5} fails.
Equivalently, the operators $M_k:V_k\rightarrow V_k$ in \eqref{4.6}
might all have real spectrum, even in cases where \eqref{4.5} fails.
Here we specialize to $M=S^2$ and make some calculations in cases
\beq
f(x)=cP_\nu(x_3),\quad \nu=2,3,4.
\label{5.1}
\eeq
The operator $M_k$ takes the form
\beq
M_k\zeta=cP'_\nu(x_3)\zeta
+(\Omega+\lambda_\nu cP'_\nu(x_3))\Delta^{-1}\zeta,\quad
\zeta\in V_k.
\label{5.2}
\eeq
The spaces $V_k$ have orthogonal bases
\beq
\{e^{ik\psi} P^k_\ell(x_3):\ell\ge |k|\},
\label{5.3}
\eeq
which can be normalized to produce orthonormal bases.  Classical identities
for spherical harmonics lead to representations of $M_k$ as infinite
matrices.  We carry out these calculations for $\nu=2$ in \S{\ref{c5s1}} and
for $\nu=3$ in \S{\ref{c5s2}}.

For short, we sometimes refer to the matrices associated to $M_k$ in
\eqref{5.2} as $P_\nu(V_k)$ models.

In \S{\ref{c5s1}} we use the matrix representation
of $M_1$ (for $\nu=2$) to prove that,
for all $\Omega\ge 0$, $M_1$ has only real spectrum.  (That $M_k$ has only real
spectrum for $|k|\ge 2$ in this situation follows from Corollary \ref{c4.3.6}.)
By contrast, the Rayleigh-type condition \eqref{4.5} guarantees $M_1$ has only
real spectrum provided $\Omega>\lambda_2=6$, but it does not apply to
$\Omega\in [0,6)$.
This extra constraint on the spectrum of $M_1$ for such small $\Omega$
was first suggested to the authors by output from a Matlab program.
Having seen the output, we were able to prove that such a constraint holds.
We also show that $M_1$ has a generalized $0$-eigenvector
at $\Omega=0$, giving rise to a weak linear instability.

In \S{\ref{c5s2}} we work out the infinite matrix representations of $M_1$
and $M_2$ (for $\nu=3$),
acting on $V_1$ and $V_2$, respectively.  (In this case, Proposition
\ref{p4.3.7} implies that $M_k$ has only real spectrum for $|k|\ge 3$.)
The Rayleigh-type condition \eqref{4.5} guarantees that $M_1$ and $M_2$ have
only real spectrum provided $\Omega>(4/5)\lambda_3=48/5$.  The analysis of
$M_1$ and $M_2$ is more difficult than that of $M_1$ in \S{\ref{c5s1}}.
At this point, we have numerical results on truncations
of these matrices that indicate linear stability for substantially smaller
values of $\Omega$ than $48/5$.

These numerical results are discussed in \S{\ref{c5s3}}.
There we take $N\times N$ matrix truncations $M_k^N$ of the operators $M_k$,
arising in \eqref{5.2}, for $\nu=3,4,\ k<\nu$.  After some discussion about
stabilization of the non-real spectrum of such matrices for moderately large
$N$, we take $N=400$.  We use Matlab to find the non-real eigenvalues and
graph their imaginary parts, as functions of $\Omega$.  These graphs indicate
linear stability for $\Omega$ somewhat less restricted than what the Arnold-type
stability analysis of \S{\ref{c4s1}} requires.  We also see numerical evidence
of how stability might not be simply a monotone function of ${\Omega}$,
for $\Omega>0$.  Taken together with the rigorous results we have
established through \S{\ref{c5s1}}, these numerical results suggest much
interesting work for the future.

\subsection{Matrix analysis for $f(x)=cP_2(x_3)$}\label{c5s1}

Here we pursue the question of when \eqref{4.3.28} holds.
We recall the setting.
\beq
M_k=M\bigr|_{V_k},\quad V_k=\{f\in L^2_b(S^2):X_3f=ikf\},
\label{5.1.1}
\eeq
and
\beq
M\zeta=A(x_3)\zeta+B(x_3)\Delta^{-1}\zeta.
\label{5.1.2}
\eeq
We take
\beq
A(x_3)=x_3,\quad B(x_3)=\Omega+\lambda_2x_3,\quad \lambda_2=6,
\label{5.1.3}
\eeq
and ask the following.

\sk
{\bf Question.} For what values of $\Omega$ does
\beq
\text{$M_1$ have a non-real eigenvalue?}
\label{5.1.4}
\eeq

$\text{}$

We assume $\Omega\ge 0$.  As we have seen, the ``Rayleigh criterion''
produces
\beq
0\le\Omega<\lambda_2
\label{5.1.5}
\eeq
as a {\it necessary} condition for \eqref{5.1.4} to hold.  We want to 
see how close \eqref{5.1.5} is to being sufficient.  In the context of
\eqref{5.1.3}, it will turn out to be far from sufficient.

To investigate  this, it is convenient to represent $M_1$ as an infinite
matrix.  An orthogonal basis of $V_1$ is given by
\beq
\tilde{\zeta}_\ell=e^{i\psi}P^1_\ell(x_3),\quad \ell\ge 1.
\label{5.1.6}
\eeq
Here $P^1_\ell$ is an associated Legendre function given in (7.12.5)
of \cite{Leb} as
\beq
P^1_\ell(t)=-(1-t^2)^{1/2} P'_\ell(t).
\label{5.1.7}
\eeq
We mention that
\beq
\Delta\tilde{\zeta}_\ell=-\lambda_\ell \tilde{\zeta}_\ell,\quad
\lambda_\ell=\ell(\ell+1).
\label{5.1.8}
\eeq
One has (\cite{Leb}, p.~201, \#10)
\beq
\int_{-1}^1 P^1_\ell(t)^2\, dt=\frac{2}{2\ell+1}\, \frac{(\ell+1)!}{(\ell-1)!}
=\frac{2\ell(\ell+1)}{2\ell+1}.
\label{5.1.9}
\eeq
Hence, up to a constant, which we can ignore, an orthonormal basis of $V_1$
is given by
\beq
\zeta_\ell=\sqrt{\frac{2\ell+1}{2\ell(\ell+1)}}\, e^{i\psi}(1-x_3^2)^{1/2}
P'_\ell(x_3),\quad \ell\ge 1.
\label{5.1.10}
\eeq
The operator $M_1$ is given by
\beq
\aligned
M_1\zeta_\ell&=A(x_3)\zeta_\ell+B(x_3)\Delta^{-1}\zeta_\ell \\
&=A(x_3)\zeta_\ell-\frac{1}{\lambda_\ell}B(x_3)\zeta_\ell \\
&=x_3\zeta_\ell-\frac{1}{\lambda_\ell}(\Omega+\lambda_2x_3)\zeta_\ell.
\endaligned
\label{5.1.11}
\eeq

To proceed, we need to write $x_3\zeta_\ell$ as a linear combination of
$\{\zeta_j\}$.  To do this, we use (7.8.4) and (7.8.2) of \cite{Leb},
\beq
\aligned
tP'_\ell(t)&=P'_{\ell-1}(t)+\ell P_\ell(t), \\
(2\ell+1)P_\ell(t)&=P'_{\ell+1}(t)-P'_{\ell-1}(t),
\endaligned
\label{5.1.12}
\eeq
which combine to give
\beq
tP'_\ell(t)=\frac{\ell+1}{2\ell+1}P'_{\ell-1}(t)
+\frac{\ell}{2\ell+1}P'_{\ell+1}(t).
\label{5.1.13}
\eeq
Plugging this into \eqref{5.1.10} yields
\beq
x_3\zeta_\ell
=\sqrt{\frac{(\ell+1)(\ell-1)}{(2\ell+1)(2\ell-1)}}\zeta_{\ell-1}
+\sqrt{\frac{\ell(\ell+2)}{(2\ell+1)(2\ell+3)}} \zeta_{\ell+1}.
\label{5.1.14}
\eeq
We use the obvious convention that $\zeta_0=0$.  It is illuminating to
write this as
\beq
x_3\zeta_\ell=a_\ell\zeta_{\ell-1}+a_{\ell+1}\zeta_{\ell+1},\quad
a_\ell=\sqrt{\frac{(\ell+1)(\ell-1)}{(2\ell+1)(2\ell-1)}}.
\label{5.1.15}
\eeq
Thus the matrix representation of $\zeta\mapsto x_3\zeta$ on $V_1$ has the
$3\times 3$ truncation
\beq
A^{(3)}=\begin{pmatrix} 0 & a_2 & {}\\ a_2 & 0 & a_3\\ {} & a_3 & 0\end{pmatrix}
=\begin{pmatrix} 0 & \sqrt{1/5} & {}\\ \sqrt{1/5} & 0 & \sqrt{8/35}\\
{} & \sqrt{8/35} & 0\end{pmatrix} .
\label{5.1.16}
\eeq
Returning to \eqref{5.1.11}, we have
\beq
\aligned
M_1\zeta_\ell
&=a_\ell\zeta_{\ell-1}+a_{\ell+1}\zeta_{\ell+1}
-\frac{\lambda_2}{\lambda_\ell}(a_\ell\zeta_{\ell-1}+a_{\ell+1}\zeta_{\ell+1})
-\frac{\Omega}{\lambda_\ell}\zeta_\ell \\
&=a_\ell\zeta_{\ell-1}+a_{\ell+1}\zeta_{\ell+1}-\frac{\lambda_2}{\lambda_\ell}
\Bigl(a_\ell\zeta_{\ell-1}+a_{\ell+1}\zeta_{\ell+1}
+\frac{\Omega}{\lambda_2}\zeta_\ell\Bigr).
\endaligned
\label{5.1.17}
\eeq
Recall that
\beq
\lambda_\ell=\ell(\ell+1),\quad \lambda_2=6,
\label{5.1.18}
\eeq
and $\ell$ runs over $\{1,2,3,\dots\}$ in \eqref{5.1.17}.  
In particular, the $3\times 3$ truncation of $M_1$ is
\beq
M_1^3=A^{(3)}-R^{(3)},
\label{5.1.19}
\eeq
with $A^{(3)}$ as in \eqref{5.1.16}, and
\beq
R^{(3)}=\lambda_2 \begin{pmatrix} 0 & a_2/\lambda_2 & {}\\
a_2/\lambda_1 & 0 & a_3/\lambda_3\\ {} & a_3/\lambda_2 & 0\end{pmatrix}
+\Omega \begin{pmatrix} 1/\lambda_1 & {} & {}\\ {} & 1/\lambda_2 & {}\\
{} & {} & 1/\lambda_3\end{pmatrix} .
\label{5.1.20}
\eeq

As it turns out, Matlab programs strongly indicate that $M_1$ has only
real spectrum, even at $\Omega=0$.  Stimulated by such programs, we have
managed to verify the results they suggest, and prove the following two
propositions.

\begin{proposition} \label{p5.1.1}
Take $\Omega=0$, and let $a_\ell$ be given by
\eqref{5.1.15}, $\lambda_\ell$ by \eqref{5.1.8}.
Then $M_1$ has only real spectrum on $V_1$.
\end{proposition}

\demo
In this case,
\beq
\aligned
&M_1\zeta_1=-2x_3\zeta_1,\quad M_1\zeta_2=0, \\
&M_1\zeta_\ell=\Bigl(1-\frac{\lambda_2}{\lambda_\ell}\Bigr)x_3\zeta_\ell,
\quad \text{for }\ \ell\ge 3.
\endaligned
\label{5.1.21}
\eeq
The formula \eqref{5.1.15} for $x_3\zeta_\ell$ gives
\beq
\text{Range }M_1\subset V_{12}=\text{ Span }\{\zeta_\ell:\ell\ge 2\}.
\label{5.1.22}
\eeq
Thus any eigenfunction of $M_1$ must lie in $V_{12}$.  Now
\beq
M_1:V_{12}\longrightarrow V_{12},\quad M_1=x_3(I+\lambda_2\Delta^{-1})
\label{5.1.23}
\eeq
implies
\beq
M_1\bigr|_{V_{12}}\ \text{ has real spectrum,}
\label{5.1.24}
\eeq
by the argument proving Proposition \ref{p4.3.5}.
This proves Proposition \ref{p5.1.1}.
\qed

Notwithstanding Proposition \ref{p5.1.1}, we do have linear instability at
$\Omega=0$. In fact, it follows from \eqref{5.1.21} and \eqref{5.1.15} that
\beq
e^{itM_1}\zeta_1=\zeta_1-\frac{2it}{\sqrt{5}}\zeta_2,
\label{5.1.25}
\eeq
so $\{e^{itM_1}:t\in\RR\}$ is not uniformly bounded on $V_1$.

We next extend Proposition \ref{p5.1.1} to cover the case $\Omega>0$.

\begin{proposition} \label{p5.1.2}
In the setting of Proposition \ref{p5.1.1} (i.e., with
$A(x_3)=x_3,\ B(x_3)=\Omega+\lambda_2x_3$), take $\Omega>0$.  Then $M_1$
has only real spectrum on $V_1$.
\end{proposition}

\demo
In place of \eqref{5.1.21}, we have
\beq
\aligned
M_1\zeta_1&=-2x_3\zeta_1+\lambda_1^{-1}\Omega \zeta_1, \\
M_1\zeta_2&=\lambda_2^{-1}\Omega \zeta_2, \\
M_1\zeta_\ell&=\Bigl(1-\frac{\lambda_2}{\lambda_\ell}\Bigr)x_3\zeta_\ell
+\lambda_\ell^{-1}\Omega\zeta_\ell,\quad \ell\ge 3.
\endaligned
\label{5.1.26}
\eeq
It is convenient to rewrite this as
\beq
\aligned
(M_1-\lambda_1^{-1}\Omega)\zeta_1&=-2x_3\zeta_1, \\
(M_1-\lambda_1^{-1}\Omega)\zeta_2
&=(\lambda_2^{-1}-\lambda_1^{-1})\Omega\zeta_2, \\
(M_1-\lambda_1^{-1}\Omega)\zeta_\ell&=\Bigl(1-\frac{\lambda_2}{\lambda_\ell}
\Bigr)x_3\zeta_\ell+(\lambda_\ell^{-1}-\lambda_1^{-1})\Omega\zeta_\ell,\quad
\ell\ge 3.
\endaligned
\label{5.1.27}
\eeq
Now \eqref{5.1.15} plus \eqref{5.1.27} yields
\beq
\text{Range }(M_1-\lambda_1^{-1}\Omega)\subset V_{12}=\Span \{\zeta_\ell:
\ell\ge 2\}.
\label{5.1.28}
\eeq
Hence any eigenfunction of $M_1$ must lie in $V_{12}$.  Also, \eqref{5.1.28}
implies $M_1:V_{12}\rightarrow V_{12}$, so
\beq
M_1-\lambda_2^{-1}\Omega:V_{12}\longrightarrow V_{12}.
\label{5.1.29}
\eeq
On the other hand (parallel to \eqref{4.3.32}--\eqref{4.3.33}, noting that
\eqref{4.3.31} holds with $\lambda_\nu=\lambda_2$),
\beq
(M_1-\lambda_2^{-1}\Omega)\zeta=\lambda_2^{-1}B(x_3)(I+\lambda_2\Delta^{-1})
\zeta,
\label{5.1.30}
\eeq
for $\zeta\in V_{12}$, and since $I+\lambda_2\Delta^{-1}$ is positive semidefinite
on $V_{12}$, it follows that
\beq
M_1\bigr|_{V_{12}}\ \text{ has real spectrum.}
\label{5.1.31}
\eeq
This proves Proposition \ref{p5.1.2}.
\qed

\sk
We conjecture linear stability when $\Omega>0$:

\sk
{\bf Conjecture.}  In the current setting, $\{e^{itM_1}:t\in\RR\}$ is
uniformly bounded on $V_1$ for each $\Omega>0$.

\subsection{Matrix analysis for $f(x)=cP_3(x_3)$}\label{c5s2}

We work with the following modification of the setting of \S{\ref{c5s1}}.
As there,
\beq
V_k=\{f\in L^2_b(S^2):X_3f=ikf\},
\label{5.2.1}
\eeq
and we set $M_k=M|_{V_k}$, with
\beq
M\zeta=A(x_3)\zeta+B(x_3)\Delta^{-1}\zeta.
\label{5.2.2}
\eeq
As before,
\beq
A(x_3)=f'(x_3),\quad B(x_3)=\Omega-w'(x_3),
\label{5.2.3}
\eeq
where $w=\Delta f$.  As before, we take $f(x_3)$ to be a zonal eigenfunction
of $\Delta$, hence a multiple of $P_\nu(x_3)$, for some $\nu$.  We saw in
\S{\ref{c4s3}} that taking $\nu=0$ or $\nu=1$ does not work,
and in \S{\ref{c5s1}} that taking $\nu=2$ does not work,
to produce examples of non-real eigenvalues of $M_1$.
Here, we take $\nu=3$.  Now
\beq
f(x_3)=\alpha P_3(x_3)\Longrightarrow \Delta f=-\lambda_3f,\quad \lambda_3=12.
\label{5.2.4}
\eeq
General formulas yield
\beq
P_3(t)=\frac{5}{2}t^3-\frac{3}{2}t,\quad \text{hence }\ P'_3(t)=
\frac{15}{2}\Bigl(t^2-\frac{1}{5}\Bigr).
\label{5.2.5}
\eeq
Hence, in \eqref{5.2.2}, we will take
\beq
A(x_3)=x_3^2-\frac{1}{5},\quad B(x_3)=\Omega+\lambda_3\Bigl(x_3^2-\frac{1}{5}
\Bigr).
\label{5.2.6}
\eeq
In light of Proposition \ref{p4.3.7}, we are interested in the behavior of
$M_1$ and $M_2$.

We start with an analysis of $M_1$.
We take the orthonormal basis $\{\zeta_\ell:\ell\ge 1\}$ of $V_1$ given by
\eqref{5.1.6}--\eqref{5.1.10}.  As noted there
\beq
\Delta\zeta_\ell=-\lambda_\ell\zeta_\ell,\quad \lambda_\ell=\ell(\ell+1),
\label{5.2.7}
\eeq
and (cf.~(5.1.15))
\beq
x_3\zeta_\ell=a_\ell\zeta_{\ell-1}+a_{\ell+1}\zeta_{\ell+1},
\label{5.2.8}
\eeq
with
\beq
a_\ell=\sqrt{\frac{(\ell+1)(\ell-1)}{(2\ell+1)(2\ell-1)}}.
\label{5.2.9}
\eeq
In the current setting,
\beq
M_1\zeta_\ell=\Bigl(x_3^2-\frac{1}{5}\Bigr)\zeta_\ell-\frac{1}{\lambda_\ell}
\Bigl(\Omega+\lambda_3\Bigl(x_3^2-\frac{1}{5}\Bigr)\Bigr)\zeta_\ell,
\label{5.2.10}
\eeq
so it is useful to note that \eqref{5.2.8} implies
\beq
x_3^2\zeta_\ell=a_\ell a_{\ell-1}\zeta_{\ell-2}
+(a_\ell^2+a_{\ell+1}^2)\zeta_\ell
+a_{\ell+1}a_{\ell+2}\zeta_{\ell+2},
\label{5.2.11}
\eeq
or equivalently
\beq
x_3^2\zeta_\ell=b_\ell\zeta_{\ell-2}+c_\ell\zeta_\ell+b_{\ell+2}\zeta_{\ell+2},
\label{5.2.12}
\eeq
with
\beq
b_\ell=a_\ell a_{\ell-1},\quad c_\ell=a_\ell^2+a_{\ell+1}^2.
\label{5.2.13}
\eeq
In \eqref{5.2.7}--\eqref{5.2.12}, $\ell\ge 1$.
We use the natural convention that $\zeta_0=
\zeta_{-1}=0$.  Putting together \eqref{5.2.10} and \eqref{5.2.11} yields
\beq
\aligned
M_1\zeta_\ell=\
&b_\ell\zeta_{\ell-2}+\Bigl(c_\ell-\frac{1}{5}\Bigr)\zeta_\ell
+b_{\ell+2}\zeta_{\ell+2} \\
&-\frac{\lambda_3}{\lambda_\ell} \Bigl(b_\ell\zeta_{\ell-2}
+\Bigl(c_\ell-\frac{1}{5}\Bigr)\zeta_\ell+b_{\ell+2}\zeta_{\ell+2}\Bigr) \\
&-\frac{\Omega}{\lambda_\ell} \zeta_\ell.
\endaligned
\label{5.2.14}
\eeq

Our goal is to investigate for what $\Omega\ge 0$ does
\beq
\text{$M_1$ have a non-real eigenvalue.}
\label{5.2.15}
\eeq
As we know, the ``Rayleigh criterion'' produces
\beq
0\le\Omega<\frac{4}{5}\lambda_3=\lambda_3\, \max\limits_{|x_3|\le 1}\,
\Bigl(x_3^2-\frac{1}{5}\Bigr),
\label{5.2.16}
\eeq
as a {\it necessary} condition for \eqref{5.2.15} to hold.
We make a numerical study of \eqref{5.2.14}
to indicate how close \eqref{5.2.16} is to being sufficient.

Numerical experiments, described in \S{\ref{c5s3}},
indicate that \eqref{5.2.15} holds for $0\le\Omega<\gamma$
with $\gamma\approx 1$.  This is a lot smaller than $(4/5)\lambda_3=9.6$.

We move along from $M_1$ to $M_2$, i.e., we take $k=2$ in \eqref{5.2.1}.
Parallel to \eqref{5.1.6}, an orthogonal basis of $V_2$ is given by
\beq
\tilde{\zeta}_\ell=e^{2i\psi}P^2_\ell(x_3),\quad \ell\ge 2.
\label{5.2.17}
\eeq
In this case, the associated Legendre function $P^2_\ell$ is given in
(7.12.5) of \cite{Leb} as
\beq
P^2_\ell(t)=(1-t^2)P''_\ell(t).
\label{5.2.18}
\eeq
We mention that
\beq
\Delta\tilde{\zeta}_\ell=-\lambda_\ell \tilde{\zeta}_\ell,\quad
\lambda_\ell=\ell(\ell+1).
\label{5.2.19}
\eeq
We emphasize that here $\ell\ge 2$.  In \eqref{5.1.6}--\eqref{5.1.8}, we
had $\ell\ge 1$.  As for the norms of these functions, one has
(\cite{Leb}, p.~201, \#10)
\beq
\int_{-1}^1 P^2_\ell(t)^2\, dt=\frac{2}{2\ell+1}\,\frac{(\ell+2)!}{(\ell-2)!}
=\frac{2(\ell+2)(\ell+1)\ell(\ell-1)}{2\ell+1}.
\label{5.2.20}
\eeq
Hence, up to an unimportant constant, an orthonormal basis of $V_2$ is given by
\beq
\zeta_\ell=\sqrt{\frac{2\ell+1}{2(\ell+2)(\ell+1)\ell(\ell-1)}}\, e^{2i\psi}
(1-x_3^2)P''_\ell(x_3),\quad \ell\ge 2.
\label{5.2.21}
\eeq
The operator $M_2$ acts on this basis as
\beq
\aligned
M_2\zeta_\ell&=A(x_3)\zeta_\ell+B(x_3)\Delta^{-1}\zeta_\ell \\
&=A(x_3)\zeta_\ell-\frac{1}{\lambda_\ell}B(x_3)\zeta_\ell \\
&=\Bigl(x_3^2-\frac{1}{5}\Bigr)\zeta_\ell-\frac{1}{\lambda_\ell}
\Bigl(\Omega+\lambda_3\Bigl(x_3^2-\frac{1}{5}\Bigr)\Bigr)\zeta_\ell.
\endaligned
\label{5.2.22}
\eeq

To proceed, we need to write $x_3\zeta_\ell$ as a linear combination of
$\{\zeta_j\}$.  From \eqref{5.1.12} and \eqref{5.1.13}, we get,
upon applying $d/dt$,
\beq
\aligned
tP''_\ell(t)+P'_\ell(t)&=\frac{\ell+1}{2\ell+1}P''_{\ell-1}(t)
+\frac{\ell}{2\ell+1}P''_{\ell+1}(t), \\
P'_\ell(t)&=-\frac{1}{2\ell+1}P''_{\ell-1}(t)+\frac{1}{2\ell+1}P''_{\ell+1}(t),
\endaligned
\label{5.2.23}
\eeq
hence
\beq
tP''_\ell(t)=\frac{\ell+2}{2\ell+1}P''_{\ell-1}(t)+\frac{\ell-1}{2\ell+1}
P''_{\ell+1}(t).
\label{5.2.24}
\eeq
Plugging this into \eqref{5.2.21} gives
\beq
x_3\zeta_\ell=\sqrt{\frac{(\ell+2)(\ell-2)}{(2\ell+1)(2\ell-1)}}\zeta_{\ell-1}
+\sqrt{\frac{(\ell+3)(\ell-1)}{(2\ell+3)(2\ell+1)}}\zeta_{\ell+1},\quad \ell\ge 2.
\label{5.2.25}
\eeq
We use the natural convention that $\zeta_j=0$ for $j<2$.  It is convenient to
write \eqref{5.2.25} as
\beq
x_3\zeta_\ell=a_\ell\zeta_{\ell-1}+a_{\ell+1}\zeta_{\ell+1},\quad
a_\ell=\sqrt{\frac{(\ell+2)(\ell-2)}{(2\ell+1)(2\ell-1)}}.
\label{5.2.26}
\eeq
This in turn implies
\beq
x_3^2\zeta_\ell=a_\ell a_{\ell-1}\zeta_{\ell-2}+(a_\ell^2+a_{\ell+1}^2)\zeta_\ell
+a_{\ell+1}a_{\ell+2}\zeta_{\ell+2},
\label{5.2.27}
\eeq
or equivalently
\beq
x_3^2\zeta_\ell=b_\ell\zeta_{\ell-2}+c_\ell\zeta_\ell+b_{\ell+2}\zeta_{\ell+2},
\label{5.2.28}
\eeq
with
\beq
b_\ell=a_\ell a_{\ell-1},\quad c_\ell=a_\ell^2+a_{\ell+1}^2.
\label{5.2.29}
\eeq
Recall that in \eqref{5.2.26}--\eqref{5.2.29},
$\ell\ge 2$, and we use the convention that $\zeta_j=0$ for $j<2$.
Putting together \eqref{5.2.22} with \eqref{5.2.28} yields
\beq
\aligned
M_2\zeta_\ell=\
&b_\ell\zeta_{\ell-2}+\Bigl(c_\ell-\frac{1}{5}\Bigr)\zeta_\ell
+b_{\ell+2}\zeta_{\ell+2} \\
&-\frac{\lambda_3}{\lambda_\ell}\Bigl(b_\ell\zeta_{\ell-2}
+\Bigl(c_\ell-\frac{1}{5}\Bigr)\zeta_\ell+b_{\ell+2}\zeta_{\ell+2}\Bigr) \\
&-\frac{\Omega}{\lambda_\ell}\zeta_\ell,
\endaligned
\label{5.2.30}
\eeq
for $\ell\ge 2$.

Our goal is to investigate for what $\Omega\ge 0$ does
\beq
\text{$M_2$ have a non-real eigenvalue.}
\label{5.2.31}
\eeq
As we know, the ``Rayleigh criterion'' produces
\beq
0\le\Omega<\frac{4}{5}\lambda_3=\lambda_3\, \max\limits_{|x_3|\le 1}\,
\Bigl(x_3^2-\frac{1}{5}\Bigr),
\label{5.2.32}
\eeq
as a necessary condition for \eqref{5.2.31} to hold.  We make a numerical
study of \eqref{5.2.30} to indicate how close \eqref{5.2.32} is to being
sufficient.

Numerical experiments, described in \S{\ref{c5s3}},
indicate that \eqref{5.2.31} holds for $0\le\Omega<\gamma$ with
$\gamma\approx 1/2$, and that \eqref{5.2.31} ceases to hold for
$\Omega>\gamma$.  Again, $1/2$ is a lot smaller than $(4/5)\lambda_3$.

In \S{\ref{c5s3}} we will also examine such matrices that arise when
$f(x_3)=cP_4(x_3)$.  More precisely, we take
\beq
A(x_3)=x_3^3-\frac{3}{7}x_3,\quad B(x_3)=\Omega+\lambda_4\Bigl(x_3^3-
\frac{3}{7}x_3\Bigr),
\label{5.2.33}
\eeq
and define
\beq
M_k\zeta_\ell=\Bigl(x_3^3-\frac{3}{7}x_3\Bigr)\zeta_\ell
-\frac{1}{\lambda_\ell}
\Bigl(\Omega+\lambda_4\Bigl(x_3^3-\frac{3}{7}x_3\Bigr)\Bigr)\zeta_\ell,
\label{5.2.34}
\eeq
on the orthonormal basis $\{\zeta_\ell\}$ of $V_k$, given by
\eqref{5.1.6}--\eqref{5.1.10} for $k=1$, by \eqref{5.2.21} for $k=2$ and by a
comparable strategy for $k=3$.
We make use of matrix formulas parallel to \eqref{5.2.14} and \eqref{5.2.30}
in this setting, which are given in \S{\ref{c5s3}}.

\subsection{Numerical study of truncated matrices}\label{c5s3}

Here we study truncated versions of the matrix operators arising from
the attack described in \S{\ref{c5s2}} on
stability of banded structures for the Euler equations with Coriolis forces
on the sphere.  Before describing how this is done, we mention one result
that leads one to believe that eigenvectors of such operators $M_k$ as
described in \eqref{5.2}, or more generally \eqref{4.3.2}--\eqref{4.3.3},
should be expected to be captured fairly accurately by such a truncation.

\begin{proposition} \label{p5.3.1}
Given $A$ and $B$ real valued and real analytic on $S^2$ in \eqref{4.3.3}, if
$\zeta\in V_k$ is an eigenvector of $M_k$, with eigenvalue $\lambda\notin\RR$,
then $\zeta$ is real analytic on $S^2$.
\end{proposition}

\demo
As seen in \eqref{4.3.6}, $\eta=\Delta^{-1}\zeta$ satisfies an elliptic
partial differential equation with analytic coefficients, so real analyticity
of $\eta$, hence of $\zeta=\Delta\eta$, follows.
\qed

Given that such an eigenfunction $\zeta\in V_k$ of $M_k$ is real analytic,
its spherical harmonic expansion is rapidly (in fact, exponentially)
convergent.  The truncation of such $\zeta$ should then  be a high order
quasimode of the associated truncation of $M_k$.  We are not currently in
a position to derive rigorous conclusions about the spectrum of $M_k$
from numerical results on truncations, but we are motivated to take such
numerical results as a strong indication of how the spectrum of $M_k$
behaves.

Numerically, we implement  in {\it Matlab} the eigenvalue
solver {\it eigs} for an $N \times N$ block truncation of matrices $M_k$,
such as arise for $k=1$ in \eqref{5.2.10} and for $k=2$ in \eqref{5.2.22}.
We will refer to these finite matrices as $M_k^N$.

\begin{figure}
\includegraphics[scale=0.35]{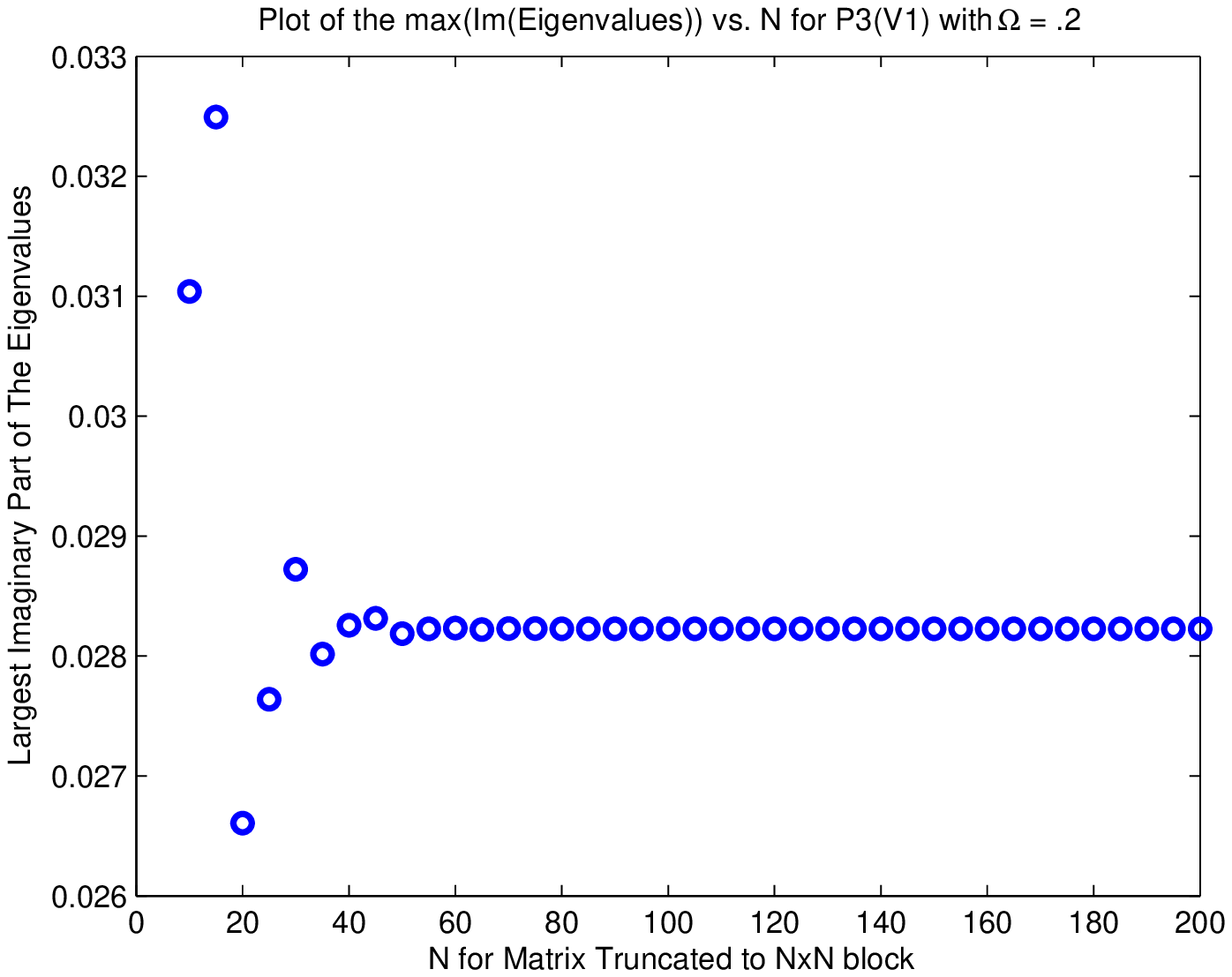}
\includegraphics[scale=0.35]{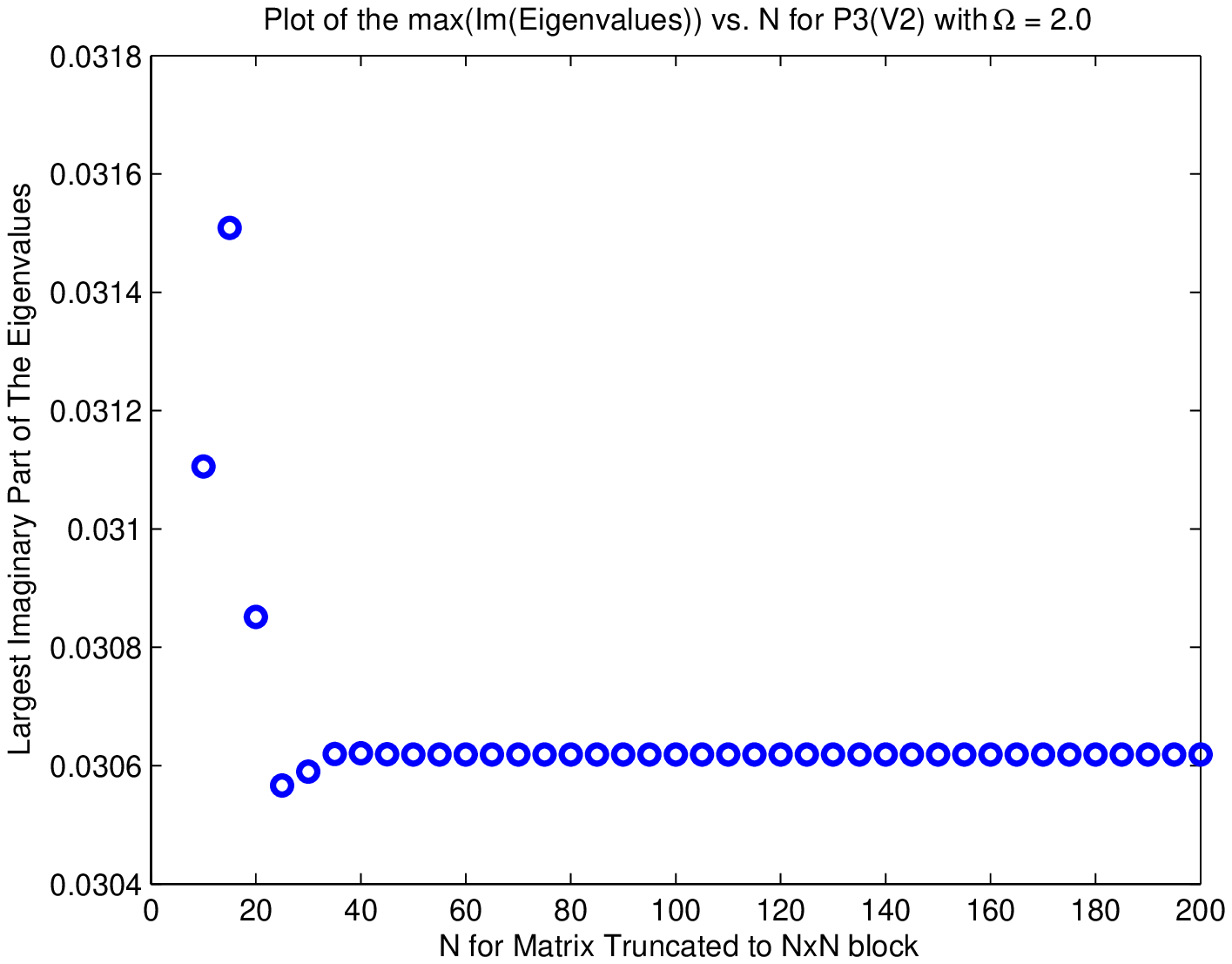}
\caption{{\bf Left}:  The size of the imaginary part of an
unstable eigenvalue, as a function of
$N$, for an $N \times N$ matrix truncation for $P_3 (V_1)$ with
$\Omega = .2$.  This is seen to stabilize strongly as soon as $N$
is sufficiently large, say $N > 100$. {\bf Right:}  The size of
such an imaginary part, as a function of $N$, for an $N \times N$ matrix
truncation for $P_3 (V_2)$ with $\Omega = 2$.  This is also seen to stabilize
strongly as soon as $N$ is sufficiently large, say
$N > 100$. }
\label{p3v1v2_Nconv}
\end{figure}

\begin{figure}
\includegraphics[scale=0.35]{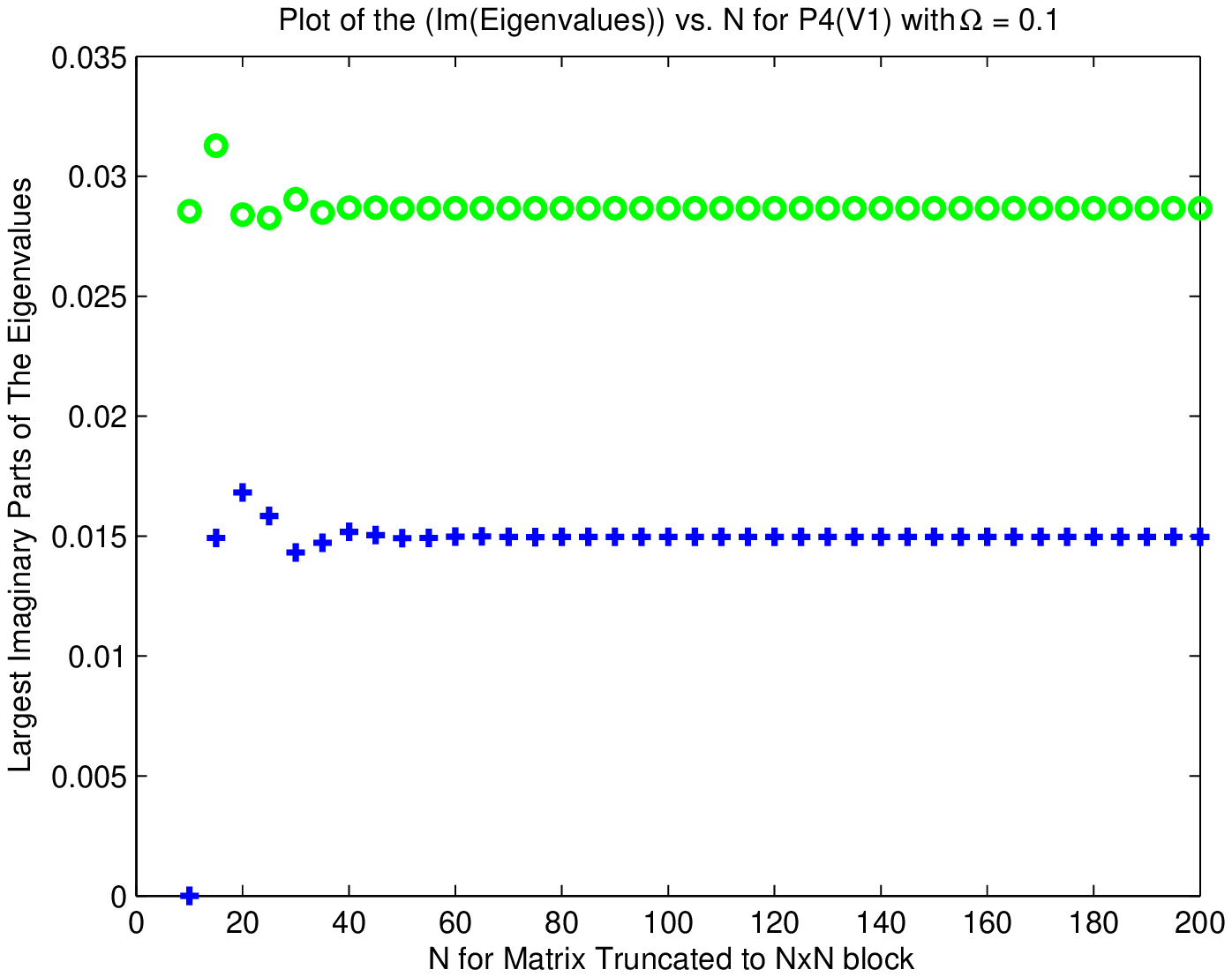}
\includegraphics[scale=0.35]{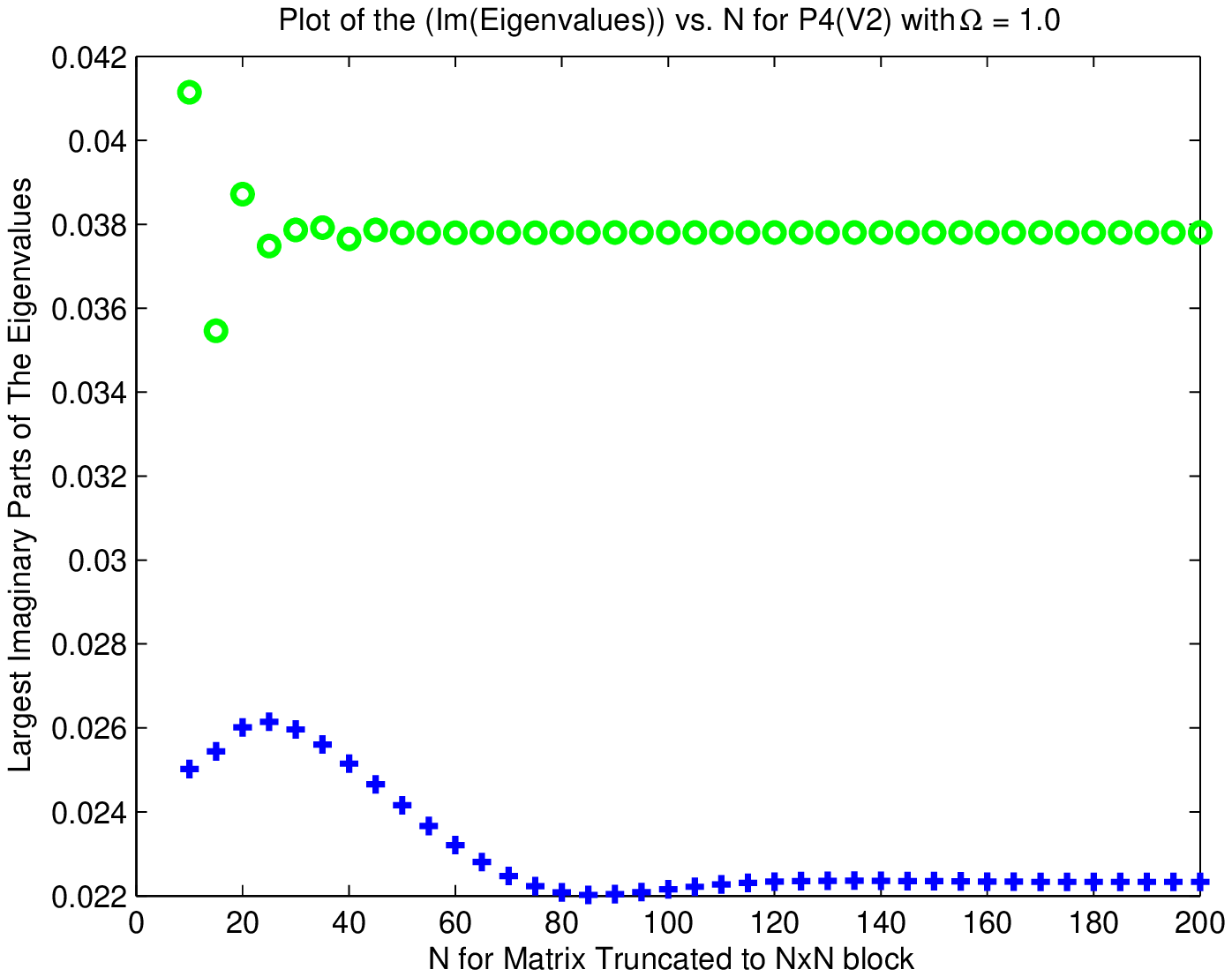}
\caption{{\bf Left}:  The size of the imaginary parts of unstable eigenvalues
as a function of $N$, for an $N \times N$ matrix truncation for
$P_4 (V_1)$ with $\Omega = .1$.  This example has multiple non-real
eigenvalues, but both are seen to stabilize strongly as soon as $N$
is sufficiently large, say $N > 100$.
{\bf Right:}  The size of the imaginary parts of unstable eigenvalues, as a
function of $N$, for an $N \times N$ matrix truncation
for $P_4 (V_2)$ with $\Omega = 1$.
This is also seen to stabilize strongly as soon as $N$ is
sufficiently large, say $N > 150$. }
\label{p4v1v2_Nconv}
\end{figure}

Our first observation is that the non-real spectrum of $M_k^N$ appears
to stabilize before $N$ becomes particularly large, and to be set in the
upper left block of such a matrix.  Figure \ref{p3v1v2_Nconv} illustrates
this for $N\times N$ matrix truncations for the $P_3(V_1)$ model, with
$\Omega=1/5$, and for the $P_3(V_2)$ model, with $\Omega=2$.  Figure
\ref{p4v1v2_Nconv} has analogous illustrations for the $P_4(V_1)$ model,
with $\Omega=1/10$, and the $P_4(V_2)$ model, with $\Omega=1$.
In three of these four cases, one sees stabilization of the non real
spectra in truncations well before $N$ reaches $100$, and in the fourth
case well before $N$ reaches 150.

Subsequent figures show the imaginary parts of the largest unstable
eigenvalues  of $M_k^N$, for various $P_\nu(V_k)$ models (cf.~\eqref{5.2})
as a function of $\Omega$.  For all these truncations, we took $N=400$.

%$$\text{}$$
%\begin{enumerate}
%\item Figure looking at the non-real part of spectrum of $M_{1,2}$ for
%increasing values of $N$ as numerical evidence of convergence.  These are
%currently plotted for $P_3$ and $P_4$ in Figures \ref{p3v1v2_Nconv},
%\ref{p4v1v2_Nconv}.

%\item  Plot of $\sigma(M_1^N)$ and $\sigma(M_2^N)$ vs. $\Omega$ exploring
%stability threshold and motion of eigenvalues.  These currently appear as
%Figures \ref{p3v1v2} and \ref{p4v1v2}.  {\bf  Compute $M_3^N$ for this case
%to complete the picture?!?}
%\end{enumerate}

For comparison to our numerically constructed matrices, it is 
convenient to write out small block representations of our desired matrices.
Given the $P_3 (V_1)$ model, we have as the $6\times 6$ truncation
\begin{align*}
M_1^6 & = \left(  \begin{array}{cccccc}
c_1 - \frac15 & 0 & b_3 & 0 & 0 & 0 \\
0 & c_2 - \frac15  & 0 & b_4 & 0 & 0 \\
b_3 & 0 & c_3 - \frac15  & 0 & b_5 & 0 \\
0 & b_4 & 0 & c_4  - \frac15  & 0 & b_6 \\
0 & 0 & b_5 & 0 & c_5 - \frac15  & 0 \\
0 & 0 & 0 & b_6 & 0 & c_6  - \frac15 
\end{array}
\right)
\\
& -  \left( \Omega -  \frac{\lambda_3}5 \right) \left(  \begin{array}{cccccc}
\lambda_1^{-1} & 0 & 0 & 0 & 0 & 0 \\
0 & \lambda_2^{-1} & 0 & 0 & 0 & 0 \\
0 & 0 & \lambda_3^{-1} & 0 & 0 & 0 \\
0 & 0 & 0 & \lambda_4^{-1} & 0 & 0 \\
0 & 0 & 0 & 0 & \lambda_5^{-1}& 0 \\
0 & 0 & 0 & 0 & 0 & \lambda_6^{-1} 
\end{array}
\right)  \\
 & \hspace{.5cm} - \lambda_3 \left(  \begin{array}{cccccc}
c_1/\lambda_1 & 0 & b_3/\lambda_3 & 0 & 0 & 0 \\
0 & c_2/\lambda_2& 0 & b_4/\lambda_4 & 0 & 0 \\
b_3/\lambda_1 & 0 & c_3/\lambda_3 & 0 & b_5/\lambda_5 & 0 \\
0 & b_4/\lambda_2 & 0 & c_4/\lambda_4 & 0 & b_6/\lambda_6 \\
0 & 0 & b_5/\lambda_3 & 0 & c_5/\lambda_5 & 0 \\
0 & 0 & 0 & b_6/\lambda_4 & 0 & c_6/\lambda_6
\end{array}
\right)
\end{align*}
for
\[
b_\ell = \sqrt{ \frac{ (\ell +1 ) (\ell - 1) }{ (2 \ell +1 )(2 \ell -1) } }
\sqrt{ \frac{ (\ell ) (\ell - 2) }{ (2 \ell -1 )(2 \ell -3) } } ,
\]
\[
c_\ell =  \frac{ (\ell +1 ) (\ell - 1) }{ (2 \ell +1 )(2 \ell - 1) }
+  \frac{ \ell (\ell +2) }{ (2 \ell +1 )(2 \ell + 3) },
\]
and $\lambda_\ell=\ell(\ell+1)$.  Compare \eqref{5.2.14}.  There is a comparable construction for $P_3 (V_2)$ using as in \eqref{5.2.26} and \eqref{5.2.30}.

\begin{figure}
\includegraphics[scale=0.35]{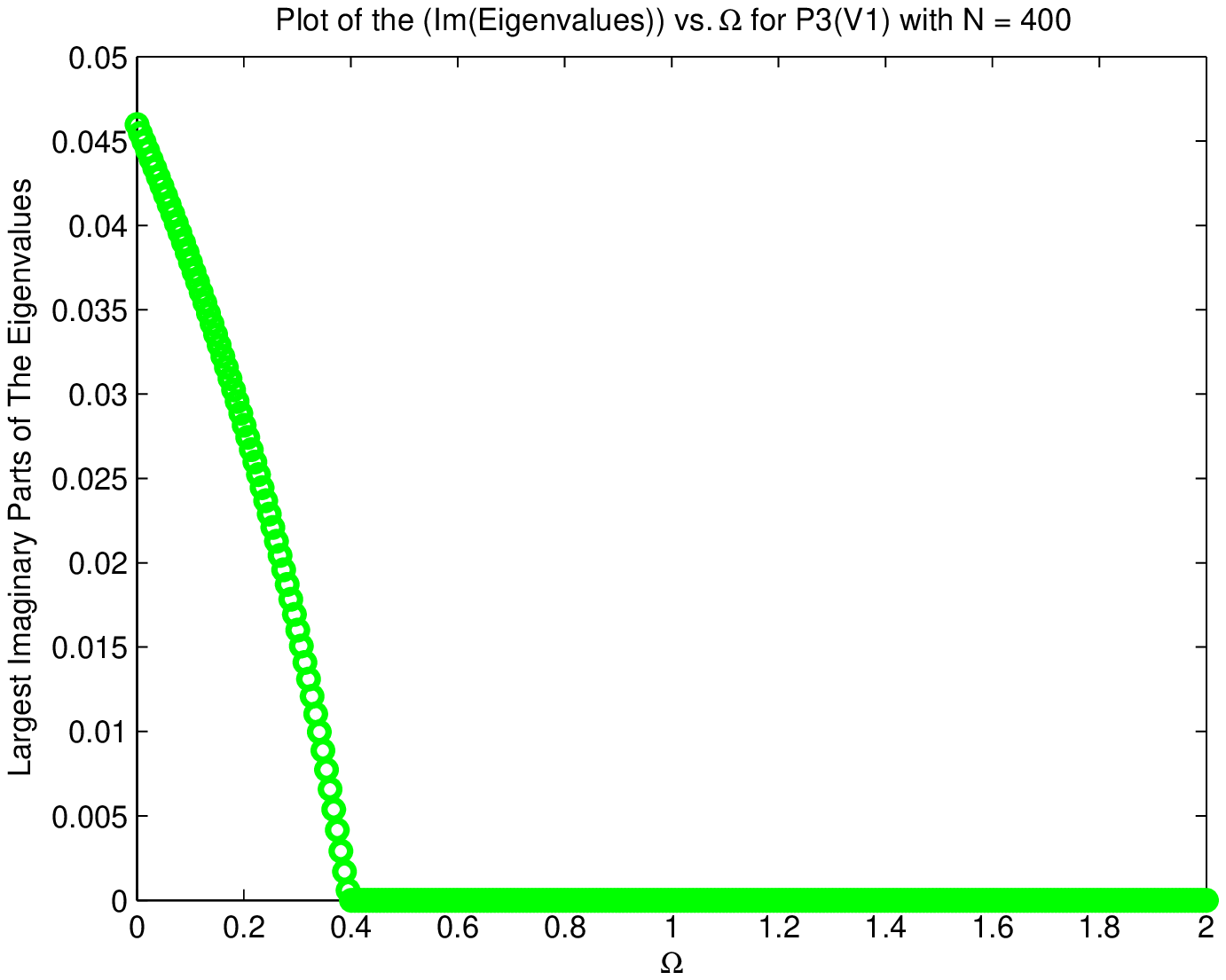}
\includegraphics[scale=0.35]{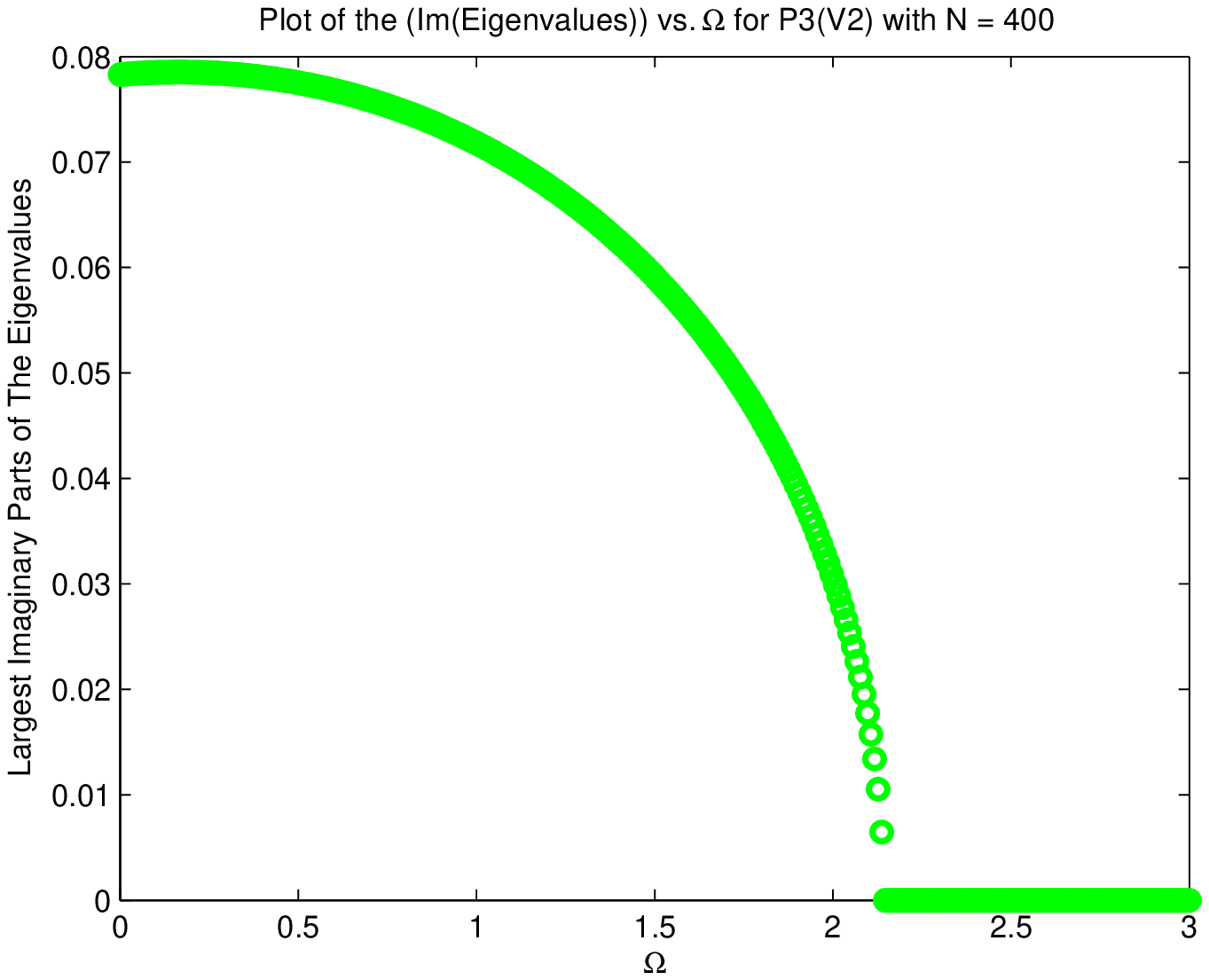}
\caption{{\bf Left:} The size of the imaginary part for the the unstable
eigenvalue for $P_3 (V_1)$, as a function of $\Omega$, with $N = 400$.
Here, we see stability for $\Omega > .4$, well below the Arnold stability
bound.  {\bf Right:} The size of the imaginary part for the the unstable
eigenvalue for $P_3 (V_2)$ as a function of $\Omega$ with $N=400$.  We see
stability for $\Omega > 2.25$, which is again well below the Arnold stability
bound.}
\label{p3v1v2}
\end{figure}

For the $P_4(V_1)$ model, we take  $A(x_3)=x_3(x_3^2-3/7)$, as in
\eqref{5.2.33}.  As a result, for the operation of $M_1$ on $\zeta_\ell$,
we have
\beq
\aligned
\label{Mp2v1}
A(x_3) \zeta_\ell =\ & \Bigl(x_3^3 - \frac37 x_3\Bigr) \zeta_\ell \\
=\ & (a_\ell a_{\ell - 1} a_{\ell-2}) \zeta_{\ell - 3}
+ a_\ell \Bigl(a_{\ell-1}^2 + a_\ell^2 + a_{\ell+1}^2 - \frac37\Bigr)
\zeta_{\ell-1} \\
& +  a_{\ell+1} \Bigl(a_{\ell}^2 + a_{\ell+1}^2 + a_{\ell+2}^2 - \frac37\Bigr)
\zeta_{\ell+1}
+  (a_{\ell +1} a_{\ell +2} a_{\ell + 3}) \zeta_{\ell + 3},
\endaligned
\eeq
with $a_\ell$ as in \eqref{5.2.9}, and 
\begin{equation}
B(x_3)\Delta^{-1}\zeta_\ell = -\frac{1}{\lambda_\ell}
( \Omega + \lambda_4 A(x_3)) \zeta_{\ell}.
\label{5.3.2}
\end{equation}

Again, for the reader's convenience and to assist with comparison to
numerically constructed matrices, we write down the $6\times 6$ truncation
for the $P_4 (V_1)$ model:
\begin{align*}
M_1^6 & = \left(  \begin{array}{cccccc}
0& b_2 & 0 & c_4 & 0 & 0 \\
b_2 & 0  & b_3 & 0 & c_5 & 0 \\
0 & b_3 & 0 &  b_4& 0 & c_6 \\
c_4 & 0 & b_4 & 0 & b_5 & 0 \\
0 & c_5 & 0& b_5 & 0 & b_6 \\
0 & 0 & c_6 & 0 & b_6 & 0
\end{array}
\right)
-\Omega\left(  \begin{array}{cccccc}
\lambda_1^{-1} & 0 & 0 & 0 & 0 & 0 \\
0 & \lambda_2^{-1} & 0 & 0 & 0 & 0 \\
0 & 0 & \lambda_3^{-1} & 0 & 0 & 0 \\
0 & 0 & 0 & \lambda_4^{-1} & 0 & 0 \\
0 & 0 & 0 & 0 & \lambda_5^{-1}& 0 \\
0 & 0 & 0 & 0 & 0 & \lambda_6^{-1} 
\end{array}
\right)  \\
 & \hspace{.5cm} - \lambda_4 \left(  \begin{array}{cccccc}
0& b_2/\lambda_2 & 0 & c_4/\lambda_4 & 0 & 0 \\
b_2/\lambda_1 & 0  & b_3/\lambda_3 & 0 & c_5/\lambda_5 & 0 \\
0 & b_3/\lambda_2 & 0 &  b_4/\lambda_4& 0 & c_6/\lambda_6 \\
c_4/\lambda_1 & 0 & b_4/\lambda_3 & 0 & b_5/\lambda_5 & 0 \\
0 & c_5/\lambda_2 & 0& b_5/\lambda_4 & 0 & b_6/\lambda_6\\
0 & 0 & c_6/\lambda_3 & 0 & b_6/\lambda_5 & 0
\end{array}
\right)\end{align*}
for
$$
\aligned
b_\ell =\ &\sqrt{ \frac{ (\ell +1 ) (\ell - 1) }{ (2 \ell +1 )(2 \ell -1)}} \\
&\times \left(  \frac{ (\ell ) (\ell - 2) }{ (2 \ell -1 )(2 \ell -3) }
+\frac{ (\ell +1 ) (\ell - 1) }{ (2 \ell +1 )(2 \ell -1) }
+\frac{ (\ell +2 ) (\ell ) }{ (2 \ell +3 )(2 \ell +1) }  -\frac37 \right),
\endaligned
$$
\[
c_\ell = \sqrt{ \frac{ (\ell +1 ) (\ell - 1) }{ (2 \ell +1 )(2 \ell -1) } }  \sqrt{ \frac{ (\ell ) (\ell - 2) }{ (2 \ell -1 )(2 \ell -3) } }  \sqrt{ \frac{ (\ell -1 ) (\ell - 3) }{ (2 \ell -3)(2 \ell -5) } } ,
\]
and, again, $\lambda_\ell=\ell(\ell+1)$.

\begin{figure}
\includegraphics[scale=0.35]{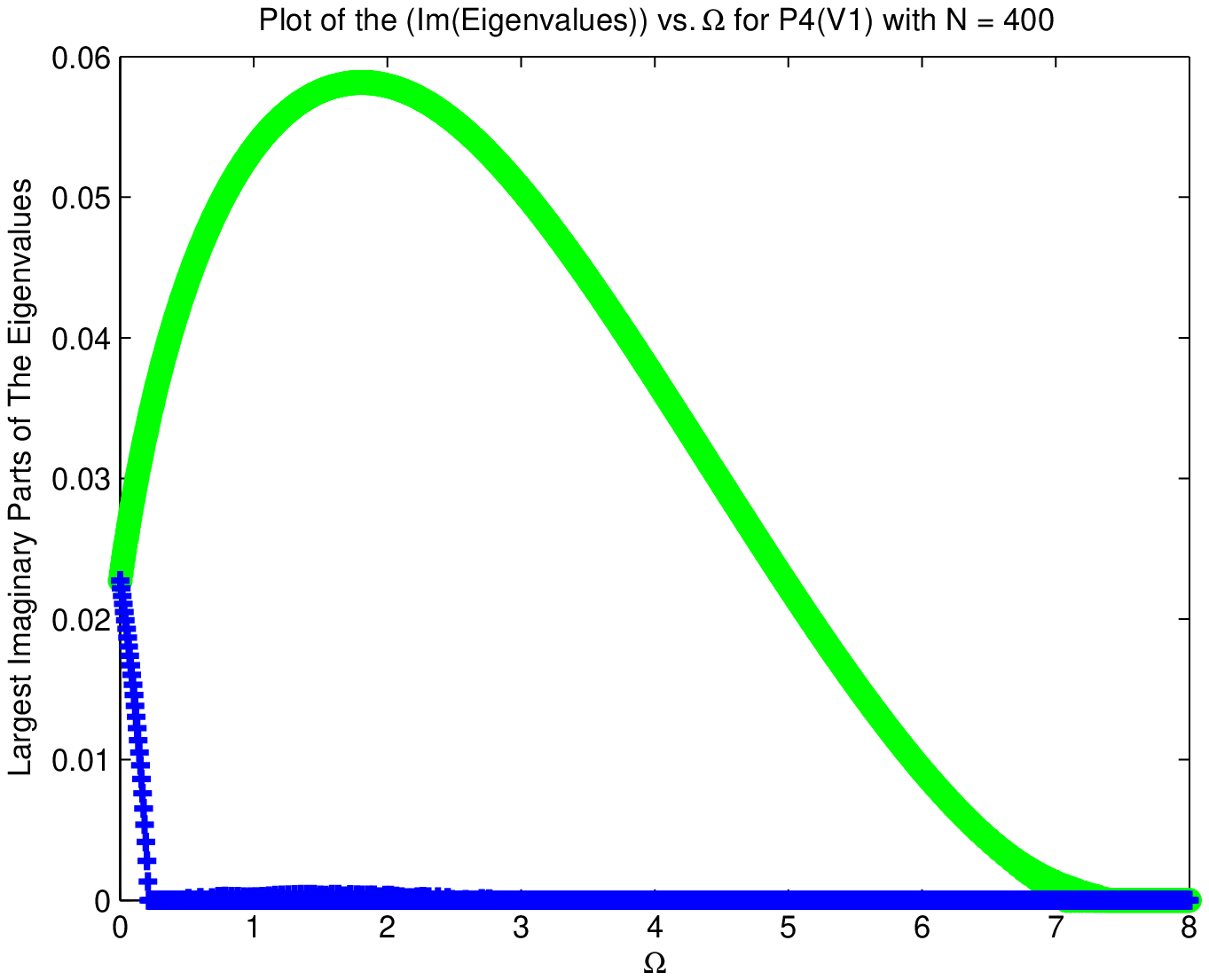}
\includegraphics[scale=0.35]{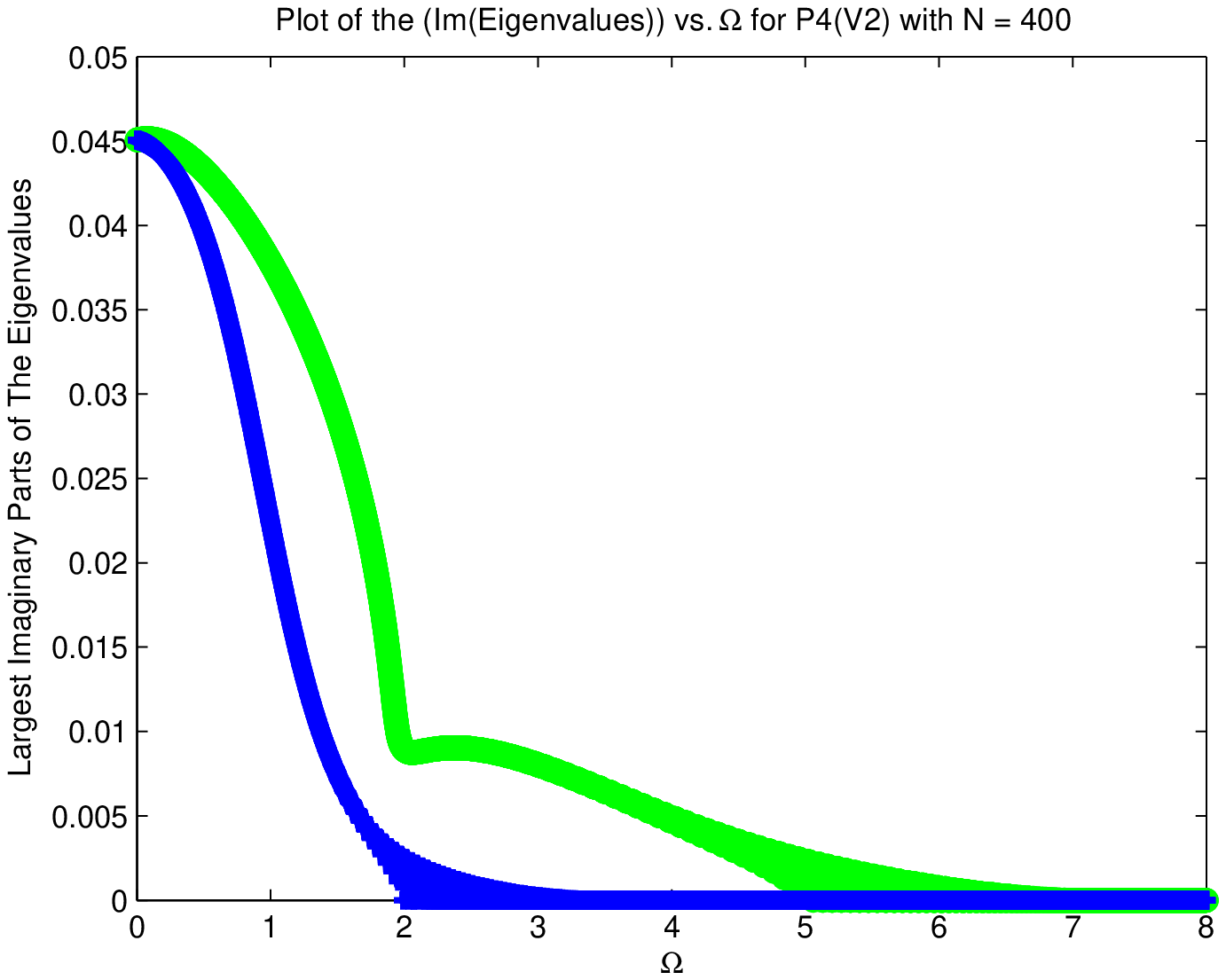} \\
\includegraphics[scale=0.35]{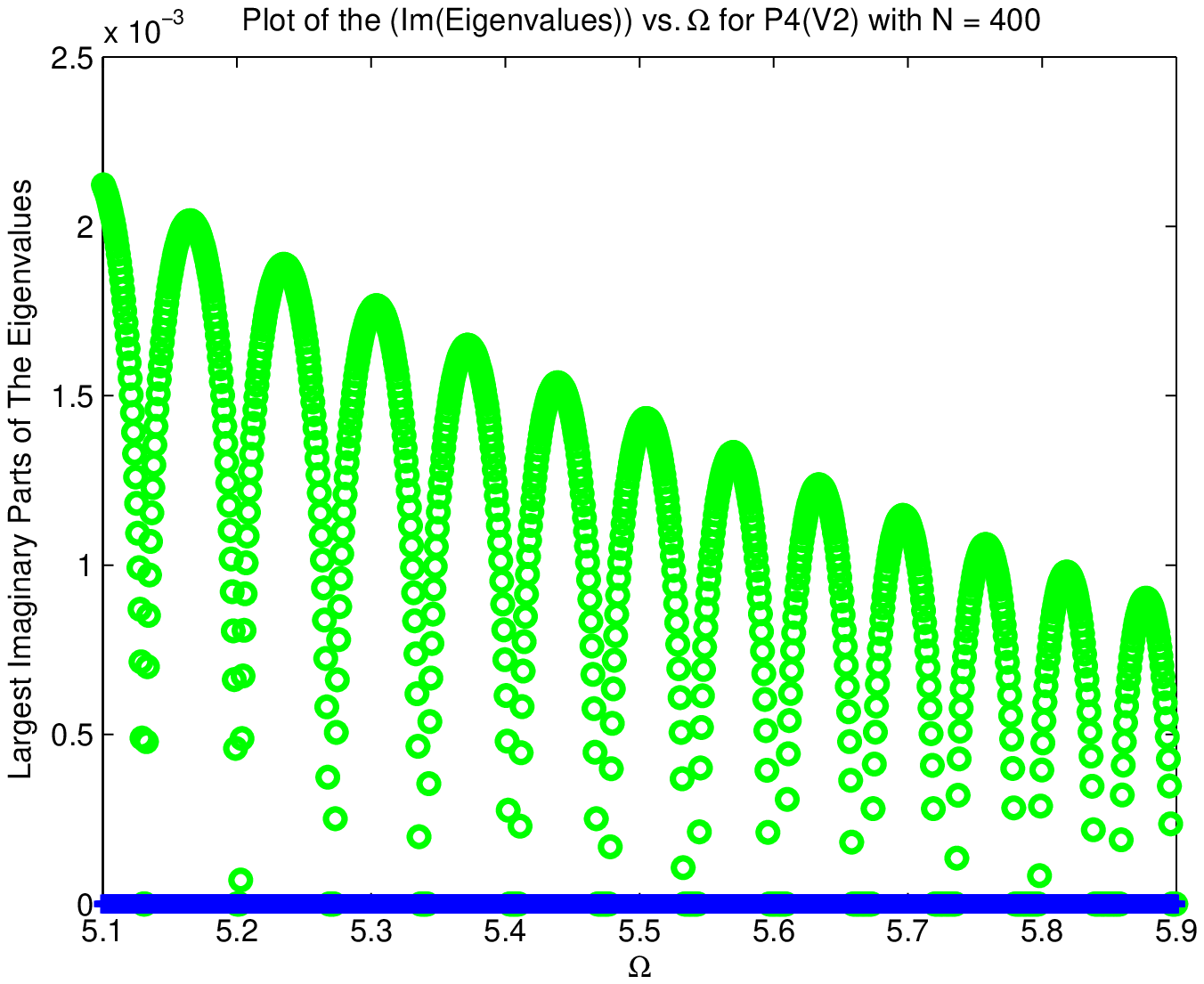}
\includegraphics[scale=0.35]{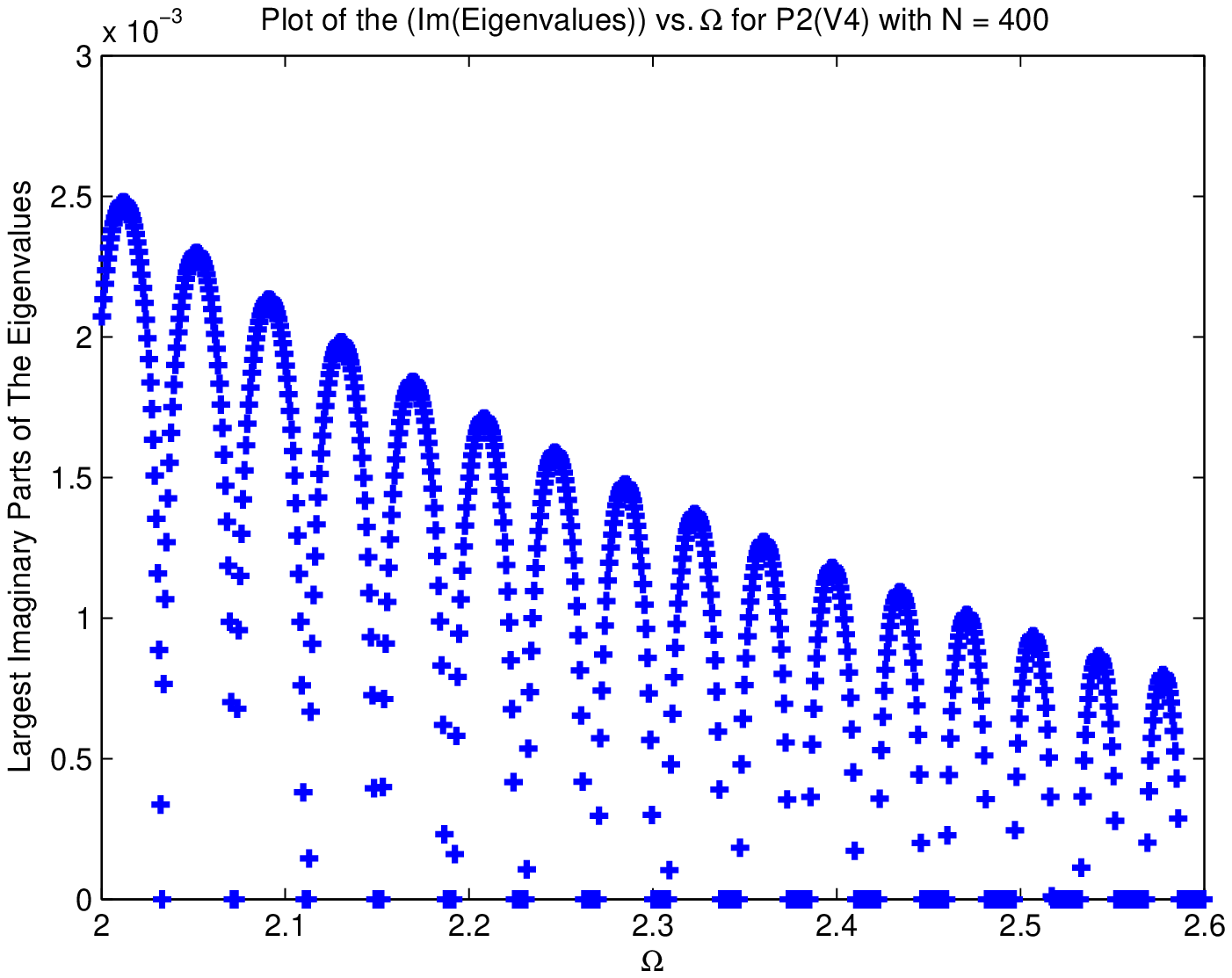} \\
\caption{{\bf Top Left:}  The size of the imaginary parts of the
two largest unstable eigenvalues for
$P_4 (V_1)$ as a function of $\Omega$ with $N = 400$.  Here, we see stability for
$\Omega > 7.5$, well below the Arnold stability bound.
{\bf Top Right:} The size of the imaginary parts of the two
largest unstable eigenvalues for
$P_4 (V_2)$ as a function of $\Omega$.   Here, we see stability for
$\Omega > 7$, well below the Arnold stability bound.
{\bf Bottom Left:}  A blow-up of the tail for the
imaginary part of the largest unstable eigenvalue
$P_4 (V_2)$ to show that the oscillations towards stable are smooth at fine
scales and not numerical errors.  {\bf Bottom Right:} A blow-up of the tail
for the imaginary part of the
smallest unstable eigenvalue $P_4 (V_2)$ to show that the
oscillations towards stable are smooth at fine scales and not numerical
errors.}
\label{p4v1v2}
\end{figure}

\begin{figure}
\includegraphics[scale=0.35]{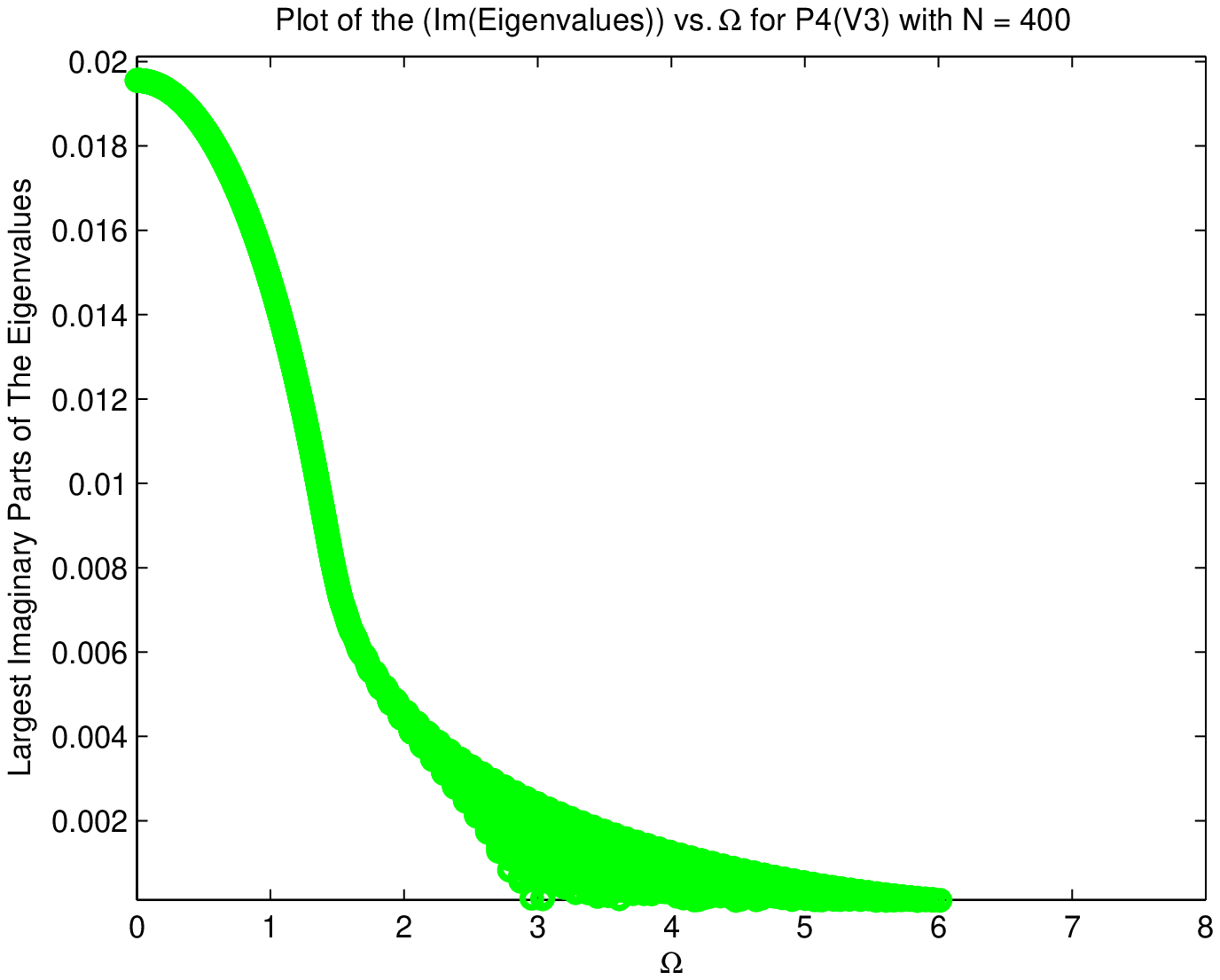}
\includegraphics[scale=0.35]{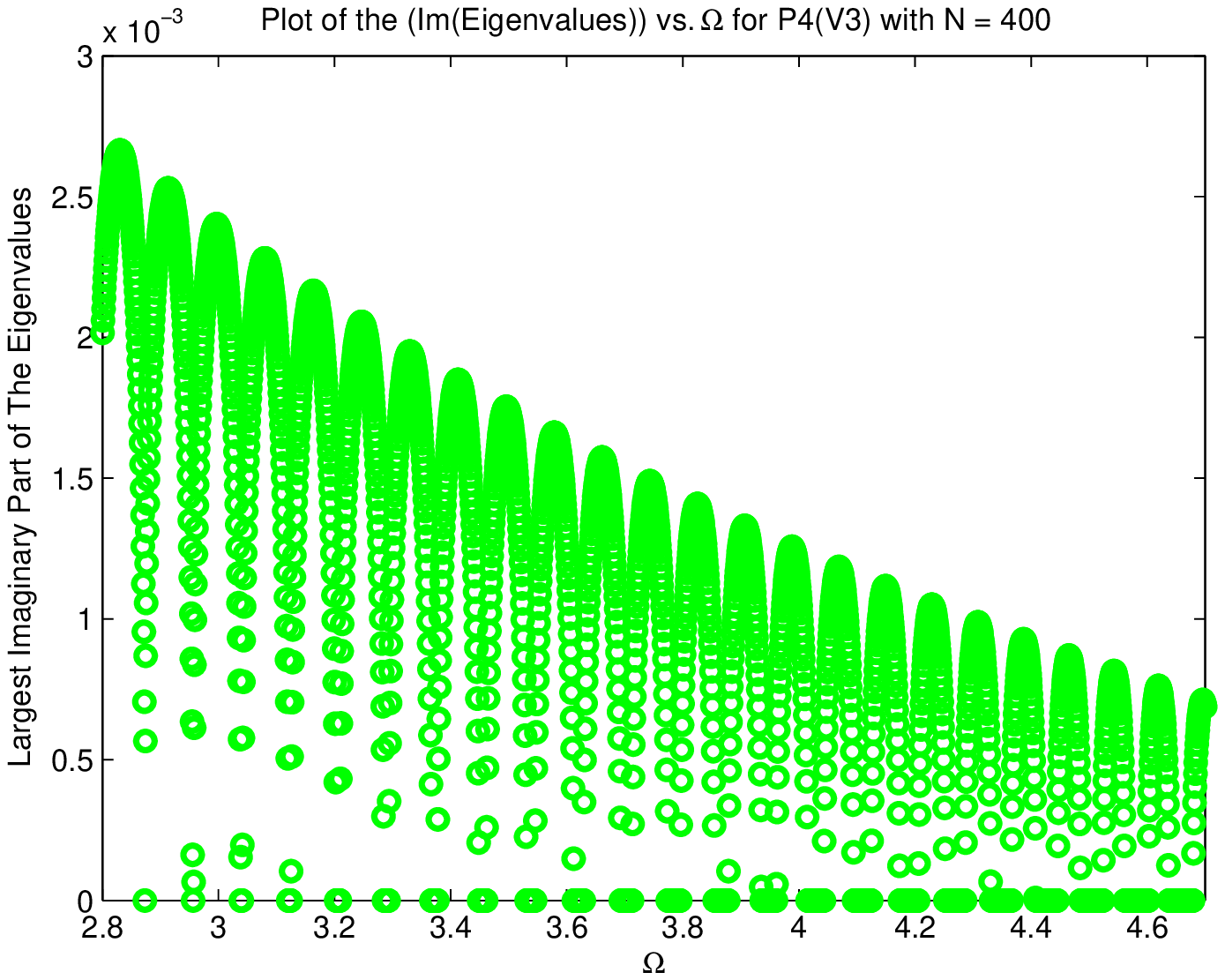}
\caption{
%{\bf Top Left:}  A blow-up of the tail for the
%imaginary part of the largest unstable eigenvalue
%$P_4 (V_2)$ to show that the oscillations towards stable are smooth at fine
%scales and not numerical errors.  {\bf Top Right:} A blow-up of the tail
%for the imaginary part of the
%smallest unstable eigenvalue $P_4 (V_2)$ to show that the
%oscillations towards stable are smooth at fine scales and not numerical
%errors.  
{\bf Left:}  The size of the imaginary part of the
 largest unstable eigenvalue for
$P_4 (V_3)$ as a function of $\Omega$.  Here, we see stability for
$\Omega > 6.5$, well below the Arnold stability bound.
{\bf Right:} The size of the imaginary part of the
largest unstable eigenvalues for
$P_4 (V_3)$ as a function of $\Omega$ blown up to see that it is smooth curve.}
\label{p4v2v3}
\end{figure}

For the $P_4(V_2)$ model, we also take $A(x_3)$ as in \eqref{5.2.33}.
Then, we observe  that for the operation of $M_2$ on $\zeta_\ell$, we have
for multiplication by $A(x_3)$ the same expression as in \eqref{Mp2v1} but with
\begin{equation}
a_\ell = \sqrt{\frac{(\ell+2)(\ell-2)}{(2\ell+1)(2\ell-1)}},
\end{equation}
as in \eqref{5.2.26} and $B(x_3)\Delta^{-1}\zeta_\ell$ as in \eqref{5.3.2}.

For the $P_4(V_3)$ model, we also take $A(x_3)$ as in \eqref{5.2.33}.
Then, we observe  that for the operation of $M_3$ on $\zeta_\ell$, we have
for multiplication by $A(x_3)$ the same expression as in \eqref{Mp2v1} but with
\begin{equation}
a_\ell = \sqrt{\frac{(\ell+3)(\ell-3)}{(2\ell+1)(2\ell-1)}},
\end{equation}
and, again, $B(x_3)\Delta^{-1}\zeta_\ell$ as in \eqref{5.3.2}.

We turn to a discussion of the spectral results recorded in Figures
\ref{p3v1v2}--\ref{p4v2v3}.  Figure \ref{p3v1v2} deals with the $P_3(V_k)$
models, for $k=1,2$.  In both cases, the graph indicates that $M^N_k$ has a
single eigenvalue with positive imaginary part, for $\Omega$ in a certain
range ($0\le\Omega\le 0.4$ for $k=1$, $0\le\Omega\le 2.25$ for $k=2$),
and this imaginary part decreases monotonically to $0$ as $\Omega$ runs over
these intervals.  In both cases, the imaginary part reaches $0$ for $\Omega$
well below the threshold specified by the Arnold method, as worked out in
\S{\ref{c4s1}}.  Figure \ref{p4v1v2} deals with $P_4(V_k)$ models, for $k=1,2$.
In this case, we have two eigenvalues of $M_k$ with positive imaginary part,
for $0\le\Omega\le 7.5$ and $0\le\Omega\le 0.3$ for $k=1$ and $0\le \Omega\le
3.4$ and $0\le\Omega\le 7.0$ for $k=2$.  Again, all imaginary parts reach $0$
for $\Omega$ well below the threshold to which the Arnold stability analysis
applies.  In these figures, we see a new phenomenon.  Namely, the imaginary
parts of the eigenvalues are no longer monotone functions of $\Omega$.
In fact, as seen in the bottom parts of Figure \ref{p4v1v2} (and also the top
part of Fig.~\ref{p4v2v3}) there is considerable oscillation of these
imaginary parts, as a function of $\Omega$, over certain ranges of $\Omega$,
particularly for $k=2$.  The bottom part of Figure \ref{p4v2v3} has
analogous graphs, for the $P_4(V_3)$ model, also illustrating such oscillation.

It is our expectation that similar results hold for the operators $M_k$, and
hence $\Gamma$, arising in the linearizaton procedure of \S{\ref{c4s2}}.
Going further, we imagine there are further stability and instability results
to be established for the Euler equation \eqref{1.1}.  We look forward to
future progress on these problems.

$$\text{}$$
Michael Taylor, corresponding author
\newline
Mathematics Dept., University of North Carolina
\newline
Chapel Hill NC 27599
\newline
Email: met @ math.unc.edu

$$\text{}$$
Jeremy Marzuola
\newline
Mathematics Dept., University of North Carolina
\newline
Chapel Hill NC 27599
\newline
Email: marzuola @ email.unc.edu

\end{document}